\documentclass[11pt]{amsart}
\usepackage{amsmath,amsthm,amssymb,xypic}

\theoremstyle{plain}  

\newtheorem{thm}{Theorem}[section]
\newtheorem{prop}[thm]{Proposition}
\newtheorem{lem}[thm]{Lemma}   
\newtheorem{cor}[thm]{Corollary}

\newtheorem{prop-defn}[thm]{Proposition--Definition}

\newtheorem*{S1a}{Condition (S1)(a)}
\newtheorem*{S1b}{Condition (S1)(b)}
\newtheorem*{S2}{Condition (S2)}
\newtheorem{cond}[thm]{Condition}

\theoremstyle{definition}

\newtheorem{defn}[thm]{Definition}
\newtheorem{notn}[thm]{Notation}

\theoremstyle{remark}

\newtheorem{rem}[thm]{Remark}
\newtheorem{ex}[thm]{Example}
\newtheorem{con}[thm]{Construction}

\DeclareMathOperator{\Aut}{Aut}
\DeclareMathOperator{\Cl}{Cl}
\DeclareMathOperator{\centre}{centre}
\DeclareMathOperator{\Spec}{Spec}
\DeclareMathOperator{\Proj}{Proj}
\DeclareMathOperator{\Def}{Def}
\DeclareMathOperator{\Emb}{Emb}
\DeclareMathOperator{\Bl}{Bl}

\DeclareMathOperator{\Hom}{Hom}
\DeclareMathOperator{\cHom}{\mathcal{H}\mathnormal{om}}
\DeclareMathOperator{\Ext}{Ext}
\DeclareMathOperator{\cExt}{\mathcal{E}\mathnormal{xt}}
\DeclareMathOperator{\cTor}{\mathcal{T}\mathnormal{or}}

\DeclareMathOperator{\dimn}{dim}
\DeclareMathOperator{\Diff}{Diff}
\DeclareMathOperator{\Pic}{Pic}
\DeclareMathOperator{\im}{im}
\DeclareMathOperator{\kernel}{ker}
\DeclareMathOperator{\cok}{coker}
\DeclareMathOperator{\Supp}{Supp}
\DeclareMathOperator{\Exc}{Exc}
\DeclareMathOperator{\ind}{index}
\DeclareMathOperator{\mult}{mult}
\DeclareMathOperator{\Ob}{Ob}

\DeclareMathOperator{\wt}{wt}
\DeclareMathOperator{\red}{red}

\DeclareMathOperator{\Spf}{Spf}

\DeclareMathOperator{\nd}{\not \: \mid}
\DeclareMathOperator{\Hilb}{Hilb}

\DeclareMathOperator{\us}{sp}
\DeclareMathOperator{\Isom}{Isom}
\DeclareMathOperator{\ch}{char}

\newcommand{\QED}{\ifhmode\unskip\nobreak\fi\quad {\rm Q.E.D.}} 

\newcommand{\bA}{\mathbb A}
\newcommand{\bC}{\mathbb C}
\newcommand{\bF}{\mathbb F}
\newcommand{\bH}{\mathbb H}
\newcommand{\bN}{\mathbb N}
\newcommand{\bP}{\mathbb P}
\newcommand{\bQ}{\mathbb Q}
\newcommand{\bR}{\mathbb R}

\newcommand{\bZ}{\mathbb Z}

\newcommand{\cB}{\mathcal B}
\newcommand{\cC}{\mathcal C}
\newcommand{\cX}{\mathcal X}
\newcommand{\cY}{\mathcal Y}
\newcommand{\sD}{\mathcal D}
\newcommand{\cE}{\mathcal E}
\newcommand{\cF}{\mathcal F}
\newcommand{\cG}{\mathcal G}
\newcommand{\sH}{\mathcal H}
\newcommand{\cI}{\mathcal I}
\newcommand{\cO}{\mathcal O}
\newcommand{\sL}{\mathcal L}
\newcommand{\cM}{\mathcal M}
\newcommand{\cN}{\mathcal N}
\newcommand{\cT}{\mathcal T}
\newcommand{\cU}{\mathcal U}

\newcommand{\cW}{\mathcal W}
\newcommand{\cZ}{\mathcal Z}

\newcommand{\map}{\rightarrow}
\newcommand{\da}{\downarrow}
\newcommand{\inj}{\hookrightarrow}

\title{A Compactification of the space of Plane Curves}
\author{Paul Hacking}
\date{}

\begin{document}

\maketitle
\pagestyle{plain}

\tableofcontents

\section{Introduction}

The aim of this thesis is to define a geometric, explicitly computable compactification of the moduli space $V_d$ of smooth
plane curves of degree $d \ge 4$. The basic idea is to regard a plane curve $C \subset \bP^2$ as a log surface $(\bP^2,C)$.
Then there is a compactification given by a moduli space $\cM^1_d$ of log surfaces $(X,D)$ where $K_X+D$ is semi-log-canonical 
and ample, the log analogue of the moduli space of surfaces of general type constructed by Koll\'{a}r and 
Shepherd-Barron \cite{KSB}. 
For the definitions of semi-log-canonical (slc) and semi-log-terminal (slt) singularities see \ref{defn-slc}.
I initially tried to compute these compactifications, with little success --- the $d=4$ case
was calculated by Hassett \cite{Has}, but as $d$ increases the problem quickly becomes unmanageable.

We consider instead the family of compactifications given by moduli spaces $\cM^{\alpha}_d$ of log surfaces $(X,D)$
where $K_X+\alpha D$ is slc and ample, where $\frac{3}{d} < \alpha \le 1$. Note that we require
$\alpha > \frac{3}{d}$ so that $K_{\bP^2}+\alpha C$ is ample. The compactification is simpler for lower $\alpha$
: roughly, the boundary has fewer components, fewer types of degenerate surfaces $X$ occur. The compactification $\cM_d$ 
we compute
below can be described as the limit of $\cM^{\alpha}_d$ as $\alpha$ approaches $\frac{3}{d}$ from above. 
However, we don't define it in this way to avoid technical difficulties and give a clearer presentation.
$\cM_d$ is a moduli space of \emph{stable pairs} $(X,D)$.
Passing from $\cM^1_d$ to $\cM_d$ affords three major simplifications:
\begin{enumerate}
\item For $(X,D)$ a stable pair, $X$ is a slc del Pezzo surface, that is, $-K_X$ is ample and $X$ has slc singularities.\\
\item For $(X,D)$ a stable pair, $dK_X+3D$ is Cartier and linearly equivalent to zero. 
This makes many computations with stable pairs easy 
(e.g., calculating the possible slt singularities on $X$ for $(X,D)$ a stable pair of given degree $d$).\\
\item For $(\cX,\sD)/S$ a relative stable pair, both $K_{\cX/S}$ and $\sD$ are $\bQ$-Cartier
 (whereas for $\cM^{\alpha}_d$ we only know that $K_{\cX/S}+\alpha\sD$ is $\bQ$-Cartier). 
In particular, we need only consider $\bQ$-Gorenstein deformations of $X$ for $(X,D)$ a stable pair.
\end{enumerate}

We now give a map of this thesis. 
In Section \ref{setup} we write down the definition of a stable pair of degree $d$, and define a 
moduli stack  
$\cM_d$ of these objects.
We prove that $\cM_d$ is separated and proper using the valuative criterion and the relative MMP as in \cite{KM}. 
Thus $\cM_d$ gives a compactification of $V_d$.
We delay the discussion of some technical points until Sections~\ref{pfbc}, \ref{QG} and \ref{construct}.
In Sections~\ref{pfbc} and \ref{QG} we explain the definition of an allowable family $(\cX,\sD)/S \in \cM_d(S)$ of stable pairs.
In particular, in Section~\ref{QG} 
we develop a theory of $\bQ$-Gorenstein deformations for slc surfaces which is of independent interest.
We apply these results in Section~\ref{construct} to prove that $\cM_d$ is a Deligne--Mumford stack,
using the methods of Artin \cite{Ar1}.
We remark that we don't need to appeal to Alexeev's results in order to bound the index of $K_X$ here --- we give an explicit
bound using elementary methods.

The aim of the remainder of this thesis is to classify stable pairs. Essentially it is enough to 
classify the degenerate surfaces $X$ occurring in stable pairs $(X,D)$. 
For, given $X$, the divisor $D$ satisfies $dK_X+3D \sim 0$ and $K_X+(\frac{3}{d}+\epsilon)D$ slc for some $0 < \epsilon \ll 1$, 
so $D$ is a member of a given linear system with specified singularities. In Section~\ref{coarse_classn}
 we give a rough classification of the 
surfaces $X$ that occur. They are of 4 types A, B, C, D. Type A are the normal surfaces --- the log terminal cases have already
been classified by Manetti \cite{Ma}. In Section~\ref{normal} 
we show that the only strictly log canonical cases 
are the elliptic cones of degree $9$. Type B are the
nonnormal slt cases --- they have two components meeting in a $\bP^1$. 
We give a classification of these in Section~\ref{B}. Types C and D are the nonnormal strictly
log canonical cases. Type C have a degenerate cusp, type D have a $\bZ / 2\bZ$ quotient of a degenerate cusp, and in both cases
the surface is slt outside these points. 
We remark that there is one more type of slc del Pezzo, 
denoted B*, which is irreducible, slt and has a folded double curve. However we show in Section~\ref{B*} that a surface of type
B* never admits a smoothing to $\bP^2$, so does not occur in a stable pair.

In Section~\ref{3_divis}, we show that if $(X,D)$ is a stable pair of degree $d$ where $3 \nd d$ then $X$ is slt.
The main point is that in the case $3 \nd d$ the condition $dK_X+3D \sim 0$ shows that $K_X$ is 3-divisible
in $\Cl X$ --- this is a strong restriction on $X$. So $X$ is either a Manetti surface or a surface of type B  and we have a 
classification in each case. In particular, $X$ has either $1$ or $2$ components.
We also show that $\cM_d$ is smooth if $3 \nd d$ using deformation theory developed in Sections~\ref{QG} and \ref{construct}.
The case $3 \mid d$ is much more involved --- for example if $d=6$ there is an $X$ with 18 components 
(see Example~\ref{ex-typeC}). Moreover $\cM_d$ is not smooth if $3 \nd d$ (see Example~\ref{ex-elliptic_cone_defns}).

In Section~\ref{small_d} we give the complete classification of the $X$ that occur in degrees 4 and 5, 
and identify the possible singularities
of $D$ on $X$. We also present partial results in the degree 6 case --- we give the complete list of surfaces of types A and B,
and give a list of candidates for the surfaces of types C and D.
We have yet to determine which of these candidates are smoothable.

I'd like to thank my supervisor Alessio Corti, for invaluable guidance throughout the course of my PhD.
Many of the new ideas contained in this thesis grew out of discussions with Alessio, I believe that 
he always had an excellent intuitive feel for the problem which was a great help to me. 
I'm indebted to Brendan Hassett, whose preprint of October 1998 (c.f. \cite{Has}) suggested that one should 
consider the moduli spaces $\cM^{\alpha}_d$ for $\alpha <1$. He also provides an outline of the $\bQ$-Gorenstein deformation 
theory of slc surfaces set out in Section~\ref{QG}. Finally, I'd like to thank Tom Fisher, Jan Wierzba and 
Nick Shepherd-Barron for various helpful conversations.

\section{Definition of the compactification}
\label{setup}

We compactify the space of smooth plane curves $V_d$ of degree $d \ge 4$.
Here $V_d = U_d / \Aut {\bP}^2$, where $U_d$ is the open subscheme of the Hilbert scheme of curves of degree $d$ in ${\bP}^2$
corresponding to smooth curves. We do this by constructing a moduli stack $\cM_d$ of stable pairs $(X,D)$ as defined below.
We use the relative MMP, see \cite{KM} for details. We always work over $\bC$. 

\begin{notn} \label{notn-T}
Let $T$ denote the spectrum of a DVR, $\eta$ the generic point, and $\bar{\eta}$ the geometric generic point.
If $0 \in S$ is a local scheme,  we use script letters to denote flat families over $S$ and 
normal letters to denote the special fibres.
We say $(\cX,\sD)/S$ is a family of pairs over $S$ if $\cX / S$ is a flat family of CM reduced surfaces and $\sD$ is a relative 
Weil divisor on $\cX/S$.
For the definition of a relative Weil divisor, see Section~\ref{pfbc}.  
\end{notn}

\begin{defn} \label{defn-slc}
Let $X$ be a proper connected surface.
Let $D$ be an effective $\bQ$-divisor on $X$.
The pair $(X,D)$ (or $K_X+D$) is \emph{semi-log-canonical} (respectively \emph{semi-log-terminal}) if
\begin{enumerate}
\item $X$ is Cohen--Macaulay and normal crossing in codimension $1$.\\
\item $K_X+D$ is $\bQ$-Cartier.\\
\item Let $\nu \colon X^{\nu} \map X$ be the normalisation of $X$.
Let $\Delta$ be the double curve of $X$ and write $\Delta^{\nu}$, $D^{\nu}$ for the inverse images on $X^{\nu}$.
 Then $K_{X^{\nu}}+\Delta^{\nu}+D^{\nu}$ is log canonical (respectively log terminal).
\end{enumerate}
We use the abbreviations slc and slt for semi-log-canonical and semi-log-terminal.
\end{defn}

\begin{rem}
Note $K_{X^{\nu}}+\Delta^{\nu}+D^{\nu}=\nu^{\star}(K_X+D)$.
\end{rem}

\begin{defn}
Let $X$ be a proper connected surface over an algebraically closed field $k$.
Let $D$ be an effective Weil divisor on $X$. Let $d \in \bN$, $d \ge 4$.
We say that $(X,D)$ is a \emph{semistable pair of degree $d$} if
\begin{enumerate}
\item $X$ is normal and log terminal.\\
\item $K_X+\frac{3}{d}D$ is slc.\\
\item $dK_X+3D \sim 0$, and moreover $\frac{d}{3}K_X+D \sim 0$ if $3 \mid d$.\\
\item There exists a smoothing $(\cX,\sD)/T$ of $(X,D)/k$ such that $\cX_{\bar{\eta}} \cong \bP^2_{\bar{\eta}}$ and 
$K_{\cX}$ and $\sD$ are $\bQ$-Cartier.
\end{enumerate}
\end{defn}

\begin{rem}
We write `$\sim$' to denote linear equivalence of Weil divisors (not $\bQ$-divisors). 
Thus if $B \sim 0$ then in particular $B$ is \emph{Cartier}. This will be important later in the proof of Theorem~\ref{thm-ind}.
We use `$\equiv$' to denote the weaker notion of numerical equivalence (of $\bQ$-Cartier $\bQ$-divisors).
\end{rem}

\begin{rem}
If $(X,D)$ is a semistable pair, $X$ is a normal degeneration of $\bP^2$ with log terminal singularities.
These have been classified by Manetti \cite{Ma}. In particular, $X$ is projective and $-K_X$ is ample.
\end{rem}

\begin{thm}
\label{thm-MMP1}
$T$, $\eta$ as in Notation \ref{notn-T}.
Let $\sD_{\eta} \subset \bP^2_{\eta}$ be a smooth curve of degree $d$ defined over $\eta$.
Then there exists a base change $T' \map T$ and a family $(\cX,\sD)/T'$ of semistable pairs completing 
$(\bP^2_{\eta'},\sD_{\eta'})$ such that $K_{\cX}$ and $\sD$ are $\bQ$-Cartier.  
\end{thm}
We give the proof of the Theorem at the end of this section.

\begin{defn} \label{def-stable}
Let $X$ be a proper connected surface over an algebraically closed field $k$.
Let $D$ be an effective Weil divisor on $X$. Let $d \in \bN$, $d \ge 4$.
$(X,D)$ is \emph{ a stable pair of degree $d$} if
\begin{enumerate}
\item There exists $\epsilon > 0$ such that $K_X+(\frac{3}{d}+\epsilon)D$ is slc and ample.\\
\item $dK_X+3D \sim 0$, and moreover $\frac{d}{3}K_X+D \sim 0$ if $3 \mid d$.\\
\item There exists a smoothing $(\cX,\sD)/T$ of $(X,D)/k$ such that $\cX_{\bar{\eta}} \cong \bP^2_{\bar{\eta}}$, 
and $K_{\cX}$ and $\sD$ are $\bQ$-Cartier.
\end{enumerate}
\end{defn}

\begin{rem}
Note that (1) and (2) imply that $-K_X$ is ample. This is the main reason that stable pairs are easy to classify.
\end{rem}

\begin{thm}
\label{thm-MMP2}
Let $\sD_{\eta} \subset \bP^2_{\eta}$ be a smooth curve of degree $d$.
Then there exists a base change $T' \map T$ and a family $(\cX,\sD)/T'$ of stable pairs completing 
$(\bP^2_{\eta'},\sD_{\eta'})$ such that $K_{\cX}$ and $\sD$ are $\bQ$-Cartier.
\end{thm}
We give the proof of the Theorem at the end of this section.

\begin{thm} \label{thm-ind}
For each $d$ there exists $N \in \bN$ such that for every stable pair $(X,D)$ of degree $d$,
$NK_X$ is Cartier.
\end{thm}

\begin{rem}
It follows by general theory that stable pairs of degree $d$ are bounded for each $d$, 
see the proof of Theorem~\ref{thm-stack}.
\end{rem}

\begin{notn}
Let $N(d)$ denote the least such $N$ for each $d$.
\end{notn}

\begin{defn} \label{defn-allowable}
Let $(\cX,\sD)/S$ be a family of stable pairs of degree $d$. 
That is, $\cX$ is flat over $S$,  $\sD$ is an effective relative Weil divisor on $\cX/S$,
and for every geometric point $s$ of $S$, the fibre $(\cX_s,\sD_{(s)})$ is a stable pair of degree $d$.
We say that $(\cX,\sD)/S$ is an \emph{allowable family} if $\omega_{\cX/S}^{[i]}$ and $\cO_{\cX}(\sD)^{[i]}$
commute with base change for all $i \in \bZ$.
\end{defn}

\begin{rem} \label{rem-QG_and_weaklyQG}
The theory of the conditions `$\omega_{\cX/S}^{[i]}$ and $\cO_{\cX}(\sD)^{[i]}$ commute with base change'
is developed in detail in Sections~\ref{pfbc} and \ref{QG}. 
In particular $K_{\cX/S}$ and $\sD$ are $\bQ$-Cartier since $K_{\cX_s}$ and $\sD_{(s)}$ 
are $\bQ$-Cartier for each closed point $s \in S$, using Lemma~\ref{lem-pfbc}. 
Moreover, if $S$ is the spectrum of a DVR, with generic point $\eta$, and $\cX_{\eta}$ is canonical,
then $(\cX,\sD)/S$ is allowable iff $K_{\cX/S}$ and $\sD$ are $\bQ$-Cartier by Proposition~\ref{prop-pfbc_DVR}.
We shall only require this case below.
\end{rem}

\begin{defn} \label{defn-smoothable_deformation}
Let $(X,D)/k$ be a stable pair of degree $d$.\par
Let $(\cX^u,\sD^u) \map (0 \in D_0)$ be a versal allowable deformation of the pair $(X,D)/k$, where $D_0$ is of finite type over
$k$ (we prove the existence of such a deformation in Section~\ref{construct}). 
Let $D_1 \subset D_0$ be the open subscheme where the geometric fibres of $\cX^u / D_0$ are isomorphic to $\bP^2$.
Let $D_2$ be the scheme theoretic closure of $D_1$ in $D_0$.\par
An allowable deformation of $(X,D)$ is \emph{smoothable} if it is obtained by pullback from the deformation
$(\cX^u, \sD^u) \times_{D_0} D_2 \map (0 \in D_2)$.
\end{defn}   

\begin{rem}
Example~\ref{ex-elliptic_cone_defns} and Example~\ref{ex-rel_smoothability}
describe two examples of this construction.
Given the smoothability assumption for stable pairs $(X,D)/k$, 
we require the relative smoothability assumption for families $(\cX,\sD)/S$ of stable pairs 
in order to obtain an algebraic stack of stable pairs   
(c.f. Example~\ref{ex-rel_smoothability}). 
\end{rem}

\begin{defn}
Let $(\cX,\sD)/S$ be an allowable family of stable pairs of degree $d$. 
$(\cX,\sD)/S$ is a \emph{smoothable family} if for every geometric point $s \in S$ 
the induced deformation of $(\cX_s,\sD_s)$ is smoothable.
\end{defn}

\begin{notn}
Let $\underline{Sch}$ be the category of noetherian schemes over $\bC$. 
\end{notn}

\begin{defn}
We define a groupoid $\cM_d \map \underline{Sch}$ as follows:
$$\cM_d(S) = \left\{ (\cX,\sD)/S \left| \begin{array}{c}	
					\mbox{$(\cX,\sD)/S$ is an allowable smoothable family}\\
			  		\mbox{of stable pairs of degree $d$}
					\end{array} \right. \right\}$$
\end{defn}

\begin{thm}
$\cM_d$ is a separated and proper Deligne--Mumford stack.
\end{thm}

Thus $\cM_d$ gives a compactification $V_d \inj \cM_d$. 
The proof that $\cM_d$ is a Deligne--Mumford stack is given in Section~\ref{construct}.  
We prove that $\cM_d$ is separated and proper in Theorems~\ref{thm-sep} and \ref{thm-proper} below.
We first give the proof of Theorem~\ref{thm-ind}. In fact we prove the following more precise result.

\begin{thm} \label{thm-ind-explicit}
Let $(X,D)$ be a stable pair of degree $d$. Then for all $P \in X$, $\ind_P K_X \le d$ if $3 \nd d$
 and  $\ind_P(K_X) < \frac{2d}{3}$ if $3 \mid d$.
\end{thm}

\begin{proof}
$K_X+(\frac{3}{d}+\epsilon)D$ is slc and $K_X$ is $\bQ$-Cartier by Definiton~\ref{def-stable}, (1) and (2).
Thus $X$ is slc, and $D=0$ at the strictly slc singularities of $X$.
Now $dK_X+3D \sim 0$ and moreover $\frac{d}{3}K_X+D \sim 0$ if $3 \mid d$, 
so the bound is clear at stricly slc points of $X$ .
Since $X$ has a $\bQ$-Gorenstein smoothing, the slt singularities are of 3 types (compare Theorem~\ref{thm-smoothability}):
\begin{enumerate}
\item $X \cong \frac{1}{kn^2}(1,na-1)$ where $(a,n)=1$.
\item $X \cong (xy=0) \subset \bA^3_{x,y,z}  / \frac{1}{r}(1,-1,a)$ where $(a,r)=1$.
\item $X \cong (x^2=zy^2) \subset \bA^3_{x,y,z}$.
\end{enumerate}
We have $\ind_P K_X=n$, $r$, and $1$ in cases 1, 2, and 3.
We now use the two properties $K_X+(\frac{3}{d}+\epsilon)D$ is slc and $dK_X+3D$ is Cartier 
to bound $n$ and $r$ in cases (1) and (2).

In case (1), we first remark that $k=1$ and $3 \nd n$ by Lemma~\ref{lem-T_1}, using the fact that $X$ smoothes to $\bP^2$.  
Consider the local smooth cover $\bA^2_{x,y} \map \bA^2_{x,y} / \frac{1}{n^2}(1,na-1)$ of $X$ at $P$,
let $\tilde{D}$ denote the inverse image of $D$, say $\tilde{P} \mapsto P$. 
Then $\mult_{\tilde{P}}(\tilde{D}) < \frac{2d}{3}$ since $K_X+(\frac{3}{d}+\epsilon)D$ is log canonical at $P$. 
Let $(f(x,y)=0)$ be the local equation of $\tilde{D}$ at $\tilde{P}$ and $x^iy^j \in f$ be a monomial of minimal multiplicity. 
Then $i+j < \frac{2d}{3}$, and $i + (na-1)j = \frac{d}{3}na \bmod n^2$,
using $dK_X+3D$ Cartier and $3 \nd n$.
Thus in particular $i = j \bmod n$. Suppose $n \ge \frac{2d}{3}$, then $i=j$, and $j < \frac{d}{3}$, $j= \frac{d}{3} \bmod n$.
It follows that if $3 \nd d$ we have $n \le d$, and if $3 \mid d $ we have $n \le \frac{d}{3}$. 
This proves the bounds in case (1).

In case (2), let $Y$ be the local analytic component $\bA^2_{x,z} /\frac{1}{r}(1,a)$  of $X$ at $P$, 
let $C \subset Y$ be the double curve, $D_Y = D \mid_Y$.
Then $K_Y+C+(\frac{3}{d}+\epsilon)D_Y$ is log canonical at $P$. Let $\tilde{Y}$ be the smooth cover 
$\bA^2_{x,z} \map \bA^2_{x,z} / \frac{1}{r}(1,a)$ of $Y$, let $\tilde{C},\tilde{D_Y}$ denote the inverse images of $C,D_Y$.
By adjunction $K_{\tilde{C}} + (\frac{3}{d}+\epsilon)\tilde{D_Y} \mid_{\tilde{C}}$ is log canonical, equivalently
$(\frac{3}{d}+\epsilon)\tilde{D_Y} \mid_{\tilde{C}}$ is reduced.
Write $\tilde{D}_Y=(f(x,z)=0)$.
Then $\tilde{D}_Y \mid_{\tilde{C}} = (f(0,z)=0)= (z^k+ \cdots =0) \subset \bA^1_z$ where $k < \frac{d}{3}$.
We have $3k=d \bmod r$ (respectively $k= \frac{d}{3} \bmod r$) if $3 \nd d$ (respectively $3 \mid d$),
using $dK_X+3D$ Cartier (respectively $\frac{d}{3}K_X+D$ Cartier).
It follows that $r \le d$ if $3 \nd d$, $r \le \frac{d}{3}$ if $3 \mid d$.
\end{proof}

We next show that $\cM_d$ is separated. This is almost immediate from the general theory of moduli of surfaces.
I've included my own proof of a foundational result (Lemma \ref{lem-iofa}) because
I could not find a proof in the literature.

\begin{defn}
Let $\cX/T$ be a normal variety, flat over $T$.
Let $\sD \subset \cX$ be an effective Weil divisor, flat over $T$.
Let $\tilde{\cX} \map \cX$ be a projective resolution of $\cX$, write $\tilde{\sD}$ for the strict transform of $\sD$.
Assume that the exceptional locus is a divisor, write $\tilde{\cE}$ for 
the sum of the exceptional divisors dominating $T$.
We say $(\tilde{\cX},\tilde{\sD}+\tilde{\cE}) \map (\cX,\sD)$ is a \emph{semistable log-resolution} of $(\cX,\sD)/T$ if
$\tilde{X}$ is reduced and $\tilde{X} \cup \tilde{\cE} \cup \Supp(\tilde{\sD})$ is a simple normal crossing divisor.
\end{defn} 

\begin{thm} \label{thm-sep}
$\cM_d$ is separated
\end{thm}

\begin{proof} We use the valuative criterion. So, suppose we are given $(\cX^1,\sD^1)/T$ and 
$(\cX^2, \sD^2)/T$ in $\cM_d$ such that $(\cX^1,\sD^1)_{\eta} \cong (\cX^2, \sD^2)_{\eta}$. 
We need to show that $(X^1,D^1) \cong (X^2,D^2)$. Note that we may assume that $\cX^i_{\bar{\eta}} \cong \bP^2_{\bar{\eta}}$ 
and $\sD^i_{\bar{\eta}}$ is smooth for $i=1$ and $2$, since the 
open subset of $\cM_d$ of pairs $(X,D)$ where $X \cong \bP^2$ and $D$ is smooth is dense. 
After base change, we may assume that $\cX^i_{\eta} \cong \bP^2_{\eta}$.

Possibly after base change, we can construct a common semistable 
log resolution  $(\tilde{\cX},\tilde{\sD})$ of $(\cX^1,\sD^1)/T$ and $(\cX^2,\sD^2)/T$ that is an isomorphism over $\eta$.
 We now claim that we can 
reconstruct the two families as the $K_{\tilde{\cX}}+\tilde{X}+(\frac{3}{d}+\epsilon)\tilde{\sD}$
canonical model of $(\tilde{\cX},\tilde{\sD})/T$ for $0 < \epsilon \ll 1$.
So $(\cX^1,\sD^1) \cong (\cX^2,\sD^2)$ and $(X^1,D^1) \cong (X^2,D^2)$ as required. 
To prove the claim, we just need to verify that
$K_{\cX^i}+X^i+(\frac{3}{d}+\epsilon)\sD^i$ is log canonical and relatively ample for $i=1,2$.
We can check ampleness on the central fibre. We have 
$K_{\cX^i}+X^i+(\frac{3}{d}+\epsilon)\sD^i \mid_{X^i} =K_{X^i}+(\frac{3}{d}+\epsilon)D^i$ by adjunction.
This is ample since $K_{X^i}+(\frac{3}{d}+\epsilon^i)D^i$ is ample for
some $\epsilon^i >0$ and $K_{X^i}+ \frac{3}{d}\sD^i \equiv 0$. 
Now, $K_{X^i}+(\frac{3}{d}+\epsilon^i)D^i$ is slc 
for some $\epsilon^i >0$ by the definition of $\cM_d$, and the same is true for $\epsilon < \epsilon^i$.
It follows that $K_{\cX^i}+X^i+(\frac{3}{d}+\epsilon)\sD^i$ is log canonical by Lemma~\ref{lem-iofa}. This completes the proof.
\end{proof}

\begin{lem}
\label{lem-iofa}
Let $\cX/T$ be a flat family of surfaces, $\cB \subset \cX$ a $\bQ$-divisor that is flat over $T$.
Suppose $\cX_{\eta}$ is normal and $K_{\cX}+\cB$ is $\bQ$-Cartier. 
Suppose also that $(\cX,\cB)/T$ admits a semistable log-resolution.
Then $(\cX,X+\cB)$ is log canonical iff $(X,B)$ is slc.
\end{lem}

\begin{proof}
If $(\cX,X+\cB)$ is log canonical, then trivially  
$(X,B)$ is slc by
adjunction. We now prove the converse. We may work locally at $P \in \cX$.
 
Take a semistable log resolution $(\tilde{\cX},\tilde{\cB}+\tilde{\cE})/T$ of $(\cX,\cB)/T$.
Let $(\hat{\cX},\hat{\cB}+~\hat{\cE})/T$ be the $K_{\tilde{\cX}}+\tilde{\cB}+\tilde{\cE}$ canonical model over 
$(\cX,\cB)$. Then $f: (\hat{\cX},\hat{\cB}+\hat{\cE}) \map (\cX,\cB)$ 
is an isomorphism over the locus where $(\cX,X+\cB)$ is log canonical.
$(\cX,X+\cB)$ is log canonical at 
any codimension $1$ point of $X$. For the implication $(X,B)$ slc $\Rightarrow$ $(\cX,X+\cB)$ log canonical is easy in the case 
$\dimn X =1$. So no exceptional divisor of $f$ has centre a codimension $1$ point of $X$.
 Now let $V$ be the normalisation of a component of $X$,
$\Delta_V$ the inverse image of the double curve of $X$ on $V$, and $B_V = \cB \mid_V$. 
Let $W$ be the normalisation of the strict transform of this component in $\hat{\cX}$, $\Delta_W$ the inverse image of the double 
curve of $\hat{X}$ on $W$, $\hat{B}_W= \hat{\cB} \mid_W$ and $\hat{E}_W=\hat{\cE} \mid_W$.
Write $g: W \map V$ for the induced map. We have
$$K_{W}+\Delta'_{V}+B'_V= g^{\star}(K_V+\Delta_V+B_V)+\sum a_iF_i,$$
where the primes denote strict transforms and the $F_i$ are the $g$-exceptional divisors.
 Then $a_i \ge -1$ for all $i$ because $(X,B)$ is slc.
Now $g_{\star}(\Delta_W+\hat{B}_W+\hat{E}_W)=\Delta_V+B_V$ since there are no $f$-exceptional divisors with centre a codimension
$1$ point of $X$.
So, we have  
$$K_{\hat{\cX}}+\hat{\cB}+\hat{\cE} \mid_W = K_{W}+\Delta_{W}+\hat{B}_W+\hat{E}_W= g^{\star}(K_V+\Delta_V+B_V)+\sum b_iF_i,$$
where $b_i < 0$ for all $i$ since $K_{\hat{\cX}}+\hat{\cB}+\hat{\cE}$ is relatively ample.
If $F_j \subset \Delta_W+\hat{E}_W$ for some $j$, we have $b_j \ge a_j+1 \ge 0$, a contradiction.
So, in particular, for each component $X_i$ of $X$, with strict transform $X_i'$ in $\hat{\cX}$, 
the map $X_i' \map X_i$ does not contract any component of the double curve of $\hat{X}$ on $X_i'$.
It follows that there are no $f$-exceptional divisors over $0 \in T$ (recall that every such divisor has centre $P$).
Moreover, for each component $X_i$ of $X$, $\hat{\cE} \mid_{X'_i}=0$.
Hence $\hat{\cE}=0$, so $f$ has no exceptional divisors. 
Then $K_{\hat{\cX}}+\hat{\cB}= f^{\star}(K_{\cX}+\cB)$ and $K_{\hat{\cX}}+\hat{\cB}$ $f$-ample implies $f$ is an isomorphism,
so $(\cX,X+\cB)$ is log canonical as required.
\end{proof}

\begin{rem}
The Lemma above is similar to \cite{KSB}, p.~325, Theorem 5.1(a). This is an inversion of adjunction type result 
(compare \cite{KM}, p.~174, Theorem~5.50), but we can't use the general theorem to prove our Lemma.
\end{rem}

We now give the proof that $\cM_d$ is proper. The main steps of the argument are in Theorem~\ref{thm-MMP1} and 
Theorem~\ref{thm-MMP2}. These theorems motivate the definition of a stable pair above.

\begin{proof}[Proof of theorem \ref{thm-MMP1}]
First complete $(\bP^2_{\eta},\sD_{\eta})$ to a family $(\bP^2_T,\sD')/T$.
Possibly after base change (which we suppress in our notation) we can take a semistable log-resolution
$(\tilde{\cX},\tilde{\sD})/T$ of $(\bP^2_T,\sD')/T$ which is an isomorphism over $\eta$. 
Then $K_{\tilde{\cX}}+\tilde{X}+\frac{3}{d}\tilde{\sD}$ is dlt.
Run a $K_{\tilde{\cX}}+\frac{3}{d}\tilde{\sD}$ MMP over $T$. Let $(\hat{\cX},\hat{\sD})/T$ be the end product.
Then $K_{\hat{\cX}}+\frac{3}{d}\hat{\sD}$ must be relatively nef, since it is zero on the generic fibre.
Now it follows that $dK_{\hat{\cX}}+3\hat{\sD} \sim 0$. For we can write $dK_{\hat{\cX}}+3\hat{\sD} \sim \sum a_i\hat{X}_i$, 
where the $\hat{X}_i$ are the components of $\hat{X}$, since $dK_{\hat{\cX}}+3\hat{\sD}$ is zero on the generic fibre.
Without loss of generality we may assume $a_1=0$ and $a_i \le 0$ for $i >1$, using the relation $\sum \hat{X}_i \sim 0$.
Restricting to $\hat{X}_1$, we see that $K_{\hat{\cX}}+\frac{3}{d}\hat{\sD}$ relatively nef implies that $a_i=0$ for all $i$
such that $\hat{X}_i \cap \hat{X}_1 \neq \emptyset$. Since $\hat{X}$ is connected, repeating the argument we obtain $a_i=0$
for all $i$,
so $dK_{\hat{\cX}}+3\hat{\sD} \sim 0$ as claimed. Compare Lemma~\ref{lem-Cl}, (1).
We have $K_{\hat{\cX}}+\hat{X}+\frac{3}{d}\hat{\sD}$ dlt and $\hat{\cX}$ is $\bQ$-factorial, so
$K_{\hat{\cX}}+\hat{X}$ is dlt.
We now run a $K_{\hat{\cX}}$ MMP over $T$ and let $(\cX,\sD)/T$ be the end product.
Then $\cX/T$ is a del Pezzo fibration, since the generic fibre is a del Pezzo surface (namely $\bP^2$).
Now, $\rho(\cX/T)=1$ and $\cX$ $\bQ$-factorial gives $X$  irreducible --- for we have an exact sequence
$$ 0 \map \frac{\oplus \bQ X_i}{(\sum X_i=0)} \map N^1(\cX) \map N^1(\bP^2_{\eta}) \map 0, $$
by Lemma~\ref{lem-Cl}, where the $X_i$ are the irreducible components of $X$.
The dlt property of $K_{\cX}+X$ now gives that $X$ is normal and log terminal.
Since $K_{\hat{\cX}}+\hat{X}+\frac{3}{d}\hat{\sD}$ is dlt and $dK_{\hat{\cX}}+3\hat{\sD} \sim 0$,
we have $K_{\cX}+X+\frac{3}{d}\sD$ log canonical and $dK_{\cX}+3\sD \sim 0$.
We obtain $K_{X}+\frac{3}{d}D$ is log canonical and $dK_{X}+3D \sim 0$ by adjunction.
The sharpening in the case $3 \mid d$ is clear.
Finally, since $K_{\cX}$ is $\bQ$-Cartier and $dK_{\cX}+3\sD \sim 0$, $\sD$ is $\bQ$-Cartier.
\end{proof}

\begin{proof}[Proof of theorem \ref{thm-MMP2}]
First complete $(\bP^2_{\eta},\sD_{\eta})$ to a family $(\cX',\sD')/T$.
By  Theorem~\ref{thm-MMP1} and its proof, we may assume that 
$K_{\cX'}+{X'}+\frac{3}{d}\sD'$ is log canonical, 
and $dK_{\cX'}+3\sD' \sim 0$. 

After base change, take a semistable log resolution $(\tilde{\cX},\tilde{\sD}) \map (\cX',\sD')$ 
which is an isomorphism over $\eta$.
Now run a $K_{\tilde{\cX}}+(\frac{3}{d}+\epsilon)\tilde{\sD}$ MMP over $(\cX',\sD')$.
Here $0 < \epsilon \ll 1$ is chosen such that for $F$ an exceptional divisor of $\tilde{\cX} \map \cX$, if
$a(F,\cX',X'+(\frac{3}{d}+\epsilon)\sD') \le -1$ then 
$a(F,\cX',X'+\frac{3}{d}\sD') = -1$.
Let $f \colon (\hat{\cX},\hat{\sD}) \map (\cX',\sD')$ be the end product. 
By the choice of $\epsilon$, we have $dK_{\hat{\cX}}+3\hat{\sD} = f^{\star}(dK_{{\cX'}}+3{\sD'}) \sim 0$.

Now take the $K_{\hat{\cX}}+(\frac{3}{d}+\epsilon)\hat{\sD}$ canonical model over $T$, denote this 
$(\cX,\sD)/T$.
Then $K_{\cX}+X+(\frac{3}{d}+\epsilon)\sD$ is log canonical and relatively ample over $T$,
 and $dK_{\cX}+3\sD \sim 0$ since $dK_{\hat{\cX}}+3\hat{\sD} \sim 0$.
We obtain $K_{X}+(\frac{3}{d}+\epsilon)D$ is slc and ample and 
$dK_{X}+3D \sim 0$  by adjunction. The sharpening in the case $3 \mid d$ is obvious.
Finally since $dK_{\cX}+3\sD \sim 0$ and 
$K_{\cX}+(\frac{3}{d}+\epsilon)\sD$ is $\bQ$-Cartier, $K_{\cX}$ and $\sD$ are $\bQ$-Cartier.
\end{proof}

\begin{thm} \label{thm-proper}
$\cM_d$ is proper
\end{thm}
\begin{proof}
We use the valuative criterion of properness.
Taking $(\cX_{\eta},\sD_{\eta})/ \eta \in \cM_d$, we need to show that, possibly after base change, 
there exists an extension $(\cX,\sD)/T \in \cM_d$. 
We may assume that $\cX_{\eta} \cong \bP^2_{\eta}$  and $\sD_{\eta}$ is smooth as in the proof of separatedness above.
Applying Theorem~\ref{thm-MMP2} we obtain our result.
\end{proof}

\section{A coarse classification of the degenerate surfaces}
\label{coarse_classn}

We now work towards a classification of the surfaces $X$ that occur in stable pairs $(X,D)$. We collect the relevant 
results from the above in the proposition below.

\begin{prop}
Let $(X,D)$ be a stable pair. Then $-K_X$ is ample, $X$ is slc, and $X$ admits a $\bQ$-Gorenstein smoothing.
\end{prop}

\begin{proof}
$-K_X$ is ample by Definition~\ref{def-stable}, (1) and (2). 
There exists a $\bQ$-Gorenstein smoothing $(\cX,\sD)/T$ of $(X,D)$ by Definition~\ref{def-stable}, (3).
Finally, $K_X+(\frac{3}{d}+\epsilon)D$ is slc, and
$K_X$ is $\bQ$-Cartier, so $X$ is slc.
\end{proof} 

\begin{notn}
Let $X$ be a surface with normal crossing singularities in codimension $1$.
Let $\Delta \subset X$ denote the double curve of $X$.
Write $X= \bigcup X_i$ for the decomposition of $X$ into its irreducible components.
Let $\Delta_i =\Delta \cap X_i$.
Let $\nu : X^{\nu} \map X$ be the normalisation of $X$. 
Write $\Delta^{\nu}$ for the inverse image of $\Delta$. 
Then $(X^{\nu},\Delta^{\nu})= \amalg (X^{\nu}_i,\Delta^{\nu}_i)$, where $X^{\nu}_i \map X_i$ is the normalisation and 
$\Delta^{\nu}_i$ is the inverse image of $\Delta_i$.
\end{notn}

First we give a rough classification of the components $(X^{\nu}_i,\Delta^{\nu}_i)$. 
Then we describe how to glue these components together to recover $X$.

\begin{thm}(\cite{KM} p119 Theorem 4.15)\\
\label{thm-lc_pairs}
Let $(0 \in Y, C)$ be a log canonical pair, where $Y$ is a surface and $C$ is an effective Weil divisor. Assume $C \neq 0$.
Then we have the following cases:\\
(a) $(\bA^2_{x,y}/ \frac{1}{r}(1,a), (x=0))$, (a,r)=1.\\
(b) $(\bA^2_{x,y} / \frac{1}{r}(1,a), (xy=0))$, (a,r)=1.\\
(c) $(\bA^2_{x,y} / \frac{1}{r}(1,a), (xy=0)) / \mu_2$, (a,r)=1, 
where the $\mu_2$-action is etale in codimension $1$ and interchanges $(x=0)$ and $(y=0)$.
\end{thm}

\begin{notn}
We will denote cases (a), (b) and (c) by $(\frac{1}{r}(1,a),\Delta)$, $(\frac{1}{r}(1,a),2 \Delta)$ and $(D,\Delta)$ respectively.
The $D$ stands for dihedral --- these singularities generalise the dihedral Du Val singularities.
\end{notn}

\begin{thm}
\label{thm-cpts}
Let $(Y,C)$ be a proper log canonical pair with $-(K_Y+C)$ ample.
Then $(Y,C)$ belongs to one of the following types:
\begin{description}
\item[(I)]	$C=0$.\\
     	      	Then $Y$ has at most one strictly log canonical singularity. 
\item[(II)]  	$C \cong {\bP}^1$ and $(Y,C)$ log terminal.\\ 
             	Then $(Y,C)$ has singularities $(\frac{1}{r_i}(1,a_i),\Delta)$ at $C$, 
		with $\sum (1-\frac{1}{r_i}) < 2$.\\
\item[(III)] 	$C \cong {\bP}^1 \cup {\bP}^1$, where the two components intersect in a node of $C$.\\
		Then $(Y,C)$ has a singularity of type $(\frac{1}{r}(1,a), 2\Delta)$
		at the node of $C$, and at most one other singularity of type 
		$(\frac{1}{r}(1,a), \Delta)$ on each component of $C$.
		Moreover $Y$ has log terminal singularities away from $C$.\\
\item[(IV)]     $C \cong {\bP}^1$ and $(Y,C)$ has a singularity 
		of type $(D,\Delta)$.\\
		Then $(Y,C)$ has at most one other singularity of type 
		$(\frac{1}{r}(1,a),\Delta)$ at $C$ and
		$Y$ has log terminal singularities away from $C$.
\end{description}
\end{thm}

To prove the theorem, we need the connectedness theorem of Shokurov:

\begin{thm}(\cite{KM} p174 Corollary 5.49)\\
\label{thm-connectedness} 
Let $X$ be a normal, proper variety, $B \subset X$ a Weil divisor, and suppose that $-(K_X+B)$ is nef.
Then the locus of log canonical singularities, i.e., $\bigcup \{ \centre_X E \}$ taken over all $E$ with discrepancy
$a(E,X,B) \le -1$, is connected.
\end{thm}

\begin{proof}[Proof of Theorem \ref{thm-cpts}]
First, by Theorem~\ref{thm-connectedness}, either $C=0$ and there is at most one strictly log canonical singularity on $Y$,
or $C \neq 0$, $C$ is connected, and $Y$ is log terminal away from $C$.
Now let $\Gamma$ be a component of $C$. We have $-(K_Y+C)$ is ample, so $(K_Y+C) \Gamma < 0$. Expanding the left hand side,
$$ 2p_a(\Gamma)-2+ (C- \Gamma)\Gamma+ \Diff(\Gamma,Y) < 0.$$
Now $(C-\Gamma)\Gamma \ge 0$ and $\Diff(\Gamma,Y) \ge 0$, thus $p_a(\Gamma)=0$, $\Gamma \cong \bP^1$, and 
$(C- \Gamma)\Gamma+ \Diff(\Gamma,Y) < 2$. We know that $(Y,C)$ is log canonical, so we have the classification given in 
Theorem~\ref{thm-lc_pairs} of the singularities of $(Y,C)$ at $C$. We calculate that each point of type (a), (b), (c) on $\Gamma$
contributes $1-\frac{1}{r}$, $1$, $1$ to $(C- \Gamma)\Gamma+ \Diff(\Gamma,Y)$ respectively. The theorem now follows easily. 
\end{proof}

We now want to give the classification of the surfaces $X$, by glueing together components as above.
We first state the classification of nonnormal slc surface singularities.

\begin{notn}
$\Delta^{\nu} \map \Delta$ is 2-to-1. Let $\Gamma \subset \Delta$ be a component,
 write $\Gamma^{\nu}$ for the inverse image in $\Delta^{\nu}$. Then either $\Gamma^{\nu}$ has two components mapping birationally to
$\Gamma$, or $\Gamma^{\nu}$ has one component mapping 2-to-1 to $\Gamma$.
In the latter case we say $\Gamma \subset \Delta$ is obtained by \emph{folding} $\Gamma^{\nu} \subset \Delta^{\nu}$. 
If we are working locally at $Q \in X$ and we are in the latter case, let $P \mapsto Q$. 
We say $\Gamma^{\nu}$ is \emph{pinched at} $P$.
\end{notn}

\begin{thm} (see \cite{KSB}, Section~4)\\
\label{thm-slc}
Let $0 \in X$ be a slc singularity, assume that $X$ is non-normal. Then we have the following cases:
\begin{enumerate}
\item $(xy=0) \subset \bA^3_{x,y,z} / \frac{1}{r}(a,b,1)$, where $(a,r)=(b,r)=1$.\\
\item $(xy=0) \subset \bA^3_{x,y,z} / \mu_r$, where, for $\zeta$ a generator of $\mu_r$, we have 
$x \mapsto \zeta^a y$, $y \mapsto x$, $z \mapsto \zeta z$, and $4 \mid r$, $(a,r)=2$.\\
\item $(x^2=zy^2) \subset \bA^3_{x,y,z} / \frac{1}{r}(1+a,a,2)$ where $r$ odd, $(a,r)=1$.\\
\item $0 \in X$ is a degenerate cusp: $X = \bigcup X_i$, with $(0 \in X^{\nu}_i,\Delta^{\nu}_i)$ of type
$(\frac{1}{r}(1,a),2\Delta)$, glued to $X^{\nu}_{i-1}$ along one component of $\Delta^{\nu}_i$ and $X^{\nu}_{i+1}$ 
along the other component, 
so that we have a cycle of components.\\
\item $0 \in X$ is a $\mu_2$-quotient of a degenerate cusp: $X=X_1 \cup \cdots \cup X_k$, where 
for $2 \le i \le k-1$, 
$(0 \in X^{\nu}_i,\Delta^{\nu}_i)$ is of 
type $(\frac{1}{r}(1,a),2\Delta)$ and is glued to $X_{i-1}$ along one component of $\Delta^{\nu}_i$ and 
$X_{i+1}$ along the other component.\par 
For $i=1$ or $k$ there are two cases. Either $(0 \in X^{\nu}_1,\Delta^{\nu}_1)$ is of type $(D,\Delta)$ and $X^{\nu}_1$ is glued to $X^{\nu}_2$ along the single component of $\Delta^{\nu}_1$. 
Or $(0 \in X^{\nu}_1,\Delta^{\nu}_1)$ is of type  $(\frac{1}{r}(1,a), 2 \Delta)$,   
$X^{\nu}_1$ is glued to $X^{\nu}_2$ along one component of $\Delta^{\nu}_1$ and the other component is pinched at $0$.   
Similarly for $i=k$. So we have a chain of components.
\end{enumerate}
\end{thm}

\begin{rem}(see \cite{KSB}, p. 326, Remark~5.2)
\label{rem-smoothability}
The following cases admit a $\bQ$-Gorenstein smoothing:\\
(1) $a+b=0 \bmod r$.\\
(2) None.\\
(3) $r = 1$.\\
(4) All.\\
(5) Unknown.\\
\end{rem}

\begin{thm}
\label{thm-glueing}
Let $X$ be a slc proper surface with $-K_X$ ample.
Then $X$ belongs to one of the following types (we use the notation of Theorem~\ref{thm-cpts} to describe the components):
\begin{description}
\item[(A)] 	$X$ normal.\\
		$(X,0)$ is a surface of type I.
\item[(B)]      $X=X_1 \cup X_2$, slt.\\ 
		$(X_1,\Delta_1),(X_2,\Delta_2)$ are surfaces of type II, 
		glued along $\Delta \cong \bP^1$.
\item[(B*)]     $X$ irreducible, non-normal, slt.\\ 
		$(X^{\nu},\Delta^{\nu})$ is a surface of type II, 
		$X$ is formed by folding $\Delta^{\nu}$.
\item[(C)]	$X$ has a degenerate cusp.\\
		$X$ is a union of components of type III which are glued together to form a cycle. 
		More precisely, $X= \bigcup X_i$, where $X^{\nu}_i$ is glued to $X^{\nu}_{i-1}$ and $X^{\nu}_{i+1}$ along the 
		two components of $\Delta^{\nu}_i$ for all $i$, in such a way that the nodes of the curves $\Delta^{\nu}_i$ 
		all coincide --- this point is the degenerate cusp.  
\item[(D)]	$X$ has a $\mu_2$ quotient of a degenerate cusp.\\
		$X$ has some components of type III glued together to form a chain. 
		At each end we glue on either a component of type IV or a component of type III with one component of 
		$C$ folded. The nodes of the curves $\Delta_i^{\nu}$ and the dihedral singularities on the components of 
		type IV all coincide --- this point is the $\mu_2$-quotient of a degenerate cusp.

\end{description}

Here, the singularities of type $(\frac{1}{r}(a,1),\Delta)$ on $(X^{\nu},\Delta^{\nu})$ are glued together together in pairs
$\frac{1}{r}(a,1)$ and $\frac{1}{r}(b,1)$ to give a singularity of type $(xy=0) \subset \frac{1}{r}(a,b,1)$ on $X$.

\end{thm}

\begin{proof}[Proof of theorem \ref{thm-glueing}]
The classification in the Theorem is immediate from the classification of the components $(X^{\nu}_i,\Delta^{\nu}_i)$
in  Theorem~\ref{thm-cpts} using the classification of slc singularities in Theorem~\ref{thm-slc}.
\end{proof}

\begin{rem}
We remark that a surface of type B* does not admit a $\bQ$-Gorenstein smoothing to $\bP^2$ --- see Theorem~\ref{thm-B*}.
\end{rem}

\begin{thm}
\label{thm-smoothability}
Notation as in theorem \ref{thm-glueing}.\\
If $X$ admits a $\bQ$-Gorenstein smoothing, we have the following additional conditions:
\begin{enumerate}
\item The normal, log terminal singularities of $X$ are cyclic quotient singularities of the form $\frac{1}{dn^2}(1,dna-1)$.\\
\item The singularities of type $(\frac{1}{r}(a,1),\Delta)$ on $(X^{\nu},\Delta^{\nu})$ are glued together together in pairs
$\frac{1}{r}(a,1)$ and $\frac{1}{r}(-a,1)$ to give a singularity of type $(xy=0) \subset \frac{1}{r}(a,-a,1)$ on $X$.\\
\item If $\Gamma^{\nu} \subset \Delta^{\nu}$ is pinched at $P$, then either $( P \in X^{\nu},\Delta^{\nu})$ is of type 
$(\frac{1}{r}(1,a),2\Delta)$ or $P \in X^{\nu}$ is smooth.
\end{enumerate}
\end{thm}

\begin{proof}
(1) is the well known classification of smoothable log terminal surface singularities 
(\cite{KSB}, p.~313 Propn~3.10).
Conditions (2) and (3) follow from Remark~\ref{rem-smoothability}.
\end{proof}

\begin{rem}
If $X$ admits a $\bQ$-Gorenstein smoothing to $\bP^2$, then in (1) we have $d=1$ and $3 \not \: \mid n$.
See Lemma~\ref{lem-T_1}.
\end{rem}

\section{Rationality of the components of $X$}
\label{rationality}

In this section we prove that all the components of $X$ are rational unless $X$ is an elliptic cone.

\begin{notn} \label{notn-cpt}
Let $X$ be a slc proper surface with $-K_X$ ample.
Let $(Y,C) \subset (X^{\nu},\Delta^{\nu})$ be a component.
Let $\pi \colon \tilde{Y} \map Y$ be the minimal resolution of $Y$.
Write $K_{\tilde{Y}}+\tilde{C}= \pi^{\star}(K_Y+C)$, where $\pi_{\star} \tilde{C} =C$.
Let $\phi \colon \tilde{Y} \map \bar{Y}$ be a minimal model of $\tilde{Y}$. Note that $Y$ is birationally ruled 
(see below), so $\bar{Y}$ is either ruled or $\bP^2$.
Write $\bar{C}=\phi_{\star} \tilde{C}$.
\end{notn}

\begin{prop} Notation as above.  \label{prop-setup}
\begin{enumerate}
\item $\tilde{C}$ is an effective  $\bQ$-divisor with components of multiplicity $\le 1$.\\
\item $-(K_{\tilde{Y}}+\tilde{C})$ is nef and big, and is zero exactly on $\Exc \pi $.\\
\item $-(K_{\bar{Y}}+\bar{C})$ is nef and big.
\end{enumerate}
\end{prop}

\begin{proof}
Write $K_{\tilde{Y}}+C' = \pi^{\star}(K_Y+C) + \sum a_iE_i$ where the $E_i$ are the exceptional divisors of $\pi$,
and $C'$ is the strict transform of $C$. Then $\tilde{C} = C' + \sum (-a_i)E_i$.
 Here $a_i \ge -1$ for all $i$ since $(Y,C)$ is log canonical.
$\pi$ is minimal, so $K_{\tilde{Y}}$ is $\pi$-nef, hence also $K_{\tilde{Y}}+C'$ is $\pi$-nef.
It follows that $a_i \le 0$ for all $i$. This proves (1). (2) is immediate since $K_{\tilde{Y}}+\tilde{C}= \pi^{\star}(K_Y+C)$
by definition and $-(K_Y+C)$ is ample. (3) now follows since $(K_{\bar{Y}}+\bar{C})=\phi_{\star}(K_{\tilde{Y}}+\tilde{C})$
\end{proof}

\begin{thm}
\label{thm-min-res}
Either $Y$ is rational, or $Y$ is an elliptic cone and $X=Y$.
\end{thm}

\begin{proof}
First, we claim $\tilde{Y}$ is birationally ruled. For $-(K_{\tilde{Y}}+\tilde{C})$ is nef and big, so
$h^0(n(K_{\tilde{Y}}+\tilde{C}))=0$ for all $n > 0$. $\tilde{C}$ is effective, so $h^0(nK_{\tilde{Y}})=0$ for all $n > 0$.
It follows that $\tilde{Y}$ is birationally ruled. We may assume that $\tilde{Y}$ is not rational.
Then $\bar{Y}$ is ruled over a curve of positive genus, let $\bar{Y} \map B$ be a ruling and let $p$ denote the composite   
$p: \tilde{Y} \map \bar{Y} \map B$.
 
Suppose that there is no horizontal component of $\tilde{C}$. Thus $\bar{C} \subset \bar{Y}$ is a sum of fibres.
Now $-(K_{\bar{Y}}+\bar{C})$  nef and big implies that  $-K_{\bar{Y}}$ is nef and big.
Then $h^1(\cO_{\bar{Y}})=h^1(K_{\bar{Y}})=0$ by Serre duality and Kodaira vanishing, a contradiction.

Thus there is an irrational component of $\Supp{\tilde{C}} \subset C' \cup \Exc(\pi)$.
By Theorem~\ref{thm-cpts} we know that $C$ has only rational components,
so by the classification of log canonical singularities it follows that $Y$ has a simple elliptic singularity.
Let $E$ denote the $\pi$-exceptional elliptic curve on $\tilde{Y}$.
Then $E$ has multiplicity $1$ in $\tilde{C}$ and $E$ is horizontal.
Now $-(K_{\tilde{Y}}+\tilde{C})$ is big so $-(K_{\tilde{Y}}+\tilde{C}).f > 0$ for $f$ a fibre of the ruling,
thus $E \cdot f=1$, $E$ is a section. Next we claim that $\tilde{Y}$ is in fact ruled. Suppose not, 
then there exists a degenerate fibre.
Let $A$ be a component meeting $E$.
Then $A$ is not contained in $\Supp \tilde{C}$.
We have $(K_{\tilde{Y}}+\tilde{C}) \cdot A \le 0$, with equality iff $A$ is contracted by $\pi$.
But also 
$$(K_{\tilde{Y}}+\tilde{C}) \cdot A \ge K_{\tilde{Y}} \cdot A +E \cdot A \ge -1 + 1 =0,$$
with equality only if $A$ is a $-1$ curve. 
Thus $A$ is a $-1$ curve and is contracted by $\pi$, a contradiction since $\pi$ is minimal.
So $\tilde{Y}$ is ruled over an elliptic curve. $Y$ is obtained from $\tilde{Y}$ by contracting the negative section 
and so $Y$ is an elliptic cone. Finally $C=0$ by Shokurov's connectedness theorem (compare Theorem~\ref{thm-cpts}) so $X=Y$.
\end{proof}

\section{Bounding $\rho$ of the components of $X$}
\label{rho}

We prove bounds on the values $\rho(Y)$ for $Y$ a component of $X$, using the existence of a smoothing to $\bP^2$.
We use these in Theorem~\ref{thm-Manetti} and Theorem~\ref{thm-B} 
to give necessary and sufficient conditions for surfaces of types A and B
to be smoothable to $\bP^2$. For types C and D we do not have necessary and sufficient criteria for smoothability as yet. 
However the bounds we give here substantially
simplify the explicit calculations we do for $d=6$ in Section~\ref{small_d} (note types C and D only occur if $3 \mid d$
--- see Theorem~\ref{thm-simp}).

\begin{thm}
\label{thm-rho}
Let $X$ be a slc proper surface with $-K_X$ ample, and $\cX/T$ a  
smoothing with $\cX_{\eta} \cong \bP^2_{\eta}$.\\
(1) $\sum \rho (X^{\nu}_i) \le V+E$, with equality only if $\cX$ is $\bQ$-factorial.\\
Here $V=$ number of components of $X$ and $E=$ number of components of $\Delta$, not 
counting components obtained by folding a component of $\Delta^{\nu}$.\\
(2) If $X_1 \cap (X - X_1) \cong \bP^1$, $\rho(X^{\nu}_1) \le 2$.\\
(3) If $X_1 \cap (X - X_1) \cong \bP^1 \cup \bP^1$, $\rho(X^{\nu}_1) \le 4$.\\
\end{thm}

\begin{cor}
\label{cor-rho}
Let $X$ be a slc proper surface with $-K_X$ ample, and $\cX/T$ a $\bQ$-Gorenstein smoothing with 
$\cX_{\eta} \cong \bP^2_{\eta}$.\\ 
In the cases enumerated in Theorem~\ref{thm-glueing}, we have the following bounds for $\rho(X^{\nu}_i)$.
\begin{description}
\item[(A)] 	$\rho(X)=1$.
\item[(B)]	Either $\rho(X_1)=\rho(X_2)=1$, or $\{ \rho(X_1),\rho(X_2) \}=\{ 1,2 \}$ and $\cX$ is $\bQ$-factorial.
\item[(C)]	$\rho(X_i) \le 4$ for all $i$, and $\sum \rho(X_i) \le 2V$, equal only if $\cX$ is $\bQ$-factorial.
\item[(D)]	$\rho(X_i) \le 3$ for $X_i$ a middle component, $\rho(X^{\nu}_i) \le 2$ for $X_i$ an end component,
		and $\sum \rho(X^{\nu}_i) \le 2V-1$, equal only if $\cX$ is $\bQ$-factorial.
\end{description}
\end{cor}

\begin{rem} 
Case B* does not occur by Theorem~\ref{thm-B*}.
\end{rem}

\begin{lem}
\label{lem-Cl}
Let $X$ be a proper surface, normal crossing in codimension $1$, and $\cX/T$ a smoothing with $\cX_{\eta} \cong \bP^2_{\eta}$.\\
(1) We have an exact sequence
$$0 \map \frac{\oplus \bZ X_i}{(\sum X_i =0)} \map \Cl(\cX) \map \Cl(\bP^2_{\eta}) \map 0.$$
(2) Assume in addition that $X$ is projective. Then $\Pic(\cX) \otimes_{\bZ}\bQ \map N^1(\cX)$ is an isomorphism.
\end{lem}

\begin{proof}
(1) Certainly 
$$\frac{\oplus \bZ X_i}{(\sum X_i =0)} \map \Cl(\cX) \map \Cl(\bP^2_{\eta}) \map 0$$ 
is exact, so we just need to show the first map is injective. Suppose $\sum a_i X_i \sim 0$, then using $\sum X_i \sim 0$,
we can replace this with a relation $\sum b_i X_i \sim 0$, where, without loss of generality, $b_1=0$ and $b_i \ge 0$ for $i>1$.
Then $(\sum b_i X_i) \mid_{X_1}$ is effective and linearly equivalent to $0$. So $(\sum b_i X_i) \mid_{X_1} =0$ since $X_1$ is 
proper. Thus $b_i=0$ for all $i$ such that $X_1 \cap X_i \neq \emptyset$. Since $X$ is connected, repeating the 
argument we obtain $b_i=0$ for all $i$, as required.\\
(2) Suppose $\sD$ is a Cartier divisor which is numerically equivalent to $0$. 
Then $\sD_{\eta} \sim 0$, thus $\sD \sim \sum a_iX_i$, some $a_i$. Now a similar argument to the above shows
$\sD \sim 0$ as required.
\end{proof}

\begin{lem}
\label{lem-MV}
Let $X$ be a slc proper surface with $-K_X$ ample, and $\cX/T$ a 
smoothing with $\cX_{\eta} \cong \bP^2_{\eta}$.\\ 
(1) $N_1(X) \map N_1(\cX)$ and $N^1(\cX) \map N^1(X)$ are isomorphisms.\\
(2) $N_1(\Delta) \map N_1(X^{\nu}) \map N_1(X) \map 0$ is exact.\\
Dually, $N^1(\Delta) \leftarrow N^1(X^{\nu}) \leftarrow N^1(X) \leftarrow 0$ is exact.\\
Here $N_1(\Delta) \map N_1(X^{\nu})$ is defined as follows: if $\Gamma$ is a component of $\Delta$ obtained by identifying
two components $\Gamma_1 \subset X^{\nu}_j$ and $\Gamma_2 \subset X^{\nu}_k$ of $\Delta^{\nu}$,
we have $[\Gamma] \mapsto (0,\ldots,0,[\Gamma_1],0,\ldots,0,-[\Gamma_2],0,\ldots,0)$, 
where we're using the decomposition $N_1(X^{\nu})=\oplus_i N_1(X^{\nu}_i)$, 
and the nonzero entries are in the $j$th and $k$th positions. If $\Gamma$ is a component of $\Delta$ obtained by folding 
a component of $\Delta^{\nu}$, then $[\Gamma] \mapsto 0$.
\end{lem}

\begin{proof}
(1) First we show that $\Pic \cX \map \Pic X$ is an isomorphism. 
From the exponential sequence for $X$ we obtain the long exact cohomology sequence
$$ \cdots \map H^1(\cO_X) \map \Pic X \map H^2(X,\bZ) \map H^2(\cO_X) \map \cdots$$
We have $h^2(\cO_X)=h^0(K_X)$ by Serre duality.
Thus $h^2(\cO_X)=0$ since $-K_X$ is ample. Now $\chi(\cO_X)=\chi(\cO_{\bP^2})=1$ gives $h^1(\cO_X)=0$. So we obtain that
$\Pic X \map H^2(X,\bZ)$ is an isomorphism. Now, since $h^1(\cO_X)=h^2(\cO_X)=0$, cohomology and base change 
(\cite{Har}, p.~290, Theorem~12.11) gives $R^1f_{\star} \cO_{\cX} = R^2f_{\star} \cO_{\cX} =0$, 
thus $h^1(\cO_{\cX})=h^2(\cO_{\cX})=0$
since $T$ is affine. Using the exponential sequence for $\cX$ we obtain that $\Pic \cX \map H^2(\cX,\bZ)$ is an isomorphism.
Finally, $H^2(\cX, \bZ) \map H^2(X,\bZ)$ is an isomorphism since $X$ is a homotopy retract of $\cX$, so $\Pic \cX \map \Pic X$ is 
an isomorphism as claimed.

Now, since $\Pic \cX \map \Pic X$ is surjective, $N^1(\cX) \map N^1(X)$ is surjective and $N_1(X) \map N_1(\cX)$ is injective.
But $N_1(X) \map N_1(\cX)$ is clearly surjective since $X$ is the only closed fibre of $\cX / T$. Hence $N_1(X) \map N_1(\cX)$
and $N^1(\cX) \map N^1(X)$ are isomorphisms.

(2) First, we show that we have isomorphisms $N^1(X^{\nu}) \cong H^2(X^{\nu},\bQ)$, $N^1(X) \cong H^2(X,\bQ)$ and dually
$N_1(X^{\nu}) \cong H_2(X^{\nu},\bQ)$, $N_1(X) \cong H_2(X,\bQ)$. Well, first from the proof of (1) above we have
$\Pic \cX \otimes_{\bZ} \bQ \stackrel{\cong}{\map} H^2(\cX,\bQ)$, and by Lemma~\ref{lem-Cl}(2) we have 
$\Pic \cX \otimes_{\bZ} \bQ \stackrel{\cong}{\map} N^1(\cX)$, so we obtain $N^1(\cX) \stackrel{\cong}{\map} H^2(\cX,\bQ)$.
This gives $N^1(X) \stackrel{\cong}{\map} H^2(X,\bQ)$. For $X^{\nu}$, let $X_i^{\nu}$ be a component, 
let $\alpha_i : \tilde{X}_i \map X^{\nu}_i$ be a resolution. 
We may assume that $X$ is not normal, otherwise we are done by the above.
Then $X_i$  is rational by Theorem~\ref{thm-min-res}. 
Thus $\Pic(\tilde{X}_i) \otimes_{\bZ} \bQ \stackrel{\cong}{\map} N^1(\tilde{X}_i)$,
 it follows that $\Pic(X^{\nu}_i) \otimes_{\bZ} \bQ \stackrel{\cong}{\map} N^1(X^{\nu}_i)$.
We have $H^1(\cO_{X^{\nu}_i})=H^2(\cO_{X^{\nu}_i})=0$.
For $h^2(\cO_{X^{\nu}_i})=h^0(K_{X^{\nu}_i})=0$ since $-(K_{X^{\nu}_i}+\Delta^{\nu}_i)$ is ample, and
$H^1(\cO_{X^{\nu}_i})=0$ since $\tilde{X_i}$ is rational (using Leray).
So $\Pic(X^{\nu}_i) \otimes_{\bZ} \bQ \stackrel{\cong}{\map} H^2(X^{\nu}_i,\bQ)$ using the exponential sequence.
Thus $N^1(X^{\nu}) \stackrel{\cong}{\map} H^2(X^{\nu},\bQ)$ as desired.

It remains to show that
$$ H_2(\Delta) \map H_2(X^{\nu}) \map H_2(X) \map 0$$ 
is exact, where the map $H_2(\Delta) \map H_2(X^{\nu})$ is as in the statement of the lemma, and we work with $\bQ$ coefficients.
We use the Mayer-Vietoris sequence inductively here, separating off one component at a time.

First we separate the double curves $\Gamma \subset \Delta$ where the two branches of $X$ at $\Gamma$ belong to the same 
component. In this case, $X$ is homotopy equivalent to $X' \cup L$, where $X'$ is $X$ with the two branches at $\Gamma$ separated,
$L$ is $\Spec_{\bP^1}(\cO_{\bP^1} \oplus \cO_{\bP^1}(-1))$, and $X'$ and $L$ are glued along $\Gamma^{\nu} \subset X'$,
where $\Gamma^{\nu} \inj L$ is obtained using the 2-to-1 map $\Gamma^{\nu} \map \Gamma$. Note that $\Gamma^{\nu} \map \Gamma$ is 
either (a) $\bP^1 \map \bP^1$, 2-to-1,  or (b) $\bP^1 \cup \bP^1 \map \bP^1$ 
where the two components are joined at a node and each maps 
isomorphically onto $\Gamma$. Then
$$H_2(\Gamma^{\nu}) \map H_2(X') \oplus H_2(L) \map H_2(X) \map 0$$  
is exact, using $H_1(\Gamma^{\nu})=0$. In case (a) we have
$$
\begin{array}{ccccccc}
\bQ  & \map & H_2(X') \oplus  \bQ & \map & H_2(X) & \map & 0 \\
 1   & \mapsto & ([\Gamma^{\nu}] , -2) &   &        &      &   
\end{array}
$$
thus $H_2(X') \map H_2(X)$ is an isomorphism.
In case (b), writing $\Gamma^{\nu}=\Gamma^{\nu}_1 \cup \Gamma^{\nu}_2$, we have
$$
\begin{array}{ccccccc}
\bQ^{ \oplus 2}  & \map & H_2(X') \oplus  \bQ & \map & H_2(X) & \map & 0 \\
 (1,0)   & \mapsto & ([\Gamma^{\nu}_1] , -1) &   &        &      & \\
 (0,1)   & \mapsto & ([\Gamma^{\nu}_2] , -1) &   &        &      & 
\end{array}
$$
thus 
$$ 
\begin{array}{ccccccc}
\bQ & \map & H_2(X') & \map & H_2(X) & \map & 0 \\
1  & \mapsto & ([\Gamma^{\nu}_1] - [\Gamma^{\nu}_2]) & & & & 
\end{array}
$$
is exact. Next we separate the components of $X$. Let $X_1$ be a component (assumed normal), and write $X = X_1 \cup X'$, 
$\Delta_1=X_1 \cap X'$. Then
$$H_2(\Delta_1) \map H_2(X_1) \oplus H_2(X') \map H_2(X) \map 0$$
is exact, using $H_1(\Delta_1)=0$ (recall $\Delta_1$ is either $\bP^1$ or $\bP^1 \cup \bP^1$).
Now, using these steps repeatedly, we obtain our result.

\end{proof}

\begin{proof}[Proof of Theorem~\ref{thm-rho}]
(1) We have 
$$N_1(\Delta) \map N_1(X^{\nu}) \map N_1(X) \map 0$$ 
exact, where $[\Gamma] \in N_1(\Delta)$ maps to zero if
$\Gamma \subset \Delta$ is obtained by folding $\Gamma^{\nu} \subset \Delta^{\nu}$.  
Now $N_1(X) \cong N^1(X)^{\vee}$, 
$$N^1(X) \cong N^1(\cX) \cong \Pic(\cX) \otimes_{\bZ} \bQ \inj \Cl(\cX) \otimes_{\bZ} \bQ,$$
and $\dimn \Cl(\cX) \otimes_{\bZ} \bQ = V$.
Thus $\dimn N_1(X) \le V$, with equality iff $\cX$ is $\bQ$-factorial. So, using the exact sequence above, we obtain
$\sum \rho(X_i^{\nu}) \le V+E$, equal only if $\cX$ is $\bQ$-factorial.\\
(2) The maps $N_1(X^{\nu}_1) \map N_1(\cX)$ and  $N^1(\cX) \map N^1(X^{\nu}_1)$ are dual.
Let $X_2$ be the component of $X - X_1$ meeting $X_1$, let $\Gamma$ denote the intersection.
By the exact sequence of Lemma~\ref{lem-MV}, (2), we find that $N_1(X^{\nu}_1) \map N_1(\cX)$ is injective
(recall that $N_1(X) \stackrel{\cong}{\map} N_1(\cX)$).
For, if $\beta$ is in the kernel,   
$\beta = \lambda [\Gamma] \in N_1(X^{\nu}_1)$
for some $\lambda \in \bQ$. Intersecting with a relatively ample divisor on $\cX/T$, we find that $\lambda=0$, so $\beta=0$. 
Now using the description of $\Cl(\cX)$ in Lemma~\ref{lem-Cl} we see that $\dimn(\im (N^1(\cX) \map N^1(X^{\nu}_1))) \le 2$. For 
locally at $X_1$, $\Cl(\cX)$ is generated by $X_2$ and $\sH$ 
(where $\sH$ is a divisor flat over $T$, of degree $1$ on $\cX_{\eta} \cong \bP^2_{\eta}$).
So we obtain $\rho(X^{\nu}_1) \le 2$ as desired.\\ 
(3) Similarly, in this case we find that $\dimn \kernel [N_1(X^{\nu}_1) \map N_1(\cX)] \le 1$, and 
$\dimn \im [N^1(\cX) \map N^1(X_1)] \le 3$, so we obtain $\rho(X^{\nu}_1) \le 4$.

\end{proof}

\begin{proof}[Proof of corollary \ref{cor-rho}]
All this follows from Theorem~\ref{thm-rho} except for the claim that
$\rho(X_i) \le 3$ for $X_i$ a middle component of a surface of type $D$.
To prove this, observe that in this case $N_1(X_i) \map N_1(\cX)$ is injective. For write $X=X_1 \cup \cdots \cup X_k$, where
$X_i$ is glued to $X_{i+1}$ by glueing $\Delta^2_i$ to $\Delta^1_{i+1}$.
Suppose $\beta$ is in the kernel of $N_1(X_i) \map N_1(\cX)$.
By the exact sequence of Lemma~\ref{lem-MV}, (2), we have
$$[(0,\ldots,\beta,0,\ldots,0)] =
(\lambda_1 [\Delta^1_1], -\lambda_1 [\Delta^1_2] + \lambda_2 [\Delta^2_2], \ldots,
  -\lambda_{k-1} [\Delta^1_k]) \in \oplus N_1(X^{\nu}_j),$$ 
for some $\lambda_1, \ldots, \lambda_{k-1} \in \bQ$. Thus $\lambda_1=\lambda_{k-1}=0$. Working inductively, we obtain 
$\lambda_1= \lambda_2 = \cdots = \lambda_{i-1}=0$
and $\lambda_{k-1}=\lambda_{k-2}=\cdots=\lambda_{i}=0$, so $\beta=0$ as claimed.
We now conclude as in the proof of Theorem~\ref{thm-rho}, (3).
\end{proof}

\section{Ruling out the type B* surfaces} \label{B*}

\begin{thm} \label{thm-B*}
A surface of type B* (see Theorem~\ref{thm-glueing}) does not admit a $\bQ$-Gorenstein smoothing to $\bP^2$.
\end{thm}

\begin{proof}
Suppose $X$ is a counter example, let $\cX/T$ be a $\bQ$-Gorenstein smoothing with $\cX_{\eta} \cong \bP^2_{\eta}$.
We have $\rho(X^{\nu})=1$ by Theorem~\ref{thm-rho}, (1).
We have $K_X^2 = K_{\bP^2_{\eta}}^2 = 9$ since $K_{\cX}$ is $\bQ$-Cartier.
Thus $(K_{X^{\nu}}+\Delta^{\nu})^2=9$, so $K_{X^{\nu}}^2 > 9$, using $-K_{X^{\nu}}$ ample and $\rho(X^{\nu})=1$.
So $K_{X^{\nu}}^2+ \rho(X^{\nu}) > 10$.

Let $\tilde{X} \map X^{\nu}$ be the minimal resolution of $X^{\nu}$. Then $\tilde{X}$ is rational by Theorem~\ref{thm-min-res}.
Applying Noether's formula we obtain $K_{\tilde{X}}^2 + \rho(\tilde{X}) = 10$. So the resolution $\tilde{X} \map X^{\nu}$ has 
strictly decreased $K^2 + \rho$. However, we calculate below that the only possible singularities on $X^{\nu}$ will increase
$K^2 + \rho$ when we take the minimal resolution, so we have a contradiction.

$(X^{\nu},\Delta^{\nu})$ has singularities of type $(\frac{1}{dn^2}(1,dna-1),0)$ and $(\frac{1}{r}(1,a),\Delta)$, with the latter
 cases occurring in pairs $\frac{1}{r}(1,a)$ and $\frac{1}{r}(1,-a)$, by Theorem~\ref{thm-smoothability}. Now, given a cyclic 
quotient singularity $\frac{1}{r}(1,a)$, let $\frac{r}{a} = [b_1,\ldots,b_k]$ be the expansion of $\frac{r}{a}$
as a Hirzebruch continued fraction.
The geometric interpretation of this is that the minimal resolution of the singularity has exceptional locus a chain of smooth
rational curves 
with self-intersections $-b_1,\ldots,-b_k$. Then on taking the minimal resolution of the singularity, the change in $K^2 + \rho$  
is given by 
$$ \beta = k+2 + \sum_{i=1}^{k} (2-b_i) - \frac{a+a'+2}{r}$$
where $a'$ denotes the inverse of $a$ modulo $r$ (\cite{Ma}, p. 111). 
We have $\beta = d-1 \ge 0$ in the case $\frac{1}{dn^2}(1,dna-1)$ (\cite{Ma}, p. 112).
Finally, we calcuate $\beta_1+\beta_2 = 4(1- \frac{1}{r}) \ge 0$ in the case of a pair of singularities 
$\frac{1}{r}(1,a),\frac{1}{r}(-1,a)$. To see this, note that if we write $\frac{r}{a}=[b_1,\ldots,b_k]$ and
$\frac{r}{r-a}=[c_1,\ldots,c_l]$, we have $\sum (b_i-1) = \sum (c_j-1) = k+l-1$.
\end{proof}

\section{Classification of the normal surfaces} \label{normal}

Manetti has classified normal log terminal degenerations of $\bP^2$ in \cite{Ma} --- we will refer
to such surfaces as Manetti surfaces. We state the basic result below.

\begin{thm} (compare \cite{Ma}, p. 90, Main Theorem) \label{thm-Manetti}
Suppose $X$ is a normal log terminal proper surface with $-K_X$ ample. 
Then $X$ admits a $\bQ$-Gorenstein smoothing to $\bP^2$ iff
\begin{enumerate}
\item $X$ has singularities of type $\frac{1}{n^2}(1,na-1)$, $(a,n)=1$.
\item $\rho(X)=1$.
\end{enumerate}
\end{thm}

We now classify the normal log canonical degenerations of $\bP^2$.

\begin{thm}
\label{thm-normal-case}
Let $X$ be a normal log canonical proper surface with $-K_X$ ample.
Then $X$ admits a $\bQ$-Gorenstein smoothing to $\bP^2$ iff $X$ is a Manetti surface or $X$ is an elliptic cone of degree $9$.
\end{thm}

\begin{rem}
If $\cX/T$ is a smoothing of a normal proper surface $X$ to $\bP^2$, 
then $\cX/T$ is projective and $\bQ$-Gorenstein, and $-K_{\cX/T}$
is relatively ample. The projectivity is proved in \cite{Ma}, p. 95, Theorem~4,
the rest follows since $\rho(\cX/T)=1$ (by Lemma~\ref{lem-Cl}).
\end{rem}

\begin{lem}\label{lem-T_1}
Let $X$ be an slc proper surface.
Suppose $X$ has a $\bQ$-Gorenstein smoothing to $\bP^2$.
Then every normal log terminal singularity of $X$ is a cyclic quotient singularity of type $\frac{1}{n^2}(1,na-1)$, where 
$(a,n)=1$. Moreover $3 \nd n$. 
\end{lem}

\begin{proof}
First, we know that every normal log terminal singularity of $X$ is of the form $\frac{1}{dn^2}(1,na-1)$, since we assume there
exists a $\bQ$-Gorenstein smoothing (compare Theorem~\ref{thm-smoothability}(1)).
We need to show that $d=1$. We sketch the proof here, for details see \cite{Ma}, p. 103, Propn~13(i) and Remark~6. 
We compute that the Milnor fibre $F$ of a smoothing of a singularity of type 
$\frac{1}{dn^2}(1,na-1)$ has $b_2(F)=d-1$ and negative definite intersection product. 
Now, since $\bP^2$ has positive definite intersection product, it follows that $b_2(F)=0$, so $d=1$ as required.

Finally, we show $3 \nd n$. Let $\cX/T$ be a $\bQ$-Gorenstein smoothing of $X$ to $\bP^2$.
Let $(P \in X) \cong \frac{1}{n^2}(1,na-1)$. Then locally at $P \in \cX$, we have $\Cl(\cX) \cong \bZ / n\bZ$,
generated by $K_{\cX}$ (see \cite{KSB}, p. 313, Propn~3.10, and 
\cite{Ko2}, p. 135, Propn~2.2.7).
But $K_{\cX} \sim -3 \sH$ locally at $P$, where $\sH$ is a divisor, flat over $T$, that restricts to a hyperplane section
on $\cX_{\eta} \cong \bP^2_{\eta}$. Thus $3 \nd n$. For another proof using the Milnor fibre, see \cite{Ma},
p. 105, Theorem~15(ii).  
\end{proof}

\begin{notn}
We call singularities of the form $\frac{1}{dn^2}(1,dna-1)$ singularities of class $T$.
We call singularities of the form $\frac{1}{n^2}(1,na-1)$ singularities of class $T_1$.
\end{notn}

\begin{proof}[Proof of Theorem~\ref{thm-normal-case}]
Suppose given a normal log canonical del Pezzo surface $X$ which admits a $\bQ$-Gorenstein smoothing to $\bP^2$.
We may assume $X$ is strictly log canonical, otherwise $X$ is a Manetti surface.
Let $\pi : \tilde{X} \map X$ be the minimal resolution of $X$. If $\tilde{X}$ is not rational, 
then  $X$ is an elliptic cone by Theorem~\ref{thm-min-res}.
Then $K_X^2=K_{\bP^2}^2=9$ gives that $X$ is an elliptic cone of degree $9$.

So we may assume $\tilde{X}$ is rational. The Leray spectral sequence  gives an exact sequence
$$0 \map H^1(\cO_X) \map H^1(\cO_{\tilde{X}}) \map H^0(R^1f_{\star}\cO_{\tilde{X}}) \map H^2(\cO_X).$$
Now $h^2(\cO_X)=h^0(K_X)=0$ since $K_X$ is ample, and $h^1(\cO_{\tilde{X}})=0$ since $\tilde{X}$ is rational.
So $H^0(R^1f_{\star}\cO_{\tilde{X}})=0$, $X$ has rational singularities.

We can now use a result of Manetti (\cite{Ma}, p. 95, Theorem~4, and p. 100, Theorem~11):
Let $\phi : \tilde{X} \map \bar{X} \cong \bF_w$ be a birational morphism, with $w$ maximal
(so $\phi$ is an isomorphism over the negative section $B$ of $\bF_w$).
Let $p \colon \tilde{X} \map \bP^1$ denote the birational ruling so obtained.
Then the exceptional locus of $\pi$ is the strict transform $B'$ of $B$ together with the 
irreducible components of the degenerate fibres of $p$ of self-intersection $\le -2$.
In particular, $w \ge 2$.
Moreover every degenerate fibre contains a unique $-1$ curve. We quickly sketch the proof of this.
First, since $X$ smoothes to $\bP^2$, we have $h^0(-K_X) \ge h^0(-K_{\bP^2})=10$, and $h^0(-K_{\tilde{X}})=h^0(-K_X)$
since $\pi : \tilde{X} \map X$ is minimal. Manetti deduces there is no horizontal curve $C$ 
on $\tilde{X}$ with $C^2 \le -2$ except possibly $B'$. Now since $\rho(X)=1$, it follows that the exceptional locus of $\pi$ is $B'$
together with all the components of the degenerate fibres of self intersection $\le -2$, 
and every degenerate fibre has a unique $-1$ curve.  

There are two types of rational strictly log canonical surface singularities --- namely a $\mu_2$ quotient of a cusp
and a quotient of a simple elliptic singularity. Consider the minimal resolutions of these singularities.
In each case the exceptional locus is a union of smooth rational curves.
For a $\mu_2$ quotient of a cusp, the exceptional locus consists of a chain of curves with two $-2$ curves off each 
end component of the chain.
For a quotient of a simple elliptic singularity, the exceptional locus consists of a curve  with three
chains of curves off it. 
We now analyse how these could possibly fit into the minimal resolution $\tilde{X}$ of $X$.

Consider the minimal model program yielding $\phi : \tilde{X} \map \bar{X} \cong \bF_w$ in the neighbourhood of a given fibre
of $p: \tilde{X} \map \bP^1$.
At each stage we contract a $-1$ curve, meeting at most $2$ components of the fibre, and disjoint from $B'$ 
(since $\phi$ is an isomorphism over $B$).
We know that 
$$\Exc(\pi)=B' \cup \{ \Gamma  \subset \tilde{Y} \mid p_{\star}(\Gamma)=0 \mbox{ and } \Gamma^2  \le -2 \}.$$
This set decomposes  into the exceptional loci of the 
minimal resolutions of one log canonical rational singularity and some $T_1$ singularities.

First concentrate on the log canonical singularity; let $E$ denote the exceptional locus of its minimal resolution.
Then $E$ contains a curve $C$ which meets $3$ other components of $E$ --- we call such a curve a \emph{fork} of $E$.
Suppose $f$ is a degenerate fibre of $p$ containing a fork $C$ of $E$. Let $A'$ denote the component meeting $B'$ 
(this is the strict transform of the corresponding fibre $A$ of $\bar{X}=\bF_w$).
Then we have a decomposition $f=P \cup \Gamma  \cup Q \cup C \cup R \cup S$, where 
\begin{enumerate}
\item $\Gamma$ is the unique -1 curve in $f$.\\
\item $P$ is a string of curves, with one end component meeting $\Gamma$, $P$ contracts to a $T_1$ singularity (or is empty).\\
\item $Q$,$R$ and $S$ are nonempty configurations of curves meeting $\Gamma$ and $C$, $C$, and $C$ and $B'$ respectively.
\end{enumerate}
Then in the MMP $\tilde{X} \map \cdots \map \bar{X}$, we contract $\Gamma$, $P \cup Q$, $C$ and $R \cup S \backslash A'$ 
in that order.
(Note: $Q$ and $R$ are nonempty because $C$ is a fork. One might think $S$ could be empty since $E$ contains $B'$, but in
that case $C=A'$ and must be contracted before we can contract $R$, a contradiction).

Next suppose that $f$ is a degenerate fibre of $p$ that does not contain a fork of $E$. Then we have a decomposition 
$f= P \cup \Gamma  \cup Q$, where 
\begin{enumerate}
\item $\Gamma$ is the unique -1 curve in $f$.\\
\item $P$ is a string of curves, with one end component meeting $\Gamma$, $P$ contracts to a $T_1$ singularity (or is empty).\\
\item $Q$ is a non-empty string of curves meeting $\Gamma$ and with one end component meeting $B'$.
\end{enumerate}

We now analyse these two cases for each of the two types of singularity. We call them fibre types I and II.
First suppose $X$ has a $\mu_2$ quotient of a cusp singularity.
So $E=F \cup G_1 \cup \cdots \cup G_4 $ where $F=F_1 \cup \cdots \cup F_k$ is a chain of $\bP^1$'s,
and  $G_1,G_2$ (respectively $G_3,G_4$) are $-2$-curves meeting $F_1$ (respectively $F_k$).  

Suppose $F_1$ is contained in a degenerate fibre $f$. Then as above we can write
$f=P \cup \Gamma  \cup Q \cup C \cup R \cup S$ where without loss of generality $Q=G_1$, $C=F_1$, $R=G_2$ and 
$S=F_2 \cup \cdots \cup F_l$ for some $l < k$. Note that $f$ cannot contain the other fork $F_k$ of $E$, since then
$F_k=A'$ contradicting the description above. We contract $P \cup \Gamma \cup Q$ first. We deduce that the curves in the string 
$P$ have self-intersections $-3,-2,\ldots,-2$. Thus $P$ contracts to a $\frac{1}{2r+1}(1,r)$ singularity, where $r$ is the 
length of the string. But this is never a $T_1$ singularity (since $(r+1,2r+1)=1$), a contradiction. So $P$ is empty.
We can now calculate that the curves in the string $S$ have self-intersections $-3,-2,\ldots,-2,-1$ if $l>2$.
Then $F_l^2=-1$, a contradiction. Hence $l=2$.

Next suppose $F_1$ is not contained in a degenerate fibre. Then $F_1$ is horizontal, hence $F_1=B'$. Then it follows that we 
have a fibre of type II with $Q=G_1$, a -2 curve. We deduce that $P$ is a single -2 curve. But then $P$ contracts to a 
$\frac{1}{2}(1,1)$ singularity, which is not $T_1$, a contradiction.

Combining, we deduce that $k=5$, and we have two fibres of the form $ \Gamma \cup G_1 \cup G_2 \cup F_1 \cup F_2$ as above,
and $F_3=B'$. There are no further degenerate fibres. It only remains to calculate $w$. We use $K_X^2=9$ to deduce $w=11$.

I claim that the surface $X$ constructed above does not admit a $\bQ$-Gorenstein smoothing.
Let $Y \map X$ be the index one cover of $X$ at the 
singular point. Then $Y$ has a cusp singularity and the exceptional locus of the minimal resolution $\tilde{Y} \map Y$
is a cycle of 
rational curves of self-intersections $-2,-2,-2,-11,-2,-2,-2,-11$. Suppose $X$ has a $\bQ$-Gorenstein smoothing, then, taking
the canonical cover of the smoothing at the singular point we obtain a smoothing of $Y$. Let $M$ denote the Milnor fibre of
the smoothing of $Y$. Consider the intersection product on $H^2(M,\bR)$, write $b_2(M)=\mu_0+\mu_+ +\mu_-$, where 
$\mu_0$, $\mu_+$ and $\mu_-$ are the number of zero, positive and negative eigenvalues of the intersection form.
Since $Y$ is normal and Gorenstein, we have \cite{St2}
$$\mu_- = 10 h^1(\cO_{\tilde{Y}})+K_{\tilde{Y}}^2+b_2(\tilde{Y})-b_1(\tilde{Y}).$$   
In our case we calculate $\mu_-=10-18+8-1=-1$, a contradiction. So $Y$ is not smoothable, hence $X$ does not have a 
$\bQ$-Gorenstein smoothing (not even locally).

Now suppose $X$ has a quotient of a simple elliptic singularity. 
So $E=F \cup G^1 \cup G^2 \cup G^3$ where $F=\bP^1$, $G^i=G^i_1 \cup \cdots \cup G^i_{k(i)}$ is a chain of smooth rational 
curves
and $G^i_1$ meets $F$, for $i=1,2$ and $3$. We first give a partial classification of these singularities.
We can contract the chains $G^i$ to obtain a partial resolution $\hat{X} \map X$.
Write $\hat{F}$ for the image of $F$ under $\tilde{X} \map \hat{X}$. Then the chains $G^i$ contract to singularities of type
$(\frac{1}{r}(1,a),\Delta)$ on $(\hat{X},\hat{F})$. Let $r_1,r_2,r_3$ be the indices of these singularities,
then $\sum \frac{1}{r_i} =1$ (because $X$ is assumed to be strictly log canonical --- the condition is equivalent to 
$K_{\hat{X}}+\hat{F} = \pi^{\star}K_X$). Thus $(r_1,r_2,r_3)=(2,3,6),(2,4,4)$ or $(3,3,3)$ after reordering.
In particular, we see that each chain $G^i$ is either a single $\bP^1$ of self-intersection $-r_i$, or
a chain of $r_i-1$ $\bP^1$'s of self-intersection $-2$.

We claim that the fork $F$ of $E$ cannot be contained in a fibre $f$. By the classification above, its enough to show that
$w \neq 2$, since this then forces $F=A'$, a contradiction.
Write $K_{\tilde{X}}+\tilde{C}= \pi^{\star}K_X$.
Let $\bar{C} = \phi_{\star} \tilde{C}$.
Then 
$$(K_{\bar{X}}+\bar{C})^2 > (K_{\tilde{X}}+\tilde{C})^2=K_X^2=9.$$
Since $\phi$ is an isomorphism over $B$, we have $(K_{\bar{X}}+\bar{C})B=(K_{\tilde{X}}+\tilde{C})B'=\pi^{\star}K_X B' =0$
because $B'$ is $\pi$-exceptional. So $K_{\bar{X}}+\bar{C} \sim \lambda (B+wA)$, writing $A$ for a fibre of 
$\bar{X} \cong \bF_w \map \bP^1$. Here $\lambda = -2 + \mu$ where $\mu$ is the multiplicity of $B'$ in $\tilde{C}$.
Now $0 \le \mu \le 1$ since $X$ is log canonical and $\pi$ is minimal. Hence 
$$9 < (K_{\bar{X}}+\bar{C})^2=\lambda^2 (B+wA)^2= \lambda^2 w \le 4w,$$
so $w>2$ as required.

Thus $F$ is horizontal, $F=B'$ and we have 3 degenerate fibres of type II.
In each case $Q$ is a single curve of self-intersection $-r_i$ or a string of $(r_i-1)$ $-2$-curves.
If the fibre $f$ is a string, 
we deduce that $P$ is a string of $(r_i-1)$ $-2$-curves or a single curve of self-intersection $-r_i$ respectively.
Now, since $P$ contracts to a $T_1$ singularity, we deduce $P$ is a single $-4$-curve and $r_i=-4$.
If $f$ is not a string, we find that $Q$ is a string of three  $-2$-curves, $\Gamma$ meets the middle component, and $P$ is empty,
hence again $r_i=4$. So $r_i=4$ for all $i$, contradicting the classification above.

It remains to show that an elliptic cone of degree $9$ admits a $\bQ$-Gorenstein smoothing to $\bP^2$. We prove this in 
Lemma~\ref{lem-elliptic_cone} below.
\end{proof}

\begin{lem} \label{lem-elliptic_cone}
Let $X$ be an elliptic cone of degree $9$. Then $X$ has a smoothing to $\bP^2$.
\end{lem}

\begin{proof}
Given an elliptic cone $X$ of degree $9$, write $\tilde{X} \map X$ for the minimal resolution of $X$.
Then $\tilde{X}$ is a ruled surface over an elliptic curve $E$, i.e., $\tilde{X}=\underline{\bP}_E(\cO_E \oplus \sL^{\vee})$,
where $\sL$ is a line bundle on $E$ of degree $9$. We claim that $X$ is determined up to isomorphism by its section $E$. 
For $\Aut E$ acts transitively on $E$, and given a line bundle $\sL$ of degree $9$ we have $\sL \sim 9P$ for some $P \in E$
(c.f. \cite{Har}, p. 337, Exercise 4.6(b)), thus $\Aut E$ acts transitively on the line bundles of degree $9$. Our claim follows.

Let $T$ be the spectrum of a DVR, and write $\cY=\bP^2_T$.
Let $E \inj Y=\bP^2$ be an elliptic curve in the special fibre.
Let $\tilde{\cY} \map \cY$ be the blowup of $\cY$ in $E$.
Then the special fibre $\tilde{Y}$ consists of the strict transform $Y'$ of $Y$ together with a ruled surface $F$ 
of degree $9$ over the elliptic curve $E$. We contract $Y'$  to obtain a family $\bar{\cY}/T$ which is a 
smoothing of an elliptic cone of degree $9$ over $E$ to $\bP^2$.
\end{proof}

\begin{rem}
Note that any smoothing of $X$ is $\bQ$-Gorenstein since $K_X$ is Cartier.
\end{rem}

\section{Push forward and base change and relative Weil divisors} \label{pfbc}

The aim of this section is to define the notion of a relative Weil divisor
and explain the conditions `$\omega_{\cX/S}^{[i]}$ and $\cO_{\cX}(\sD)^{[i]}$ commute with base change'
in the definition (\ref{defn-allowable}) of an allowable family of stable pairs.
See also Section~\ref{QG}. 
We first recall some of Koll\'{a}r's theory of push forward and base change for open immersions. 
These results are proved in  \cite{Ko1} (unpublished), I at least provide the statements here.

\begin{notn}
Let $f \colon \cX  \map S$ be a morphism of schemes.
Let $i \colon \cU \inj \cX$ be an open subscheme, $\cW=\cX-\cU$, and $\cZ \subset \cX-\cU$ a subscheme proper over $S$.
Let $\cF$ be a coherent sheaf on $\cU$ which is flat over $S$.
Given a morphism $g \colon S' \map S$, write $\cX^{g}=\cX \times_S S'$, $\cU^{g}=\cU \times_S S'$ etc., 
$g_{\cX} \colon {\cX}^g \map {\cX}$ for the induced morphism, $\cF^g=g_{\cX}^{\star}\cF$.
\end{notn}

\begin{defn}
We say that the push forward of $\cF$ commutes with $g \colon S' \map S$ if the natural map
$g_{\cX}^{\star}i_{\star}\cF \map i^{g}_{\star}\cF^g$ is an isomorphism in a neighbourhood of $\cZ$.
We say that the push forward of $\cF$ commutes with arbitrary base change if 
the push forward of $\cF$ commutes with any $g \colon S' \map S$.
\end{defn}

\begin{rem}
In our applications, $\cZ$ contains all closed points of $\cW$, so an isomorphism in a neighbourhood of $\cZ$ is a 
global isomorphism. More specifically, we are only interested in the following special case:
$\cF$ is an invertible sheaf,
$\cX/S$ is a family of CM reduced surfaces, $\cX_s -\cU_s$ is finite for every $s \in S$, 
and either $\cX/S$ is proper and $\cZ=\cX-\cU$ or $P \in \cX$
is local and $\cZ=P$. Then $i_{\star}\cF$ is a `divisorial sheaf' as defined below. 
We write `$i_{\star}\cF$ commutes with base change' 
to mean the push forward of $\cF$ commutes with arbitrary base change in this case (since the choice of $\cU$ is immaterial by 
Lemma~\ref{lem-S_2}(1)).
\end{rem}

\begin{lem} \label{lem-pfbc_ff}
Notation as above.
\begin{enumerate}
\item The push forward of $\cF$ commutes with any flat $S' \map S$
\item Let $h \colon T \map S$ be faithfully flat. Then the push forward of $\cF$ commutes with arbitrary base change 
iff the push forward of $\cF^h$ commutes with arbitrary base change.
\item Assume that  the push forward of $\cF$ commutes with arbitrary base change and let $h \colon T \map S$ be a morphism,
then the push forward of $\cF^h$ commutes with arbitrary base change.
\end{enumerate}
\end{lem}

\begin{lem} \label{lem-pfbc_closed_pts}
The following are equivalent:
\begin{enumerate}
\item 	The push forward of $\cF$ commutes with arbitrary base change.
\item 	The push forward of $\cF$ commutes with every $j \colon s \inj S$ where $s \in S$ is a closed point.
\item 	The natural morphism $j^{\star}_X i_{\star}\cF \map i^{j}_{\star}\cF^j$ is surjective for
	every $j \colon s \inj S$ where $s \in S$ is a closed point.
\end{enumerate}
\end{lem}

\begin{lem} \label{lem-pfbc_flat}
If the push forward of $\cF$ commutes with arbitrary base change then $i_{\star}\cF$ is flat over $S$.
\end{lem}

\begin{lem} \label{lem-pfbc_inv_lim}
Notation as above. Assume that $0 \in S$ is local, 
$\cF^0$ satisfies Serre's condition $S_2$, and that 
$i^0_{\star}(\cF^0)$ is coherent. Then the following are equivalent:
\begin{enumerate}
\item The push forward of $\cF$ commutes with arbitrary base change.
\item The push forward of $\cF$ commutes with $j \colon 0 \inj S$.
\item For every local Artinian subscheme $h \colon A \inj S$ the push forward of $\cF^h$ commutes with $j \colon 0 \inj A$.
\end{enumerate}
\end{lem}

\begin{thm} \label{thm-pfbc_stratification} 
Notation as above. Assume that $S$ is Noetherian, $\cX/S$ is projective, and that $\cF^s$ is $S_2$ and $i^s_{\star}\cF^s$
is coherent for all $s \in S$. Then there exists a locally closed stratification $\amalg S_i \map S$, such that if
$T$ is a reduced scheme and $h \colon T \map S$ is morphism, then the push forward of $\cF^h$ commutes with 
arbitrary base change iff $h$ factors through $\amalg S_i \map S$.
\end{thm}

\begin{rem}
The assumption that $\cX/S$  is projective is necessary for general $S$.
However, if we assume that $S$ is the spectrum of a complete local ring, then the conclusion holds for arbitrary $\cX/S$.
\end{rem}

\begin{rem} \label{rem-stratification}
We would like to remove the requirement that $T$ is a \emph{reduced} scheme
--- we do not know if this is possible.
\end{rem}

We now define the notion of a relative Weil divisor for a family of CM reduced surfaces $\cX/S$
over an arbitrary base $S \in \underline{Sch}$. This is a generalisation of Mumford's notion of a relative Cartier divisor
(\cite{Mu}, Lecture 10).

\begin{defn}
Let $\cX/S$ be a family of CM reduced surfaces.
We say a codimension $1$ closed subscheme $\sD$ of $\cX$ is a \emph{relative effective 
Weil divisor} if there exists an open subscheme $i \colon \cU \inj \cX$ and an effective Cartier divisor $\sD^0$ on $\cU/S$, 
flat over $S$, such that
\begin{enumerate}
\item $\cX_s - \cU_s$ is finite for each $s \in S$.
\item $\sD = \overline{\sD^0}$, the scheme theoretic closure of $\sD^0$ in $\cX$.
\end{enumerate}
We define a \emph{relative Weil divisor} to be a formal difference $\sD^{+} - \sD^{-}$ of relative effective Weil divisors. 
Given a relative Weil divisor $\sD$, we define an associated sheaf $\cO_{\cX}(\sD)=i_{\star}\cO_{\cU}(\sD |_{\cU})$,
where $i \colon \cU \inj \cX$ is the open subscheme where $\sD^{+}$ and $\sD^{-}$ are Cartier,
and $\cO_{\cU}(\sD |_{\cU})$ is the invertible sheaf corresponding to the Cartier divisor $\sD |_{\cU}$ as usual.
We say that $\sD$ is \emph{Cartier} if $\cO_{\cX}(\sD)$ is invertible.
Given $T \map S$, 
we define the pullback $\sD_{(T)}$ of a relative effective Weil divisor $\sD$ to $\cX_T=\cX \times_S T$ via
$\sD_{(T)}=\overline{\sD^0 \times_S T}$. We define the pullback of a general relative Weil divisor by linearity. 

We say a coherent sheaf $\cF$ on $\cX/S$ is a \emph{divisorial sheaf} if $\cF=\cO_{\cX}(\sD)$ 
for some relative Weil divisor $\sD$.
Equivalently, there exists  an open subscheme $i \colon \cU \inj \cX$ such that $\cF \mid_{\cU}$ is invertible,
$\cX_s - \cU_s$ is finite for each $s \in S$ and $\cF=i_{\star}i^{\star}\cF$.

Given a divisorial sheaf $\cF$ and $N \in \bN$, let $\cF^{[N]}$ denote the sheaf $i_{\star}((i^{\star}\cF)^{\otimes N})$ 
(corresponds to multiplication of the  divisor by $N$). 
\end{defn}

\begin{rem} 
The assumption that $\sD^0$ is flat over $S$ is equivalent to the following: 
for all $s \in S$, $\Supp \sD$ does not contain any component of $\cX_s$ 
(\cite{Mu}, Lecture 10, p. 72).
\end{rem}

\begin{rem} 
Given a relative Weil divisor $\sD$ on $\cX/S$, the sheaf $\cO_{\cX}(\sD)$ is coherent 
(using Remark~\ref{rem-rel_S_2} below).
\end{rem}

\begin{rem}
If $\cF$ is a divisorial sheaf and $j \colon \cU \inj \cX$ is \emph{any} 
open subscheme such that $\cX_s - \cU_s$ is finite for each $s$ then $j_{\star}j^{\star}\cF=\cF$ 
(using Lemma~\ref{lem-S_2}(1)).
\end{rem}

\begin{rem}
Note that $\sD_{(T)}$ is \emph{not} the same as $\sD_T=\sD \times_S T$ in general, and moreover $\sD$ is not necessarily flat over
$S$.
See Lemma~\ref{lem-rwd_ideal_pfbc} and Example~\ref{ex-rel_W_div} below. 
\end{rem} 

Our next result is a technical lemma which, given a family $\cX/S$ and a sheaf $\cF$ on $\cX$, flat over $S$,
relates the $S_2$ property for the fibres $\cF_s$ of $\cF$ to a relative $S_2$-type property for $\cF$. 
Note that it is \emph{not} true that if every fibre $\cF_s$ is $S_2$ then   
the sheaf $\cF$ is $S_2$ --- we can easily construct a counter example where the base $S$ is not $S_2$.

\begin{lem} \label{lem-S_2}
Let $\cX/S$ be a family of CM reduced surfaces and $\cF$ a coherent sheaf on $\cX$ which is flat over $S$.
\begin{enumerate}
\item Suppose that for each closed point $s \in S$ the sheaf $\cF_s$ on $\cX_s$ satisfies Serre's condition $S_2$.
Let $i \colon \cU \inj \cX$ be an open subscheme such that the set $\cX_s - \cU_s$ is finite for each $s \in S$.
Then the natural map $\cF \map i_{\star}i^{\star}\cF$ is an isomorphism.
\item Suppose that for each closed point $s \in S$ the sheaf $\cF_s$ on $\cX_s$ is invertible in codimension $1$.
Then there exists $i \colon \cU \inj \cX$ such that the set $\cX_s - \cU_s$ is finite for each $s \in S$ and 
$i^{\star}\cF$ is invertible.
\end{enumerate}
In particular, if $\cF_s$ is invertible in codimension $1$ and 
$S_2$ for each closed point $s \in S$, then $\cF$ is a divisorial sheaf.
\end{lem}

\begin{proof}
(1) Write $\cZ = \cX - \cU$. We work locally at a closed point $P \in \cZ \subset \cX$, say $P \mapsto s \in S$.
Then $\cZ_s \subset \cX_s$ is a closed subscheme with support $P$.
The sheaf $\cF_s$ is $S_2$ by assumption, so there exists a regular sequence $x_s,y_s \in m_{\cX_s,P}$ for $\cF_s$
at $P$. Replacing $x_s,y_s$ by $x_s^k,y_s^k$ if necessary, we may assume that $x_s,y_s \in \cI_{\cZ_s}$.
Now lift $x_s,y_s$ to $x,y \in \cI_{\cZ}$, then $x,y$ is a regular sequence for $\cF$ at $P$
 (\cite{Mat}, p. 177, Corollary of Theorem~22.5).
Equivalently, we have an exact sequence
$$0 \map \cF \stackrel{(y,-x)}{\map} \cF \oplus \cF \stackrel{(x,y)}{\map} \cF.$$
Consider the natural map  $\cF \map i_{\star}i^{\star} \cF$, write $K$ for the kernel and $C$ for the cokernel.
$K$ and $C$ have support contained in the set $\cZ$, so any given element of $K$ or $C$ is annihilated by some power of $\cI_Z$.
So, if $K \neq 0$, there exists $0 \neq g \in K$ such that $\cI_{\cZ} g = 0$, so in particular $xg=yg=0$, 
contradicting the exact sequence above. Similiarly if $C \neq 0$, there exists $g \in i_{\star}i^{\star}\cF - \cF$ such that
$\cI_{\cZ}g \subset \cF$. Again using the exact sequence above, since $(yg,-xg) \mapsto 0$ we obtain $(yg,-xg)=(yg',-xg')$
for some $g' \in \cF$, it follows that $g=g'$, a contradiction. Thus $K=C=0$, so the map $\cF \map i_{\star}i^{\star} \cF$
is an isomorphism as claimed.

(2) If $\cF_s$ is invertible at $P \in \cX_s$ then, working locally at $P$, lifting an isomorphism $\cO_{\cX_s} \map \cF_s$
we obtain a surjection $\cO_{\cX} \map \cF$, by Nakayama's Lemma. Now, by flatness of $\cF$, it follows that $\cO_{\cX} \map \cF$
is an isomorphism, so $\cF$ is invertible at $P \in \cX$.
\end{proof}

\begin{rem} \label{rem-rel_S_2}
In particular, if $\cX/S$ is a family of CM reduced surfaces, and $i \colon \cU \inj \cX$ is an open inclusion
such that $\cX_s-\cU_s$ is finite for each $s \in S$, then $i_{\star}i^{\star}\cO_{\cX}=\cO_{\cX}$.
\end{rem}

\begin{ex}
The sheaf $\omega_{\cX/S}$ is divisorial for a family $\cX/S$ of slc surfaces.
For $\omega_{\cX/S}$ is flat over $S$, and, for each $s \in S$, 
the natural map $\omega_{\cX/S}\otimes k(s) \map \omega_{\cX_s}$ is an 
isomorphism and $\omega_{\cX_s}$ is invertible in codimension $1$ and $S_2$. 
We write $K_{\cX/S}$ for a relative Weil divisor such that $\omega_{\cX/S}=\cO_{\cX}(K_{\cX/S})$.
\end{ex}

\begin{lem} \label{lem-rwd_ideal_pfbc}
Let $\cX/S$ be a family of CM reduced surfaces,  $\sD \subset \cX$ a relative effective Weil divisor.
Then $\sD_{(T)}=\sD_T$ for all $T \map S$ iff $\cO_{\cX}(-\sD)$ commutes with base change.
Moreover, in this case $\sD$ is flat over $S$.
\end{lem}
\begin{proof}
We may assume that $S$ and $T$ are affine, write $S=\Spec A$, $T=\Spec B$. We have a commutative diagram:
\begin{eqnarray*}
\begin{array}{cccccccccc}
0 & \map & \cTor_1^{A}(\cO_{\sD},B) & \map & \cO_{\cX}(-\sD) \otimes_A B & \map & \cO_{\cX_B} & \map  & \cO_{\sD_B} & \map 0\\
  &      &                             &      & \da                       &      &  \da =      &       & \da         & \\
  &      &  0                          & \map & \cO_{\cX_B}(-\sD_{(B)})        & \map & \cO_{\cX_B} & \map  & \cO_{\sD_{(B)}} & \map 0
\end{array}
\end{eqnarray*}
Thus $\sD_B =\sD_{(B)}$ iff $\cO_{\cX}(-\sD) \otimes_A B \map \cO_{\cX_B}(-\sD_{(B)})$ is surjective, by the snake lemma.
This proves the first part using Lemma~\ref{lem-pfbc_closed_pts}(3). 
Moreover, we see that in this case $\cTor_1^{A}(\cO_{\sD},B)=0$
for all $A \map B$, hence $\sD$ is flat over $S$ as required.
\end{proof}

\begin{ex} \label{ex-rel_W_div}
Let $\cX/T$ be a family of surfaces over the spectrum of a DVR $T$ with generic fibre $\bF_0/ k(\eta)$ and special fibre $\bF_4 /k$
(this can be realised as a family of scrolls in a projective space). We can contract the negative section $B$ of $\bF_4$ 
to obtain a family $\bar{\cX}/T$. Let $P \in \bar{\cX}$ be the image of $B$, then $K_{\bar{\cX}}$ is not $\bQ$-Cartier at $P$.
However the special fibre $\bar{X}=\bP(1,1,4)$ has $2K_{\bar{X}}$ Cartier.
Let $\sD \in |-2K_{\bar{\cX}}|$ be a relative Weil divisor,
then necessarily $P \in \Supp \sD$, let $D$ be the restriction of the divisor 
$\sD$ to the special fibre. Then the natural map $\cO_{\cX}(-\sD) \otimes k \map \cO_X(-D)$ is not surjective. 
For otherwise $\cO_{\cX}(-\sD)$ is invertible by Nakayama's Lemma, a contradiction. Thus the scheme theoretic fibre 
$\sD \otimes k$ has an embedded point at $P$, in particular $\sD \otimes k \neq D$.
\end{ex}

\begin{ex} \label{ex-rwd}
Let $\cX/S$ be a family of CM reduced surfaces and $\sD \subset \cX/S$ a codimension 1 closed subscheme, flat over $S$.
Then for $s \in S$ we can define the restriction $\sD_{(s)}$ of $\sD$ to the fibre $\cX_s$
(we take the double dual of the ideal sheaf of the fibre $\sD_s$ as the ideal sheaf of $\sD_{(s)}$).
Consider the locus $S'$ of $s \in S$ such that $\sD_s=\sD_{(s)}$ (equivalently $\cO_{\cX}(-\sD) \otimes k(s)$ is $S_2$). 
Since $\cO_{\cX}(-\sD)$ is flat over $S$, this locus is open.
Finally, let $S''$ be the open locus of points $s \in S'$ such that $\sD_{(s)}$ is Cartier in codimension $1$. 
Then $\sD''=\sD \times_S {S''}$ is a relative Weil divisor on $\cX''=\cX \times_S S''/ S''$ by Lemma~\ref{lem-S_2}.
\end{ex}

We prove below some foundational results concerning $\bQ$-Cartier relative Weil divisors.

\begin{lem}
\label{lem-pfbc}
Let $\cX/S$ be a family of CM reduced surfaces, with $0 \in S$ local.
Let $\sD$ be a relative Weil divisor. 
Suppose that $ND$ is Cartier.
Then $\cO_{\cX}(N\sD)$ commutes with base change if and only if $N\sD$ is Cartier.
\end{lem}
\begin{proof}
If $\cO_{\cX}(N\sD)$ commutes with base change then $\cO_{\cX}(N\sD) \otimes k \map \cO_X(ND)$ is an isomorphism,
and $\cO_X(ND)$ is invertible by assumption. Since $\cO_{\cX}(N\sD)$ is coherent, it follows that $\cO_{\cX}(N\sD)$ is invertible
by Nakayama's Lemma and the flatness of $\cO_{\cX}(N\sD)$ (using Lemma~\ref{lem-pfbc_flat}). 

Conversely, suppose $N\sD$ is Cartier. By Lemma~\ref{lem-pfbc_closed_pts}(2), 
we need only show that $\cO_{\cX}(N\sD) \otimes k \map \cO_X(ND)$ is an isomorphism.
Both sides are invertible by assumption and the map is an isomorphism in codimension $1$, 
thus it is an isomorphism since $X$ is $S_2$.
\end{proof}

\begin{prop-defn}
\label{pd-index_one_cover}
Let $\cX/S$ be a family of CM reduced surfaces with $P \in \cX$ local,
and $\sD$ a $\bQ$-Cartier relative Weil divisor on $\cX/S$ of index $N$. 
We define a $\mu_N$ quotient $\pi \colon \cZ \map \cX$ as follows:
$$\cZ=\underline{\Spec}_{\cX}(\cO_{X} \oplus \cO_{\cX}(\sD) \oplus \cdots \oplus \cO_{\cX}((N-1)\sD))$$
where the multiplication is given by fixing an isomorphism $\cO_{\cX}(N\sD) \cong \cO_{\cX}$.
Let $i \colon \cX^0 \inj \cX$ be the locus where $\sD$ is Cartier and $j \colon \cZ^0 \inj \cZ$ the
inverse image in $\cZ$. We have
\begin{enumerate}
\item $\pi^0 \colon \cZ^0 \map \cX^0$ is etale.
\item The relative Weil divisor $\sD_{\cZ}=\pi^{\star} \sD$ is Cartier.
\item $j_{\star}j^{\star}\cO_{\cZ}=\cO_{\cZ}$.
\end{enumerate}
Conversely, any $\mu_N$ quotient $\pi \colon \cZ \map \cX$ satisfying these criteria is of the form above.
We refer to such a cover as a cyclic cover of $\cX/S$ defined by $\sD$, or, in the case $\sD=K_{\cX/S}$,
an index one cover of $\cX/S$.

The construction is unique up to the choice of an element of $H^0(\cO_{\cX}^{\times}) /H^0(\cO_{\cX}^{\times})^N$.
In particular if $P \in \cX$ is a local analytic germ, the construction is unique. 
\end{prop-defn}

\begin{proof} 
This is a straightforward generalisation of the usual cyclic covering trick \cite{YPG}.
\end{proof}

\begin{rem}
Note that $\pi^{-1}P$ is a single point (using $\ind \sD =N$).
\end{rem}

\begin{lem} 
\label{lem-rel_index}
Let $\cX/S$ be a family of CM reduced surfaces, with $0 \in S$ local, and $\sD$ a relative Weil divisor on $\cX$.
Suppose $\sD$ is $\bQ$-Cartier. Then $\ind(\sD)=\ind(D)$.
\end{lem} 

\begin{proof}
Let $\ind\sD =N$, $\ind D = M$, then $N\sD$ Cartier gives $ND = N\sD \mid_X$ Cartier, so $M \mid N$.
We claim that $M =N$.
We may work locally at $P \in \cX$.

First suppose that $S=\Spec A$ where $A$ is a local Artinian $k$-algebra.
Let $\pi \colon \cZ \map \cX$ be a cyclic cover defined by $\sD$. 
Let $i \colon \cX^0 \inj \cX$ be the locus where $\sD$ is Cartier 
(then $\cX^0$ is the open subscheme of $\cX$ with underlying space
$X-P$) and let $j \colon \cZ^0 \inj \cZ$ be the inverse image in $\cZ$. We have $j_{\star}\cO_{\cZ^0}=\cO_{\cZ}$.
Thus $\cZ$ connected implies  $\cZ^0$ is connected.
Hence $Z^0$ is connected, since $Z^0$ and $\cZ^0$ have the same underlying space. 
But $$Z^0=\underline{\Spec}_{X^0}(\cO_{X^0} \oplus \cO_{X^0}(-D^0) \oplus \cdots \oplus \cO_{X^0}(-(N-1)D^0)),$$
thus the cover 
$$Z'=\underline{\Spec}_X(i_{\star}\pi^0_{\star}\cO_{Z^0}) = 
\underline{\Spec}_{X}(\cO_{X} \oplus \cO_{X}(-D) \oplus \cdots \oplus \cO_{X}(-(N-1)D))$$ 
of $X$ is connected, it follows that the index $M$ of $D$ is equal to $N$.

The general case now follows using Lemma~\ref{lem-pfbc} and Lemma~\ref{lem-pfbc_inv_lim} --- we have $M\sD_{(A)}$ Cartier
for every Artinian subscheme $A \inj S$, hence $M\sD$ is Cartier, so $M=N$.
\end{proof}

\section{A review of versal deformations according to Artin} \label{artin}
 
The aim of this section is to give necessary conditions for a groupoid $F$ over $\underline{C}$ to admit
algebraic, everywhere versal deformations (Theorem~\ref{thm-patches}). This theory was developed by Artin in \cite{Ar1}, 
starting from the work of Schlessinger \cite{Sch}. In Section~\ref{construct} we use our result to prove that $\cM_d$ is an 
algebraic stack.

\begin{notn}
Let $\underline{C}$ be the category of noetherian $\bC$-algebras. 
We say a morphism $A' \map A$ in $\underline{C}$ is an \emph{extension} if it is surjective,
we say it is an $\emph{infinitesimal extension}$ if it has nilpotent kernel.
\end{notn}

\begin{notn}
Let $F$ be a groupoid over $\underline{C}$.
We write $F(A)$ for the fibre over $A \in \underline{C}$.
Given $A \map B$ in $\underline{C}$ and $b \in F(B)$, we write $F_b(A)$ for the groupoid of maps $a \map b$ in $F$ 
lying over $A \map B$.
We write $\bar{F}(A)$, $\bar{F}_b(A)$ for the isomorphism classes of $F(A)$, $F_b(A)$, these define functors 
$$\bar{F} \colon \underline{C} \map (\emph{Sets}),$$ 
$$\bar{F}_b \colon \underline{C} \backslash B \map (\emph{Sets}).$$
\end{notn}

\begin{defn}
Let $F$ be a groupoid over $\underline{C}$.
We say that $F$ is \emph{limit preserving} if the natural functor 
$$\stackrel{\lim}{\map} F(A_i) \map F(\stackrel{\lim}{\map} A_i)$$
is an equivalence of categories for every direct system $\{ A_i \}$ in $\underline{C}$ such that 
$\stackrel{\lim}{\map} A_i \in \underline{C}$.
\end{defn}

\begin{notn} \label{notn-art}
Let $A_0 \in \underline{C}$, let $A \map A_0$ be an extension, and $A' \map A$ an infinitesimal extension with kernel $M$, where $M$
is a finite $A_0$-module (i.e., writing $N$ for the kernel of $A' \map A_0$, we have $MN=0$ in $A'$).
We write $A+M$ for the trivial extension of $A$ by $M$, namely the $A$-module $A \oplus M$ with multiplication given by $M^2=0$.
\end{notn}

We define conditions (S1)(a), (S1)(b) and (S2) for a groupoid $F$ over $\underline{C}$ as follows:

\begin{S1a} Let
\begin{eqnarray*}
\begin{array}{ccc}
   &      & B \\
   &      & \da \\
A' & \map & A
\end{array}
\end{eqnarray*}
be a diagram in $\underline{C}$, where $A' \map A$ is as in Notation~\ref{notn-art}.
Assume that the composed map $B \map A_0$ is surjective. Let $a \in F(A)$.
Then the canonical map 
$$\bar{F}_a(A' \times_A B) \map \bar{F}_a(A') \times \bar{F}_a(B)$$
is surjective.
\end{S1a}

\begin{S1b} Let $B \map A_0$ be surjective, and let $M$ be a finite $A_0$-module. Let $b \in F(B)$ have direct image 
$a_0 \in F(A_0)$. Then the canonical map 
$$\bar{F}_b(B+M) \map \bar{F}_{a_0}(A_0+M)$$
is bijective.
\end{S1b}

\begin{rem}
If $F$ satisfies (S1)(b), $\bar{F}_{a_0}(A_0+M)$ has a natural $A_0$-module structure, and the underlying additive group 
acts on $\bar{F}_a(A')$ \cite{Sch}.
\end{rem}

\begin{notn} 
Write $D_{a_0}(M)= \bar{F}_{a_0}(A_0+M)$.
\end{notn}

\begin{S2}
$D_{a_0}(M)$ is a finite $A_0$-module.
\end{S2}

\begin{defn}
By an obstruction theory $\cO$ for $F$ we mean the following data:
\begin{enumerate}
\item For each infinitesimal extension $A \map A_0$ and element $a \in F(A)$, a functor
$$\cO_a \colon ( \emph{ finite } A_0-\emph{modules } ) \map ( \emph{ finite } A_0-\emph{modules } ).$$
\item For each $A' \map A$ as in Notation~\ref{notn-art} and $a \in F(A)$, an element $o_a(A') \in \cO_a(M)$
which is zero iff $\bar{F}_a(A') \neq \emptyset$.
\end{enumerate}
\end{defn} 

\begin{rem}
The data is required to be functorial. That is, given a morphism of extensions
\begin{eqnarray*}
\begin{array}{ccccccccc}
0 & \map & M   & \map & A'  & \map & A      & \map & 0\\
  &      & \da \phi &      & \da &      & \da \theta   &      &  \\
0 & \map & N   & \map & B'  & \map & B      & \map & 0 
\end{array}
\end{eqnarray*}
where $A \map B \map A_0$ is surjective and $M$, $N$ are finite $A_0$-modules,
and given $a \in F(A)$ with direct image $b \in F(B)$, we have 
a natural map $\cO_a(M) \map \cO_b(N)$ determined by $\theta$ and $\phi$ such that
$o_a(A') \mapsto o_b(B')$.

We also require that the data is linear in $(A_0,M)$, i.e., given $M$ and $N$ finite $A_0$-modules, the map
$$\Hom_{A_0}(M,N) \map \Hom_{A_0}(\cO_{a}(M),\cO_{a}(N)),\mbox{ } \phi \mapsto \cO_a(\phi)$$
is an $A_0$-module homomorphism.
\end{rem}

\begin{defn}
We say $A \in \underline{C}$ is \emph{algebraic} if it is of finite type over $\bC$.
We say $a \in F(A)$ is \emph{algebraic} if $A$ is algebraic.
\end{defn}

We state some further conditions we require for $D$ and $\cO$.

\begin{cond} \label{cond-D_and_O}
\begin{enumerate}
\item $D$ and $\cO$ are compatible with etale localisation, i.e., given $p \colon A \map B$ etale, $a \in F(A)$ with direct image
$b \in F(B)$ etc., we have
$$D_{b_0}(M \otimes B_0) \cong D_{a_0}(M) \otimes B_0,$$
and
$$\cO_b(M \otimes B_0) \cong \cO_a(M) \otimes B_0.$$
\item $D$ is compatible with completions, i.e., for $m$ a maximal ideal of $A_0$, we have
$$D_{a_0}(M) \otimes \hat{A_0} \cong \stackrel{\lim}{\leftarrow} D_{a_0}(M / m^n M).$$
\item Constructibility: For $A_0$ reduced, there is an open dense set of points $p \in \Spec A_0$ 
such that 
$$D_{a_0}(M) \otimes k(p) \cong D_{a_0}(M \otimes k(p)),$$
and
$$\cO_a(M) \otimes k(p) \subseteq \cO_a(M \otimes k(p)).$$
\end{enumerate}
\end{cond}

\begin{defn}
The \emph{lifting property} for $v \in \Ob(F)$ is the following: Given $\theta \colon v \map a$ and $\phi \colon a' \map a$  in 
$F$ such that the direct image $A' \map A$ of $\phi$ in $\underline{C}$ is surjective, there exists 
$\theta' \colon v \map a'$  such that $\phi \circ \theta' = \theta$.  

An element $v \in F(R)$ is \emph{formally versal} at $p \in \Spec R$ if the lifting property holds whenever $A'$ is a finite 
length extension of the residue field $k(p)$. We say $v \in F(R)$ is \emph{formally smooth} over $F$ 
if the lifting property holds whenever $A' \map A$ is an infinitesimal extension.
\end{defn}

\begin{thm} \label{thm-patches}
Let $F$ be a limit preserving groupoid over $\underline{C}$.
Assume that we are given an obstruction theory $\cO$ for $F$.
Suppose that $F$ satisfies (S1)(a),(b) and (S2) and that $D$ and $\cO$ satisfy the conditions of \ref{cond-D_and_O} 
for algebraic $A$, $B \in \underline{C}$. Suppose also that, if $\hat{A}$ is a complete local ring in $\underline{C}$,
the map  
$$\bar{F}(\hat{A}) \map \stackrel{\lim}{\leftarrow} \bar{F}(\hat{A}/m^n)$$
has dense image. Then, given $u \in F(\bC)$ 
there exists an algebraic ring $R$ with a closed point $0 \in \Spec R$ and $v \in F(R)$, formally smooth over $F$,
such that $v_0=u$.
\end{thm}
 
\begin{proof}
This follows from \cite{Ar1}, p. 170, Corollary~3.2 and p. 175, Theorem~4.4.
\end{proof}

\section{$\bQ$-Gorenstein deformation theory} \label{QG}

We define a workable theory of $\bQ$-Gorenstein deformations of slc surfaces as suggested in \cite{Ko6} ---
we say $\cX/S$ is $\bQ$-Gorenstein if $\omega_{\cX/S}^{[i]}$ commutes with base change for all $i \in \bZ$ 
(Definition~\ref{def-QG}).
In particular, for $(\cX,\sD)/S \in \cM_d(S)$, the family $\cX/S$ is $\bQ$-Gorenstein by definition.  
This theory is used in Section~\ref{construct} to show that $\cM_d$ is an algebraic stack using the methods of Artin. 
We first review some standard deformation theory.

\begin{notn} \label{notn-scheme}
We consider families of schemes $\cX/A$, where $A \in \underline{C}$ and $\cX$ is a noetherian scheme, flat over $A$, 
which is either of finite type over $A$ or affine, and is separated for the $m_P$-adic topology for each $P \in \Spec A$. 
Given an infinitesimal extension $A' \map A$ in $\underline{C}$, 
write $\Def_{\cX/A}(A')$ for the set of deformations of $\cX/A$ over $A'$.
Given a family $\cX'/A'$ extending $\cX/A$,
write $\Aut_{\cX/A}(\cX'/A')$ for the group of automorphisms of $\cX'/A'$ over $A'$ which restrict to the identity on $\cX/A$.
Given $A \map B$ in $\underline{C}$, we also write $\Aut_{\cX/A}(B)$ for $\Aut_{\cX/A}(\cX \otimes_A B /B)$. 
\end{notn}

\begin{rem} 
We want to allow the fibres to be \emph{local}, in particular we cannot assume that $\cX/A$ is of 
finite type. We need to insist that $\cX/A$ is separated for the $m_P$-adic topology for each $P \in \Spec A$ to ensure 
that, e.g., there are no empty fibres of $\cX/A$ (c.f. \cite{Ko4}, p. 21).
\end{rem}

\begin{defn} Let $A \in \underline{C}$ and let $\cX/A$ be a family of schemes over $A$ as in Notation~\ref{notn-scheme}. 
Let $\sL_{\cdot}$ be a cotangent complex for $\cX/A$, in the derived category of coherent sheaves on $\cX$
(\cite{LS}, p.~44, Definition~2.1.3). 
Given a coherent sheaf $\cF$ on $\cX$, define $\cM(\cF)^{\cdot}=\cHom_{\cO_{\cX}}(\sL_{\cdot},\cF)$.
Define $\cT^i(\cX/A,\cF) =\sH^i(\cM(\cF)^{\cdot})$ and $T^i(\cX/A,\cF)=\bH^i(\cM(\cF)^{\cdot})$.
Here, given a complex of sheaves, we use $\sH$ to denote the cohomology sheaves and $\bH$ to denote the hypercohomology groups.
\end{defn}

\begin{rem} 
If $\cX/A$ is of finite type we may assume that the sheaves $\sL_j$ are coherent, so in particular
the $\cT^i(\cX/A,\cF)$ are coherent.
\end{rem}

\begin{rem} We only need the cases $\cF=\cO_{\cX} \otimes_A B$ for some finite $A$-module $B$.
\end{rem}

\begin{rem} There is a local-to-global spectral sequence relating the $\cT^i$ and the $T^i$ :
$$E^{pq}_2 =H^p(\cX,\cT^q(\cX/A,\cF)) \Rightarrow T^{p+q}(\cX/A,\cF).$$ 
This is the usual hypercohomology spectral sequence, see \cite{GH}, p. 445.
\end{rem}

\begin{thm} \label{thm-calc_cT}
Let $\cX/A$ be a family of schemes and $\cF$ a coherent sheaf on $\cX$. Then
\begin{enumerate}
\item $\cT^0(\cX/A,\cF) = \cHom_{\cX}(\Omega_{\cX/A},\cF)$.
\item $\cT^1(\cX/A,\cF)$ is supported on the locus where $\cX/A$ is not a smooth morphism.
\item $\cT^2(\cX/A,\cF)$ is supported on the locus where $\cX/A$ is not an lci morphism.
\end{enumerate}
\end{thm}
\begin{proof}
See \cite{LS}, 2.3, Theorem~3.1.5 and Corollary~3.2.2.
\end{proof}

\begin{thm} \label{thm-T^i}
Notation as in \ref{notn-art}. Let $\cX_0/A_0$ be a family of schemes and $\cX/A$ a family extending  $\cX_0/A_0$.
\begin{enumerate}
\item There exists a canonical element $o_{\cX/A}(A') \in T^2(\cX_0/A_0, \cO_{\cX_0} \otimes_{A_0} M)$ such that 
$\Def_{\cX/A}(A') \neq \emptyset$ iff $o_{\cX/A}(A') = 0$.\\
\item If $o_{\cX/A}(A')=0$, $\Def_{\cX/A}(A')$ is a principal homogeneous space under 
	$T^1(\cX_0/A_0, \cO_{\cX_0} \otimes_{A_0} M)$.\\
\item Given $\cX'/A'$ extending $\cX/A$, $\Aut_{\cX/A}(\cX'/A')$ is naturally isomorphic to 
	$T^0(\cX_0/A_0, \cO_{\cX_0} \otimes_{A_0} M)$.
\end{enumerate}
\end{thm}

\begin{proof}
For the local case, we use \cite{LS}, p. 66, Theorem~4.3.3, and p. 50, 2.3.2 
(which shows that $T^i(\cX/A,\cO_{\cX_0/A_0}\otimes_{A_0} M) \cong T^i(\cX_0/A_0,\cO_{\cX_0/A_0} \otimes_{A_0} M)$ for each $i$).
The global case follows formally.
\end{proof}

We now develop a $\bQ$-Gorenstein deformation theory for slc surfaces.
We first state the most obvious definition of a $\bQ$-Gorenstein deformation (we refer to this as `weakly $\bQ$-Gorenstein'),
and then refine our definition.

\begin{defn}
Let $\cX/A$ be a family of slc surfaces.
We say $\cX/A$ is \emph{weakly} $\bQ$-\emph{Gorenstein} if the relative Weil divisor $K_{\cX/A}$ is $\bQ$-Cartier.
\end{defn}

This is the definition suggested in \cite{Ko2}, p. 185, Remark~6.27 and is equivalent to 
the definition used in \cite{KSB}. However, it is not well understood over an Artinian base, 
in particular there is no known obstruction theory for weakly $\bQ$-Gorenstein deformations.
To remedy this we make the following definition.

\begin{defn} \label{def-QG}
Let $\cX/A$ be a family of slc surfaces.
We say $\cX/A$ is $\bQ$-\emph{Gorenstein} if $\omega_{\cX/A}^{[i]}$ commutes with base change for all $i \in \bZ$.
Given a $\bQ$-Gorenstein family of slc surfaces $\cX/A$ and an infinitesimal extension $A' \map A$ in $\underline{C}$, 
write $\Def^{QG}_{\cX/A}(A')$ for the set of $\bQ$-Gorenstein deformations of $\cX/A$ over $A'$.
\end{defn}

\begin{rem}
Note that $\bQ$-Gorenstein implies weakly $\bQ$-Gorenstein by Lemma~\ref{lem-pfbc}.
\end{rem}

\begin{rem} \label{rem-finite}
Equivalently, we require that $\omega_{\cX/A}^{[i]}$ commutes with base change for $1 \le i \le N$, where $NK_X$ is Cartier for
each geometric fibre $X$. 
For, by Lemma~\ref{lem-pfbc}, if $\omega_{\cX/A}^{[N]}$ commutes with base change then $\omega_{\cX/A}^{[N]}$ invertible,
thus locally on $\cX$ we have $\omega_{\cX/A}^{[i]} \cong \omega_{\cX/A}^{[i \mod N]}$.
\end{rem}

This definition was given in \cite{Ko6}, p. 260, Definition~5.2. 
Geometrically, the $\bQ$-Gorenstein deformations of a local slc surface $X$   
are precisely those deformations  which lift to deformations of the index one cover $Z \map X$. This is made more precise in
Proposition~\ref{prop-QG_def} below. This description enables us to prove an analogous 
result to Theorem~\ref{thm-T^i} for $\bQ$-Gorenstein deformations (Theorem~\ref{thm-T^i_{QG}}). In particular, we define an 
obstruction theory for $\bQ$-Gorenstein deformations. 
An outline of this theory was given in a preprint of Hassett (c.f. \cite{Has}).
We also show that, for a smoothing $\cX/T$ of an  slc surface over the 
spectrum of a DVR, $\cX/T$ is $\bQ$-Gorenstein iff it is weakly $\bQ$-Gorenstein.

\begin{prop} \label{prop-QG_def}
Given $\cX/A$ a $\bQ$-Gorenstein family of slc surfaces, with $\cX$ and $A$ local,
fix an index one cover $\cZ \map \cX /A$. This is a $\mu_N$ quotient, where $N=\ind \cX$.
Then $\cZ/A$ is flat, and for any $A \map B$ in $\underline{C}$, $\cZ \otimes_A B$ is an index one cover of
$\cX \otimes_A B$. 

Let $A' \map A$ be a local  
extension in $\underline{C}$, and consider the set of families of slc surfaces $\cZ'/A'$
 extending $\cZ/A$. We have an action of $\mu_N$ on this set coming from the $\mu_N$ action on $\cZ/A$,
let $\cZ'/A'$ be an invariant element. Then $\cX'=(\cZ'/\mu_N)/A'$ is a $\bQ$-Gorenstein family of slc surfaces extending $\cX/A$
and $\cZ' \map \cX'/A'$ is an index one cover. 

In particular, if $A' \map A$ is an infinitesimal extension, we have a natural 
isomorphism
$$ \Def_{\cZ/A}(A')^{\mu_N} \map \Def^{QG}_{\cX/A}(A')$$
$$ [\cZ'/A'] \mapsto [(\cZ'/\mu_N)/A']$$
with inverse $[\cX'/A'] \mapsto [\cZ'/A']$ where $\cZ' \map \cX'/A'$ is the unique index one cover extending $\cZ \map \cX/A$.
\end{prop} 

\begin{proof}
We are given a family $\cX/A$ such that $\omega_{\cX/A}^{[i]}$ commutes with base change for all $i$.
Then in particular $K_{\cX/A}$ is $\bQ$-Cartier and $\ind K_{\cX/A}=N$ by Lemma~\ref{lem-rel_index}. 
Let $\pi \colon \cZ \map \cX$ be an index one cover. 
Then $\cZ$ is flat over $A$, because 
$$\pi_{\star}\cO_{\cZ}=\cO_{\cX} \oplus \omega_{\cX/S} \oplus \cdots \oplus \omega_{\cX/S}^{[N-1]}$$
and $\omega_{\cX/A}^{[i]}$ is flat over $A$ for each $i$ by Lemma~\ref{lem-pfbc_flat}.
Given $A \map B$ in $\underline{C}$, 
the isomorphisms $\omega_{\cX/A}^{[i]} \otimes_A B \cong \omega_{\cX \otimes_A B/B}^{[i]}$ show that 
$\cZ \otimes_A B$ is an index one cover of $\cX \otimes_A B$. 
 
Given $\cZ'/A'$ a $\mu_N$ invariant extension of $\cZ/A$ we write $\cX'=(\cZ'/\mu_N)/A'$, then $\cX'/A'$ is a flat family
extending $\cX/A$.
We claim that the quotient $\pi' \colon \cZ' \map \cX'$ is an index one cover of $\cX'$.
First, $\omega_{\cZ/A}$ is invertible, and $\omega_{\cZ'/A'}$ commutes with base change, thus $\omega_{\cZ'/A'}$ is invertible. 
Thus $\pi^{\star}K_{\cX'/A'}=K_{\cZ'/A'}$ is Cartier. In particular $K_{\cX'/A'}$ is $\bQ$-Cartier, so 
$\ind K_{\cX'/A'}=N$ by Lemma~\ref{lem-rel_index}.
Let $i: {\cX'}^0 \inj \cX'$ denote the open subscheme of $\cX'$ where $\omega_{\cX'/A'}$ is invertible,
and write $j: {\cZ'}^0 \inj \cZ'$ for the corresponding open subscheme of $\cZ'$.
Then $\pi^0 \colon {\cZ'}^0 \map {\cX'}^0$ is an etale $\mu_N$ quotient (because the map is etale over $A$). 
Finally $j_{\star}j^{\star}\cO_{\cZ'}=\cO_{\cZ'}$ by Lemma~\ref{lem-S_2}(1).
So $\cZ' \map \cX'$ is an index one cover, by Proposition~\ref{pd-index_one_cover}.
Moreover, we claim $\omega_{\cX'/A'}^{[i]}$ commutes with base change for all $i$. 
First $\omega_{\cX'/A'}^{[N]}$ commutes with base change since it is invertible, using Lemma~\ref{lem-pfbc}.
Also, by the above, $\cZ'$ is an index one cover of $\cX'$, and $\cZ' \otimes_{A'} A$  is an index one cover $\cZ$ of $\cX$
by assumption. It follows that $\omega_{\cX'/A'}^{[i]} \otimes_{A'} A \map \omega_{\cX/A}^{[i]}$ is an isomorphism for 
$1 \le i \le N-1$.
Now $\omega_{\cX/A}^{[i]}$ commutes with base change for all $i$ by assumption, thus $\omega_{\cX'/A'}^{[i]}$ commutes with 
base change for $1 \le i \le N-1$ using Lemma~\ref{lem-pfbc_closed_pts}(2).
So $\omega_{\cX'/A'}^{[i]}$ commutes with base change for all $i$ by Remark~\ref{rem-finite}.

Finally, if $A' \map A$ is a infinitesimal extension, there is a unique index one cover of $\cX'/A'$ extending $\cZ \map \cX/A$
--- for such an index one cover $\cZ' \map \cX'/A$ is determined by the choice of an isomorphism 
$\omega_{\cX'/A'}^{[N]} \cong \cO_{\cX'}$ extending a given isomorphism $\omega_{\cX/A}^{[N]} \cong \cO_{\cX}$ and thus by
a unit $u \in K=\ker(H^0(\cO_{\cX'}^{\times}) \map H^0(\cO_{\cX}^{\times}))$. Moreover multiplying  $u$ by $v^N$
for some $v \in K$ does not change the isomorphism type of the extension
$\cZ' \map \cX'/A'$ of $\cZ \map \cX/A$. But we can always take $N$th roots in $K$ since $A' \map A$ is infinitesimal, thus
$\cZ' \map \cX'/A'$ is uniquely determined as claimed. The last part of the Proposition follows. 
\end{proof}

\begin{prop} \label{prop-QG_DVR}
Let $T$ be the spectrum of a DVR with generic point $\eta$. 
Let $\cX/T$ be a weakly $\bQ$-Gorenstein family of slc surfaces such that $\cX_{\eta}$ is canonical.
Then $\cX/T$ is $\bQ$-Gorenstein.
\end{prop}
\begin{proof}
We may work locally at $P \in \cX$. Let $\pi \colon \cZ \map \cX$ be an index one cover.
It is enough to show that the map $\pi_0 \colon Z \map X$ of the special fibres is an index one cover of $X$, 
by Proposition~\ref{prop-QG_def}.
Now $Z-\pi^{-1}(P) \map X-P$ extends to an index one cover $Z' \map X$, we need to show that $Z$ is $S_2$
to deduce $Z \cong Z'$. First, $X$ slc and $\cX_{\eta}$ canonical implies that $\cX$  is canonical.
For there exists a finite base change $T' \map T$ such that $\cX'=\cX \times_T T'$ 
admits a semistable resolution, then $\cX'$ is canonical 
by Lemma~\ref{lem-iofa}. Using \cite{KSB}, p. 310, Lemma~3.3 we deduce that $\cX$ is canonical.
Now $\cX$ canonical implies that $\cZ$ is canonical, and canonical singularities are rational so in particular CM.
Hence $Z$ is CM since $Z=(t=0) \subset \cZ$ (where $t$ is a uniformising parameter). This completes the proof.
\end{proof}

\begin{defn}  Let $A \in \underline{C}$ and let $\cX/A$ be a $\bQ$-Gorenstein family of slc surfaces. 
Let $\cF$ be a coherent sheaf on $\cX$.
We define a complex $\cM_{QG}(\cF)^{\cdot}$ in the derived category of
coherent sheaves on $\cX$ as follows : Let $\pi \colon \cZ \map \cX$ be a local index one cover of $\cX$,
a $\mu_N$ quotient, where $N$ is the local index of $\cX$.
Let $\sL_{\cdot}$ be a cotangent complex for $\cZ/A$. 
Define $\cM_{QG}(\cF)^{\cdot}=(\pi_{\star}\cHom_{\cO_{\cZ}}(\sL_{\cdot}, \pi^{\star}\cF))^{\mu_N}$ locally.
Now define $\cT_{QG}^i(\cX/A,\cF) =\sH^i(\cM_{QG}(\cF)^{\cdot})$ and $T_{QG}^i(\cX/A,\cF)=\bH^i(\cM_{QG}(\cF)^{\cdot})$.
\end{defn}

\begin{rem} \label{rem-T^i_{QG}}
Note that, since the functor $\{ \cO_{\cZ} \mbox{ modules} \} \map \{\cO_{\cX} \mbox{ modules} \}$, 
$\cF \mapsto (\pi_{\star}\cF)^{\mu_N}$ 
is exact, we have
$$\cT_{QG}^i(\cX/A,\cF) = (\pi_{\star} \cT^i(\cZ/A, \pi^{\star}\cF))^{\mu_N}$$
locally.
\end{rem}

\begin{rem} \label{rem-QG_ss}
We have a local-to-global spectral sequence
$$E^{pq}_2 =H^p(\cX,\cT_{QG}^q(\cX/A,\cF)) \Rightarrow T_{QG}^{p+q}(\cX/A,\cF)$$
as above.
\end{rem}

\begin{lem} \label{lem-T^0_{QG}}
Let $A \in \underline{C}$, let $\cX/A$ be a $\bQ$-Gorenstein family of slc surfaces, and $\cF=\cO_{\cX}\otimes_{A} M$,
some finite $A$-module $M$. Then $\cT^0_{QG}(\cX/A,\cF)=\cT^0(\cX/A,\cF)$.
\end{lem}
\begin{proof}
Let $i \colon \cX^0 \inj \cX$ be the inclusion of the locus where $\omega_{\cX/A}$ is invertible. 
Let $\pi \colon \cZ \map \cX$ be a local index one cover, a $\mu_N$ quotient say, 
and write $j \colon \cZ^0 \inj \cZ$ for the inverse image of $\cX^0$, then $\cZ^0 \map \cX^0$ is etale.
We have $\cT^0(\cX/A,\cF)= \cHom_{\cO_{\cX}}(\Omega_{\cX/A},\cF)$ and $i_{\star}i^{\star}\cF=\cF$ using Lemma~\ref{lem-S_2},
thus $i_{\star}i^{\star}\cT^0(\cX/A,\cF)=\cT^0(\cX/A,\cF)$. Similiarly, 
$j_{\star}j^{\star}\cT^0(\cZ/A,\pi^{\star}\cF)=\cT^0(\cZ/A,\pi^{\star}\cF)$, so,  using 
$\cT^0_{QG}(\cX/A,\cF)=(\pi_{\star}\cT^0(\cZ/A,\pi^{\star}\cF))^{\mu_N}$, we find 
$i_{\star}i^{\star}\cT^0_{QG}(\cX/A,\cF)=\cT^0_{QG}(\cX/A,\cF)$.
But $\cT^0_{QG}(\cX/A,\cF)$ and $\cT^0(\cX/A,\cF)$ agree on $\cX^0$, so we obtain our result.
\end{proof}

\begin{thm} \label{thm-T^i_{QG}}
Notation as in \ref{notn-art}. Let $\cX_0/A_0$ be a $\bQ$-Gorenstein family of slc surfaces and 
$\cX/A$ a $\bQ$-Gorenstein family of slc surfaces extending $\cX_0/A_0$.
\begin{enumerate}
\item There exists a canonical element $o^{QG}_{\cX/A}(A') \in T_{QG}^2(\cX_0/A_0, \cO_{\cX_0} \otimes_{A_0} M)$ such that
$\Def^{QG}_{\cX/A}(A') \neq \emptyset$ iff $o^{QG}_{\cX/A}(A') \neq 0$.\\
\item If $o^{QG}_{\cX/A}(A')=0$, $\Def^{QG}_{\cX/A}(A')$ is a principal homogeneous space under 
	$T_{QG}^1(\cX_0/A_0, \cO_{\cX_0} \otimes_{A_0} M)$.\\
\item 	Given a $\bQ$-Gorenstein family $\cX'/A'$ extending $\cX/A$, $\Aut_{\cX/A}(\cX'/A')$ is naturally isomorphic to 
	$T^0(\cX_0/A_0, \cO_{\cX_0} \otimes_{A_0} M)$.
\end{enumerate}
\end{thm}

\begin{proof}
For the local case we use Theorem~\ref{thm-T^i} together with Proposition~\ref{prop-QG_def}.
We may assume $A_0$ is a local ring.
Given a $\bQ$-Gorenstein family of slc surfaces $\cX/A$ extending $\cX_0/A_0$, 
take an index one cover $\cZ \map \cX /A$ (a $\mu_N$ quotient, say). 
This is flat over $A$ and restricts to an index one cover
$\cZ_0 \map \cX_0 /A_0$ by Proposition~\ref{prop-QG_def}.
Now, by Theorem~\ref{thm-T^i}(1), there is a canonical element $o_{\cZ/A}(A') \in T^2(\cZ_0/A_0,\cO_{\cZ_0} \otimes_{A_0} M)$ 
such that
$o_{\cZ/A}(A')=0$ iff $\Def_{\cZ/A}(A') \neq \emptyset$. By functoriality of the obstruction map, we have
$o_{\cZ/A}(A') \in T^2(\cZ_0/A_0,\cO_{\cZ_0} \otimes_{A_0} M)^{\mu_N}=T^2_{QG}(\cX_0/A_0,\cO_{\cX_0} \otimes_{A_0} M)$ using
 Remark~\ref{rem-T^i_{QG}}. 
We define $o^{QG}_{\cX/A}(A') = o_{\cZ/A}(A') \in T^2_{QG}(\cX_0/A_0,\cO_{\cX_0} \otimes_{A_0} M)$.
In order to show property (1) holds, we just need to verify that if $\Def_{\cZ/A}(A') \neq \emptyset$ then 
$\Def_{\cZ/A}(A')^{\mu_N} \neq \emptyset$, for by Proposition~\ref{prop-QG_def} we have a bijection 
$\Def_{\cZ/A}(A')^{\mu_N} \map \Def^{QG}_{\cX/A}(A')$ --- we use Lemma~\ref{lem-G_invt} below.
Property (2) follows immediately from Proposition~\ref{prop-QG_def} and Theorem~\ref{thm-T^i}(2) since 
$T^1_{QG}(\cX_0/A_0,\cO_{\cX_0} \otimes_{A_0} M)=T^1(\cZ_0/A_0,\cO_{\cZ_0} \otimes_{A_0} M)^{\mu_N}$.
Finally, (3) is a special case of Theorem~\ref{thm-T^i}(3).

The global case now follows formally exactly as for ordinary deformation theory using the spectral sequence of
Remark~\ref{rem-QG_ss} (we use Lemma~\ref{lem-T^0_{QG}} to identify
$\cT^0_{QG}$ with $\cT^0$, the sheaf of infinitesimal automorphisms). 
\end{proof}

\begin{lem} \label{lem-G_invt}
Let $A' \map A$ be an infinitesimal extension in $\underline{C}$ with $A$ local and $M^2=0$,
where $M = \kernel(A' \map A)$. 
Let $\cX/A$ be a family of local schemes, $\cX'/A'$ a family extending $\cX/A$, and 
$\cZ \subset \cX/A$ a family of closed subschemes. Suppose given a $\mu_N$ action on $(\cX,\cZ)/A$ which extends to
a $\mu_N$ action on $\cX'/A'$.
Then, if there exists an extension $\cZ' \subset \cX'/A'$ of $\cZ \subset \cX/A$, there exists a $\mu_N$-invariant extension.
\end{lem}
\begin{proof}
Write  $\cX=\Spec P$, $\cX'= \Spec P'$, and $\cI_{\cZ}=J=(F_1,\ldots,F_N) \subset P$.
Let $J' \subset P'$ be an ideal defining a closed subscheme $\cZ' \subset \cX'/A'$ (not necessarily flat over $A'$) extending
$\cZ \subset \cX/A$. Then $\cZ'/A'$ is flat iff $J'=(F_1',\ldots,F_N')$ for some liftings $F_i' \in P'$ of $F_i \in P$ and 
every relation $\sum R_i F_i =0$, $R_1, \ldots, R_N \in P$, between the $F_i$ 
lifts to a relation $\sum R'_i F'_i =0$, $R'_1, \ldots, R'_N \in P'$, between the $F'_i$ (c.f. \cite{Ar2}).

In our case, pick generators $F_1,\ldots, F_N$ of $J$ which are $\mu_N$ eigenfunctions (using \cite{KM},
p. 219, Lemma~7.29). 
Pick $\mu_N$ eigenfunctions $F'_1, \ldots, F'_N \in P'$ lifting $F_1,\ldots, F_N$. 
Let
$$0 \map R \map P^N \stackrel{(F_1,\ldots,F_N)}{\map} J \map 0,$$
be exact. We define a $\mu_N$ action on $P^N$ (in the obvious way) such that $P^N \map J$ is $\mu_N$ equivariant, 
and hence obtain a $\mu_N$ action on the module of relations $R \subset P^N$.
We have a well defined map
$$\phi \colon R \map M \otimes_A P/J$$
$$ (R_1,\ldots,R_N) \mapsto \sum R'_iF'_i,$$
where $R'_i \in P'$ are \emph{some} lifts of the $R_i$.
Here we first regard $\sum R'_iF'_i$ as an element of $M \otimes_A P= \kernel (P' \map P)$
and then take its image in $M \otimes_A P/J$. Then $\phi=0$ iff $\cZ'=(F'_1,\ldots,F'_N=0) \subset \cX'/A'$ is flat over $A'$,
by the criterion above.
Thus there exists a flat extension $\tilde{\cZ'} \subset \cX'/A'$ of $\cZ \subset \cX/A$ iff 
$$\phi \in \im( \Hom (P^N, M \otimes_A P/J) \map \Hom(R, M \otimes_A P/J) )$$
--- for this is the condition that we can change the $F'_i$ so that $\phi$ becomes the zero map.
Now, by construction, $\phi$ is a $\mu_N$ invariant element of $\Hom(R, M \otimes_A P/J)$, so equivalently we require that 
$$\phi \in \im( \Hom (P^N, M \otimes_A P/J)^{\mu_N} \map \Hom(R, M \otimes_A P/J)^{\mu_N} ).$$  
In this case we can replace the $F'_i$ by $\mu_N$ eigenfunctions $\tilde{F'_i}$ to obtain a $\mu_N$ invariant flat extension
$\tilde{\cZ'} \subset \cX'/A'$. This completes the proof.
\end{proof}

\section{Construction of the stack $\cM_d$} \label{construct}

The aim of this section is to prove that the groupoid $\cM_d$ defined in Section~\ref{setup} is an algebraic stack.
We use the theory of Artin \cite{Ar1} reviewed in Section~\ref{artin}.
We do not work directly with $\cM_d$, instead we define a related groupoid
$F_d$ which (roughly) has the same definition as $\cM_d$ except that we drop the smoothability assumption. 
We show that $F_d$ is an algebraic stack (not necessarily proper in general) and then obtain $\cM_d$ as a closed substack.

We first verify that $F_d$ satisfies the conditions of Theorem~\ref{thm-patches}, and hence obtain local patches of the stack. 
The bulk of the work is giving a concrete description of $D$ and constructing an obstruction theory $\cO$ 
(Theorem~\ref{thm-T^i_{QG}(X,D)}).
Given the local patches it's easy to show that $F_d$ is an algebraic stack 
(we just need to show that $F_d$ is `relatively representable').     

One might think that we could construct $\cM_d$ 
as a quotient of some locally closed subscheme of a Hilbert scheme of pairs $(X,D)$.
However, we do not know that the base change conditions for $(\cX,\sD)/S \in \cM_d(S)$ are locally closed 
(c.f. Theorem~\ref{thm-pfbc_stratification} and Remark~\ref{rem-stratification}).
Thus we cannot obtain $\cM_d$ in this way.

\begin{defn} \label{def-quasistable}
Let $X$ be a proper connected surface.
Let $D$ be an effective Weil divisor on $X$. Let $d \in \bN$, $d \ge 4$.
$(X,D)$ is a \emph{quasistable pair of degree $d$} if
\begin{enumerate}
\item There exists $\epsilon > 0$ such that $K_X+(\frac{3}{d}+\epsilon)D$ is slc and ample.\\
\item $dK_X+3D \sim 0$, and moreover $\frac{d}{3}K_X+D \sim 0$ if $3 \mid d$.\\
\item $3 \nd \ind D$  if $3 \nd d$, $H^1(\cO_X(D))=0$, and $\chi(\cO_X)=1$.
\end{enumerate}
\end{defn}

\begin{lem} \label{lem-stable_and_quasistable}
A stable pair is quasistable.
\end{lem}
\begin{proof}
We just need to show that the smoothability property for a stable pair $(X,D)$ of degree $d$ gives  
$3 \nd \ind D$  if $3 \nd d$, $H^1(\cO_X(D))=0$, and $\chi(\cO_X)=1$.
First note that trivially $\chi(\cO_X)=\chi(\cO_{\bP^2})=1$ if $X$ admits a smoothing to $\bP^2$.

If $3 \nd d$ the only possible singularities of $X$ are of local analytic types $\frac{1}{n^2}(1,na-1)$ where $3 \nd n$
and $(xy=0) \subset \frac{1}{r}(1,-1,a)$ (using Theorem~\ref{thm-simp}, Theorem~\ref{thm-smoothability} and 
Lemma~\ref{lem-T_1}). In the second case $dK_X+3D \sim 0$ shows $3 \nd r$. 
Thus $3 \nd \ind D$.

We now prove $H^1(\cO_X(D))=0$. 
First note that $h^1(\cO_X(D))=$ \mbox{$h^1(\cO_X(K_X-D))$} by Serre duality, and $-(K_X-D)$ is ample.
We have an exact sequence
$$0 \map \cO_X(K_X-D) \map \cO_{X^{\nu}}(\nu^{\star}(K_X-D)) \map \cO_{\hat{\Delta}}(\lfloor K_X-D |_{\hat{\Delta}} \rfloor)$$
where $\nu \colon X^{\nu} \map X$ is the normalisation of $X$, and $\hat{\Delta} \map \Delta$ is the normalisation of the double 
curve $\Delta$ of $X$. We obtain a short exact sequence
$$0 \map \cO_X(K_X-D) \map \cO_{X^{\nu}}(\nu^{\star}(K_X-D)) \map \cF \map 0$$
where $\cF \inj \cO_{\hat{\Delta}}(\lfloor K_X-D |_{\hat{\Delta}} \rfloor)$. In particular, since $-(K_X-D)$ is ample,
we have $H^0(\cF)=0$. Thus the long exact sequence of cohomology associated to the short exact sequence above gives
$H^1(\cO_X(K_X-D)) \inj$ \mbox{$H^1(\cO_{X^{\nu}}(\nu^{\star}(K_X-D)))$}. If $X^{\nu}$ is klt, then 
$H^1(\cO_{X^{\nu}}(\nu^{\star}(K_X-D)))=0$ by Kodaira vanishing. Otherwise, $X$ is an elliptic cone by 
Theorem~\ref{thm-glueing} and Theorem~\ref{thm-normal-case}, and $3 \mid d$, $D \sim -\frac{d}{3}K_X$.
An easy calculation shows that $H^1(\cO_X(D))=0$ in this case.
\end{proof}

\begin{defn}
We say that $(\cX,\sD)/S$ is a family of quasistable pairs of degree $d$ over $S$ if $\cX$ is a flat family over $S$,  
$\sD \subset \cX$ is a relative Weil divisor over $S$,
and for every geometric point $s$ of $S$, the fibre $(\cX_s,\sD_{(s)})$ is a quasistable pair of degree $d$.
We say that $(\cX,\sD)/S$ is an \emph{allowable family} if 
$\omega_{\cX/S}^{[i]}$ and $\cO_{\cX}(i\sD)$ commute with base change for all $i \in \bZ$.
\end{defn}

\begin{defn}
We define a groupoid $F_d$ over $\underline{Sch}$ as follows: 
For $S \in \underline{Sch}$ let
$$F_d(S)=\{ \mbox{ Allowable families of quasistable pairs of degree $d$ over } S \}.$$
\end{defn}

\begin{rem} 
We also regard $F_d$ as a groupoid over the category $\underline{C}$ of noetherian $\bC$-algebras 
without further comment.
\end{rem}

\begin{lem}
$\cM_d$ is a subgroupoid of $F_d$
\end{lem}
\begin{proof}
This is immediate from Lemma~\ref{lem-stable_and_quasistable} --- $\cM_d \subset F_d$ is the subgroupoid of smoothable families.
\end{proof}

In what follows, we suppress the degree $d$ to simplify our notation.

\begin{prop} \label{prop-pfbc_DVR}
Let $T$ be the spectrum of a DVR with generic point $\eta$.
A family $(\cX,\sD)/T$ of quasistable pairs such that $\cX_{\eta}$ is canonical 
is allowable iff $K_{\cX/T}$ and $\sD$ are $\bQ$-Cartier.
\end{prop}
\begin{proof}
The 'only if' part follows from Lemma~\ref{lem-pfbc}.
So, suppose given a family $(\cX,\sD)/T$ of quasistable pairs such that $\cX_{\eta}$ is canonical and  $K_{\cX/T}$ and $\sD$ are 
$\bQ$-Cartier. Then $\cX/T$ is a $\bQ$-Gorenstein family by Proposition~\ref{prop-QG_DVR}.
Let $\cZ \map \cX$ be an index one cover, with special fibre $Z \map X$ an index one cover of $X$.
Since $\sD$ is $\bQ$-Cartier, 
and $D_Z=\pi_0^{\star}D$ is Cartier (using $dK_X+3D \sim 0$ and $3 \nd \ind D$), we have $\sD_{\cZ}=\pi^{\star}\sD$
Cartier, using Lemma~\ref{lem-rel_index}. It follows that $\sD \sim jK_{\cX/T}$ locally, some $j \in \bZ$, hence
$\cO_{\cX}(i\sD)$ commutes with base change for all $i$. Thus $(\cX,\sD)/T$ is allowable as required.
\end{proof}

\begin{lem} \label{lem-Emb_and_F}
For $A \in \underline{C}$, $(\cX,\sD)/A \in F(A)$ iff the following conditions are satisfied:
\begin{enumerate}
\item $\cX/A$ is a $\bQ$-Gorenstein family of slc surfaces.
\item $\sD \subset \cX$ is a codimension $1$ closed subscheme, flat over $A$.
\item $(\cX_s,\sD_{(s)})$ is a quasistable pair of degree $d$ for every geometric point $s$ of $S$.
\item $\sD_{s}=\sD_{(s)}$ for every geometric point $s$ of $S$.
\end{enumerate}
Here $\sD_s$ is the scheme theoretic fibre of $\sD$ over $s$, and $\sD_{(s)}$ is the restriction of $\sD$ to the fibre $\cX_s$
defined by taking the double dual of the ideal sheaf of $\sD_s$. 
\end{lem}

\begin{proof}
Given $(\cX,\sD)/A \in F(A)$, $\sD$ is flat over $A$ and $\sD_s=\sD_{(s)}$ for all $s$ by Lemma~\ref{lem-rwd_ideal_pfbc}.
The other conditions are satisfied by the definition of $F$.

Conversely let $(\cX,\sD)/A$ satisfy the conditions (1) to (4). We claim that $(\cX,\sD)/A \in F(A)$.
Note that $\sD$ is a relative Weil divisor (c.f. Example~\ref{ex-rwd}).
It is enough to show that $\cO_{\cX}(i\sD)$ commutes with base change for all $i \in \bZ$.
By Lemma~\ref{lem-pfbc_inv_lim}, we may assume that $A$ is a local Artinian ring.
So, by induction, it is enough to show the following:  
Let $A$ be a local ring with residue field $k$ and $A' \map A$ a small extension 
(i.e., the kernel $J$ of $A' \map A$ is annihilated by the maximal ideal of $A$).
Suppose given $(\cX,\sD)/A \in F(A)$ and $(\cX',\sD')/A'$ extending $(\cX,\sD)/A$ which satisfies  conditions (1) to (4). 
Then $(\cX',\sD')/A' \in F(A')$.

We work locally on $\cX'$.
Let $\pi' \colon \cZ' \map \cX'$ be an index one cover (a $\mu_N$ quotient, say), 
$\pi \colon \cZ \map \cX$ the index one cover obtained by restriction 
to $A$, and $\sD_{\cZ}=\pi^{\star}\sD$. 
Given an $\mu_N$ invariant extension $\sD_{\cZ'} \subset \cZ'$ of $\sD_{\cZ} \subset \cZ$ we obtain an extension 
$\sD'=(\sD_{\cZ'} / \mu_N) \subset \cX'$ of $\sD \subset \cX$. We claim that every extension $\sD' \subset \cX'$ of
$\sD \subset \cX$ occurs in this way.
Let $\pi_0 \colon (Z,D_Z) \map (X,D)$ be the special fibre of $\pi \colon (\cZ,\sD_{\cZ}) \map (\cX,\sD)$.
Then $Z \map X$ is an index one cover, and since $dK_X+3D \sim 0$ and $3 \nd \ind D$ we have $D_Z$ Cartier.
Thus $D_Z$ has unobstructed embedded deformations (locally) and the extensions $\sD_{\cZ'} \subset \cZ'$  
of $\sD_{\cZ} \subset \cZ$ form a prinicipal homogeneous space under $\Hom_Z(\cI_{D_Z},\cO_{D_Z}) \otimes_k J$ (using 
Theorem~\ref{thm-Emb}). There exists a $\mu_N$ invariant extension by Lemma~\ref{lem-G_invt}, hence also
$D$ has unobstructed embedded deformations (locally) and the extensions $\sD' \subset \cX'$  
of $\sD \subset \cX$ form a prinicipal homogeneous space under $\Hom_X(\cI_D,\cO_D) \otimes_k J$.
To prove our claim, we just need to identify the sheaves $({\pi_0}_{\star}\cHom_Z(\cI_{D_Z},\cO_{D_Z}))^{\mu_N}$
and $\cHom_X(\cI_D,\cO_D)$. We have the exact sequence
$$0 \map \cI_D \map \cO_X \map \cO_D \map 0,$$
applying the functor $\cHom_X(\cI_D,\cdot)$ we obtain 
$$0 \map \cO_X \map \cO_X(D) \map \cHom_X(\cI_D,\cO_D) \map \cExt^1_X(\cI_D,\cI_D) \map \cdots.$$
Now $\cExt^1_X(\cI_D,\cI_D)=\cExt^1_X(\cO_X,\cO_X)=0$, since $\cI_D$ is invertible in codimension $1$ and $S_2$.
Hence we have a short exact sequence
$$0 \map \cO_X \map \cO_X(D) \map \cHom_X(\cI_D,\cO_D) \map 0,$$
and similiarly
$$0 \map \cO_Z \map \cO_Z(D_Z) \map \cHom_Z(\cI_{D_Z},\cO_{D_Z}) \map 0.$$
Applying the exact functor $({\pi_0}_{\star}( \cdot))^{\mu_N}$ to the last exact sequence we obtain
$$0 \map \cO_X \map \cO_X(D) \map ({\pi_0}_{\star}\cHom_Z(\cI_{D_Z},\cO_{D_Z}))^{\mu_N} \map 0,$$
thus $({\pi_0}_{\star}\cHom_Z(\cI_{D_Z},\cO_{D_Z}))^{\mu_N}=\cHom_X(\cI_D,\cO_D)$ as required.
So, given an extension $\sD' \subset \cX'$ of $\sD \subset \cX$, it is obtained from an extension
$\sD_{\cZ'} \subset \cZ'$ of $\sD_{\cZ} \subset \cZ$ by taking the $\mu_N$-quotient. In particular,
$\sD_{\cZ'}={\pi'}^{\star}\sD'$ is Cartier. Thus locally $\sD' \sim jK_{\cX'/A'}$ for some $j$, hence
$\cO_{\cX'}(i\sD')$ commutes with base change for all $i \in \bZ$ as required.
\end{proof}

\begin{thm} \label{thm-F_to_Def^QG_smooth}
Let $A \in \underline{C}$ and $(\cX,\sD)/A \in F(A)$.
The forgetful map of functors (of infinitesimal extensions of $A$)
$$\bar{F}_{(\cX,\sD)/A} \map \Def^{QG}_{\cX/A}$$ 
is smooth.
\end{thm}

\begin{proof}
Given an infinitesimal extension $A' \map A$ in $\underline{C}$ and a $\bQ$-Gorenstein family $\cX'/A'$ extending $\cX/A$,
we need to show that there exists an extension $\sD' \subset \cX'$ of $\sD \subset \cX$ (using Lemma~\ref{lem-Emb_and_F}).
By the proof of Lemma~\ref{lem-Emb_and_F} there are no local obstructions.
We may assume that $A$ is local and $A' \map A$ is a small extension. 
Let $A$ have residue field $k$ and write $(X,D)/k$ for the special fibre of $(\cX,\sD)/A$.
Then the obstruction to extending $\sD \subset \cX$ to some
$\sD' \subset \cX'$ lies in $\Ext^1_X(\cI_D,\cO_X)$ by Theorem~\ref{thm-Emb}(1).
Moreover, since there are no local obstructions, using the local-to-global spectral sequence for $\Ext$ we see that the 
obstruction lies in $H^1(\cHom_X(\cI_D,\cO_X))$.
We have an exact sequence
$$0 \map \cO_X \map \cO_X(D) \map \cHom_X(\cI_D,\cO_D) \map 0$$
(see the proof of Lemma~\ref{lem-Emb_and_F}). Now $H^1(\cO_X(D))=0$ by assumption, and $h^2(\cO_X)=h^0(K_X)=0$ since $-K_X$ 
is ample, thus $H^1(\cHom_X(\cI_D,\cO_X))=0$. So the obstruction is zero, and there is an extension $\sD' \subset \cX'$
of $\sD \subset \cX$ as required.
\end{proof}

\begin{notn} \label{notn-emb}
Let $A' \map A$ be an infinitesimal extension in $\underline{C}$, $\cX/A$ a family of schemes over $A$, $\cZ \subset \cX/A$
a family of closed subschemes, and $\cX'/A'$ a family extending $\cX/A$. 
Write $\Emb_{(\cX,\cZ)/A}(\cX'/A')$ for the set of embedded deformations of $\cZ/A$ over $A'$ inside $\cX'/A'$.
Write $\overline{\Emb}_{(\cX,\cZ)/A}(\cX'/A')$ for the quotient of $\Emb_{(\cX,\cZ)/A}(\cX'/A')$ by the action of 
$\Aut_{\cX/A}(\cX'/A')$. Given $\cZ' \subset \cX' /A'$ extending $\cZ \subset \cX/A$, write
$\Aut_{(\cX,\cZ)/A}((\cX',\cZ')/A')$ for the group of automorphisms of $(\cX',\cZ')/A'$ over $A'$ 
which restrict to the identity on $(\cX,\cZ)/A$.
Given $A \map B$ in $\underline{C}$, we also write $\Emb_{(\cX,\cZ)/A}(B)$ for $\Emb_{(\cX,\cZ)/A}(\cX \otimes_A B /B)$,
similiarly for $\overline{\Emb}_{(\cX,\cZ)/A}$ and ${\Aut}_{(\cX,\cZ)/A}$. 
\end{notn}

\begin{thm} \label{thm-Emb}
Notation as in \ref{notn-art}. Let $\cX_0/A_0$ be a family of schemes, $\cX/A$ a family extending $\cX_0/A_0$, and $\cX'/A'$
a family extending $\cX/A$. Let $\cZ_0 \subset \cX_0/A_0$ be a family of closed subschemes, and $\cZ \subset \cX/A$ a family
of closed subschemes extending $\cZ_0 \subset \cX_0 /A_0$.
\begin{enumerate}
\item There is a canonical element $o_{(\cX,\cZ)/A}(A') \in \Ext^1_{\cX_0}(\cI_{\cZ_0}, \cO_{\cZ_0} \otimes_{A_0} M)$ such that
$o_{(\cX,\cZ)/A}(A')=0$ iff $\Emb_{(\cX,\cZ)/A}(\cX'/A') \neq \emptyset$.\\
\item 	If $o_{(\cX,\cZ)/A}(A')=0$, $\Emb_{(\cX,\cZ)/A}(\cX'/A')$ is a principal homogeneous space under 
	$\Hom_{\cX_0}(\cI_{\cZ_0},\cO_{\cZ_0}\otimes_{A_0}M)$.\\
\item 	We have a natural map $$\phi \colon \Aut_{\cX_0/A}(A_0+M) \map \Hom_{\cX_0}(\cI_{\cZ_0},\cO_{\cZ_0} \otimes_{A_0} M)$$
	which identifies the action of $\Aut_{\cX_0/A}(A_0+M)$ on \mbox{$\Emb_{(\cX_0,\cZ_0)/A_0}(A_0+M)$} 
	in terms of the action of 
	\mbox{$\Hom_{\cX_0}(\cI_{\cZ_0},\cO_{\cZ_0} \otimes_{A_0} M)$}.
	
	If $\overline{\Emb}_{(\cX,\cZ)/A}(\cX'/A') \neq \emptyset$, it is a principal homogeneous space under $\cok\phi$. 
	Given  a family of closed subschemes $\cZ' \subset \cX'/A'$ extending $\cZ \subset \cX/A$, 
	$\Aut_{(\cX,\cZ)/A}((\cX',\cZ')/A')$ is naturally isomorphic to $\ker\phi$.
\end{enumerate}
\end{thm}

\begin{proof}
For (1) and (2) see \cite{Ko4}, p. 28, Proposition~2.5 (in fact only the case $A_0=k$, $A$ local is treated in \cite{Ko4}, 
but the same argument proves the general case). 
We now prove (3). If $\Emb_{(\cX,\cZ)/A}(A') \neq \emptyset$, it is a principal homogeneous space under
$\Hom_{\cX_0}(\cI_{\cZ_0},\cO_{\cZ_0} \otimes_{A_0} M)$ by (2). Now $\Aut_{\cX/A}(\cX'/A')$ acts on  $\Emb_{(\cX,\cZ)/A}(A')$,
so we obtain a homomorphism 
$$\Aut_{\cX/A}(\cX'/A') \map \Hom_{\cX_0}(\cI_{\cZ_0},\cO_{\cZ_0} \otimes_{A_0} M).$$
In the case $A=A_0$, $A'=A_0+M$, $\cX'=\cX \otimes_A A'$ we obtain the map $\phi$ above.
We have a natural isomorphism $\Aut_{\cX/A}(\cX'/A') \cong \Aut_{\cX_0/A_0}(A_0+M)$ by Theorem~\ref{thm-T^i}(3),
which is compatible with the maps to \mbox{$\Hom_{\cX_0}(\cI_{\cZ_0},\cO_{\cZ_0} \otimes_{A_0} M)$}.
The result now follows.
\end{proof}

\begin{defn} \label{def-T^i_{QG}(X,D)}
Given $A \in \underline{C}$, $M$ a finite $A$-module and $(\cX,\sD)/A \in F(A)$,
define 
$$T^0((\cX,\sD)/A,\cO_{\cX}\otimes_{A} M)=\Aut_{(\cX,\sD)/A}(A+M),$$
and 
$$T^1_{QG}((\cX,\sD)/A,\cO_{\cX} \otimes_{A} M)=\bar{F}_{(\cX,\sD)/A}(A+M).$$
\end{defn}

\begin{lem} \label{lem-T^i_{QG}(X,D)}
Notation as Definition~\ref{def-T^i_{QG}(X,D)}. The sets 
$T^0((\cX,\sD)/A,\cO_{\cX}\otimes_{A} M)$ and $T^1_{QG}((\cX,\sD)/A,\cO_{\cX}\otimes_{A} M)$ have natural $A$-module structures, 
and, writing $\cF=\cO_{\cX}\otimes_{A} M$, we have an exact sequence of $A$-modules
$$0 \map T^0((\cX,\sD)/A,\cF) \map T^0(\cX/A,\cF) \map \Hom_{\cX}(\cI_{\sD},\cO_{\sD} \otimes_A M) \map$$  
$$\map T^1_{QG}((\cX,\sD)/A,\cF) \map T^1_{QG}(\cX/A,\cF) \map 0.$$
\end{lem}

\begin{proof}
We have $T^0(\cX/A,\cF)=\Aut_{\cX/A}(A+M)$ and $\Hom_{\cX}(\cI_{\sD},\cO_{\sD}\otimes M)=\Emb_{(\cX,\sD)/A}(A+M)$, so we have
a natural map $$\phi \colon T^0(\cX/A,\cF) \map \Hom_{\cX}(\cI_{\sD},\cO_{\sD}\otimes M)$$ with kernel 
$T^0((\cX,\sD)/A,\cF)=\Aut_{(\cX,\sD)/A}(A+M)$ and cokernel \mbox{$\overline{\Emb}_{(\cX,\sD)/A}(A+M)$} 
(compare Theorem~\ref{thm-Emb}(3)).
In particular, $T^0((\cX,\sD)/A,\cF)$ is an $A$-module, moreover $T^1((\cX,\sD)/A,\cF)=\bar{F}_{(\cX,\sD)/A}(A+M)$ is an $A$-module 
since the functor $\bar{F}$ satisfies Artin's criterion S1(b).
Now $T^1_{QG}((\cX,\sD)/A,\cF)=$ \mbox{$\bar{F}_{(\cX,\sD)/A}(A+M)$} and $T^1_{QG}(\cX/A,\cF)=\Def^{QG}_{\cX/A}(A+M)$ so we have a natural map
$$\psi \colon T^1_{QG}((\cX,\sD)/A,\cF) \map T^1_{QG}(\cX/A,\cF),$$ 
with kernel \mbox{$\overline{\Emb}_{(\cX,\sD)/A}(A+M)$}
by Lemma~\ref{lem-Emb_and_F}.
Finally, $\psi$ is surjective by Theorem~\ref{thm-F_to_Def^QG_smooth}.
\end{proof}

\begin{thm} \label{thm-T^i_{QG}(X,D)}
Notation as in \ref{notn-art}. Let $(\cX_0,\sD_0)/A_0 \in F(A_0)$ and $(\cX,\sD)/A \in F_{(\cX_0,\sD_0)/A_0}(A)$.
\begin{enumerate}
\item We have $F_{(\cX,\sD)/A}(A') \neq \emptyset$ iff $o^{QG}_{\cX/A}(A') \neq 0$, where 
$o^{QG}_{\cX/A}(A') \in T_{QG}^2(\cX_0/A_0, \cO_{\cX_0} \otimes_{A_0} M)$ is the canonical element constructed above.\\
\item If $o^{QG}_{\cX/A}( A')=0$, $\bar{F}_{(\cX,\sD)/A}(A')$ is a principal homogeneous space under 
	$T_{QG}^1((\cX_0,\sD_0)/A_0, \cO_{\cX_0} \otimes_{A_0} M)$.\\
\item 	Given $(\cX',\sD')/A' \in F_{(\cX,\sD)/A}(A')$, $\Aut_{(\cX,\sD)/A}((\cX',\sD')/A')$ is naturally isomorphic to 
	$T^0((\cX_0,\sD_0)/A_0, \cO_{\cX_0} \otimes_{A_0} M)$.
\end{enumerate}
\end{thm}

\begin{proof}
Property (1) is immediate from Theorem~\ref{thm-T^i_{QG}} and Theorem~\ref{thm-F_to_Def^QG_smooth} --- because the forgetful map 
$\bar{F}_{(\cX,\sD)/A} \map \Def^{QG}_{\cX/A}$ is smooth, the obstruction theory for $\Def^{QG}$ gives an obstruction theory 
for $F$.

We have a natural map $\psi \colon \bar{F}_{(\cX,\sD)/A}(A') \map \Def^{QG}_{\cX/A}(A')$ which is surjective
by Theorem~\ref{thm-F_to_Def^QG_smooth}.
Assuming $o^{QG}_{\cX/A}(A')=0$, $\Def^{QG}_{\cX/A}(A')$ is a principal homogeneous space under 
$T^1_{QG}(\cX_0/A_0,\cO_{\cX_0}\otimes_{A_0}M)$ by Theorem~\ref{thm-T^i_{QG}}(2). 
Moreover, given $[\cX'/A'] \in \Def^{QG}_{\cX/A}(A')$, $\psi^{-1}([\cX'/A'])$ is the set 
$\overline{\Emb}_{(\cX,\sD)/A}(\cX'/A')$. This is a principal homogeneous space under
$\cok\phi$, where $$\phi \colon T^0(\cX_0/A_0,\cO_{\cX_0}\otimes_{A_0}M) \map \Hom_{\cX_0}(\cI_{\sD_0},\cO_{\sD_0})$$ 
as in Theorem~\ref{thm-Emb}(3).
Now, the set $\bar{F}_{(\cX,\sD)/A}(A')$ has a natural action of 
$\bar{F}_{(\cX_0,\sD_0)/A_0}(A_0+M)=T^1_{QG}((\cX_0,\sD_0)/A_0,\cO_{\cX_0} \otimes_{A_0} M)$, 
since the functor $\bar{F}$ satisfies Artin's criterion S1(b).
This action is compatible with the actions just described via the exact sequence of Lemma~\ref{lem-T^i_{QG}(X,D)}, 
it follows that $\bar{F}_{(\cX,\sD)/A}(A')$ is a principal homogeneous space under 
\mbox{$T^1_{QG}((\cX_0,\sD_0)/A_0,\cO_{\cX_0} \otimes_{A_0} M)$} as required.
Finally, property (3) is a special case of Theorem~\ref{thm-Emb}(3).
\end{proof}

We can now identify the functor $D$ and an obstruction theory $\cO$ for our groupoid $F$.
Given $(\cX_0,\sD_0)/A_0 \in F(A_0)$ and a finite $A_0$-module $M$ we have
$$D_{(\cX_0,\sD_0)/A_0}(M)=\bar{F}_{(\cX_0,\sD_0)/A_0}(A_0+M) = T^1_{QG}((\cX_0,\sD_0)/A_0,\cO_{\cX_0} \otimes_{A_0} M)$$
using Theorem~\ref{thm-T^i_{QG}(X,D)}(2). Next, given an extension $A \map A_0$ in $\underline{C}$, 
$(\cX,\sD)/A \in F(A)$ and a finite $A_0$-module $M$, define
$$\cO_{(\cX,\sD)/A}(M)=T^2_{QG}(\cX_0/A_0,\cO_{\cX_0} \otimes_{A_0} M).$$
Given an extension $A' \map A$ of $A$ by $M$, define
$$o_{(\cX,\sD)/A}(A')=o^{QG}_{\cX/A}(A') \in \cO_{(\cX,\sD)/A}(M).$$
Then, by Theorem~\ref{thm-T^i_{QG}(X,D)}(1), these data give an obstruction theory for $F$.
We are now ready to prove the existence of (algebraic) formally smooth deformations of $F$ --- 
roughly, these provide the local patches of an algebraic stack.

\begin{thm} \label{thm-F_patches}
Given $(X,D) \in F(\bC)$, there exists an algebraic ring $R$ with a closed point $0 \in Spec(R)$ and 
$(\cX,\sD)/R \in F(R)$, formally smooth over $F$, such that $(\cX_0,\sD_0)=(X,D)$.
\end{thm}

\begin{proof}
The theorem is obtained by applying Theorem~\ref{thm-patches} to our functor $F$. We verify the conditions of the theorem in  Lemmas
\ref{lem-F-limit_preserving}, \ref{lem-F-S1,2}, \ref{lem-F-inv_lim}, and \ref{lem-F-D_and_O} below.
\end{proof}

\begin{lem} Let $(\cX,\sD)/A \in F(A)$, $A \in \underline{C}$ 
then $\omega^{[-N]}_{\cX/A}$ defines a projective embedding $\cX \inj \bP^M_A /A$ 
for $N \in \bN$ sufficiently large and divisible.  
\end{lem}
\begin{proof}
We have that $-K_{\cX/A}$ is $\bQ$-Cartier and relatively ample
using Lemma~\ref{lem-pfbc} and the base change property for $\omega_{\cX/A}^{[i]}$, $i \in \bZ$.
So, taking a sufficiently large and divisible multiple of $-K_{\cX/A}$, we obtain a projective embedding of $\cX/A$.
\end{proof}

\begin{rem}
Note that it is \emph{not} necessarily true that there exists $N$ such that for \emph{every} 
$(\cX,\sD)/A \in F$, the sheaf $\omega_{\cX/A}^{[-N]}$ defines a projective embedding of $\cX/S$.
The point is that $F$ defines a stack which is only \emph{locally} of finite type, i.e., we may require infinitely many patches.
However, if we restrict ourselves to smoothable pairs, i.e., if we consider the stack $\cM_d \subset F_d$, we can show 
that there is such an $N$ (using the bound on the index provided by Theorem~\ref{thm-ind}), 
and deduce that $\cM_d$ is of finite type (Theorem~\ref{thm-stack}).
\end{rem}

\begin{lem} \label{lem-Isom}
Given $(\cX,\sD)/A$ and  $(\cX',\sD')/A \in F(A)$, $A \in \underline{C}$, the functor
$$\Isom_{A}((\cX,\sD),(\cX',\sD')) \colon \underline{C} \backslash A \map (Sets)$$
$$ B \mapsto \{ \mbox{ Isomorphisms } \phi \colon (\cX,\sD) \otimes_A B \map (\cX',\sD') \otimes_A B \}$$
is represented by a quasiprojective scheme $\underline{\Isom}_{A}((\cX,\sD),(\cX',\sD'))/A$.
\end{lem}
\begin{proof}
We have a canonical polarisations $\omega_{\cX/A}^{[-N]}$ and $\omega_{\cX'/A'}^{[-N]}$ on $\cX/A$ and $\cX'/A'$
for some $N \in \bN$, and so by \cite{Gr} we know that the 
functor 
$$\Isom_{A}(\cX,\cX') \colon \underline{C} \backslash A \map (\emph{Sets})$$
$$ B \mapsto \{ \mbox{ Isomorphisms } \phi \colon \cX \otimes_A B \map \cX' \otimes_A B \}$$
is represented by a quasiprojective scheme $\underline{\Isom}_{A}(\cX,\cX')/A$.
It is then easy to construct $\underline{\Isom}_{A}((\cX,\sD),(\cX',\sD'))/A$ as a locally closed subscheme of 
$\underline{\Isom}_{A}(\cX,\cX')/A$.
\end{proof}

\begin{lem}(Open loci results) \label{lem-open_loci}
\begin{enumerate}
\item Let $\cX/S$ be a projective family of surfaces and $\sD \subset \cX$ a codimension $1$  closed subscheme, flat over $S$. 
Let $S' \subset S$ be the locus of points $s \in S$ such that
\begin{enumerate}
\item $\cX_s$ is CM, reduced and Gorenstein in codimension $1$.
\item $\sD_{s}=\sD_{(s)}$ and $\sD_{(s)}$ is Cartier in codimension $1$
\end{enumerate}
Then $S' \subset S$ is open, $\sD'=\sD \times_S S'$ is a relative Weil divisor on $\cX'=\cX \times_S S' / S'$, 
and $\omega_{\cX'/S'}$ corresponds to a relative Weil divisor $K_{\cX'/S'}$ on $\cX'/S'$.

\item Let $\cX/S$ be a projective family of CM reduced surfaces, Gorenstein in codimension $1$, 
and $\sD$ a relative effective Weil divisor on $\cX$.
Suppose that $\omega_{\cX/S}^{[i]}$ and $\cO_{\cX}(i\sD)$ commute with base change for all $i \in \bZ$.
Then the locus $S' \subset S$ where the geometric fibres of $\cX/S$ are quasistable of degree $d$ is open.
\end{enumerate} 
\end{lem}

\begin{proof}
(1) The locus where the fibres $\cX_s$ are CM is open by \cite{Mat}, p. 177, Corollary to Theorem 22.5.
$\cX_s$ is reduced iff it is regular in codimension $0$ --- this is an open condition.
Since $\omega_{\cX/S}$ commutes with base change, the requirement that $\cX_s$ is Gorenstein in codimension $1$ is open.
The condition $\sD_{s}=\sD_{(s)}$ is equivalent to requiring that the sheaf $\cO_{\cX}(-\sD) \otimes k(s)$ is $S_2$, which is an 
open condition, again by \cite{Mat} (note that $\cO_{\cX}(-\sD)$ is flat over $S$). 
Assuming this is satisfied, the natural map $\cO_{\cX}(-\sD) \otimes k(s) \map \cO_{\cX}(-\sD_{(s)})$ is an isomorphism, 
thus the condition $\sD_{(s)}$ Cartier in codimension $1$ is open. Hence $S' \subset S$ is open as required.
Using Lemma~\ref{lem-S_2} we deduce that $\sD'$ and $K_{\cX'/S'}$ are relative Weil divisors.

(2) It is enough to show the following:
Let $T$ be the spectrum  of a DVR with generic point $\eta$ and closed point $0=\Spec(k)$.
Let $(\cX,\sD)/T$ be a family of pairs such that $\omega_{\cX/T}^{[i]}$ and $\cO_{\cX}(i\sD)$ 
commute with base change for all $i \in \bZ$, and such that the special fibre $(X,D)/k$ is quasistable of degree $d$.
Then $(\cX_{\eta},\sD_{\eta})/ \eta$ is quasistable of degree $d$. 
Clearly $K_{\cX_{\eta}}+(\frac{3}{d}+\epsilon)\sD_{\eta}$
is ample, and it is also slc by Lemma~\ref{lem-iofa}. We are given $dK_X+3D \sim 0$ , we claim that
this implies $dK_{\cX_{\eta}}+3\sD_{\eta} \sim 0$ --- the essential point here is that $H^1(\cO_X)=0$, so $\Pic X$ is
discrete. To prove the claim, observe that 
$dK_{\cX}+3\sD$ is Cartier by Lemma~\ref{lem-pfbc} and Lemma~\ref{lem-rel_index}, and 
the restriction map $\Pic(\cX) \map \Pic(X)$ is an isomorphism (c.f. Proof of Lemma~\ref{lem-MV}(1)).
Thus $dK_{\cX}+3\sD \sim 0$, and restricting to the generic fibre we obtain our result.
We similiarly obtain $\frac{d}{3}K_{\cX_{\eta}}+\sD_{\eta} \sim 0$ in the case $3 \mid d$.
If $3 \nd d$ then we are given $3 \nd \ind D$, now $\ind \sD_{\eta} \mid  \ind \sD =\ind D$ using Lemma~\ref{lem-rel_index},
it follows that $3 \nd \ind \sD_{\eta}$. Given $H^1(\cO_X(D))=0$, it follows that $H^1(\cO_{\cX_{\eta}}(\sD_{\eta}))=0$
by semicontinuity (using the base change property for $\cO_{\cX}(\sD)$, note that $\cO_{\cX}(\sD)$ is flat over $T$ by
Lemma~\ref{lem-pfbc_flat}). Finally, $\chi(\cO_{X})=1$ is trivially an open condition. This completes the proof.
\end{proof}

\begin{lem} \label{lem-F-limit_preserving}
$F$ is limit preserving.
\end{lem}
\begin{proof}
Let $\{ A_i \}_{i \in I}$ be a direct system in $\underline{C}$ such that the limit 
$A = \stackrel{\lim}{\map} A_i$ lies in $\underline{C}$.
Suppose given $(\cX,\sD)/A \in F(A)$, we need to show that this is obtained from some $(\cX^i,\sD^i)/A_i \in F(A_i)$ by pullback.

We know that $\cX/A$ is projective, fix an embedding $\cX \inj \bP^N_A$.
Then $\cX /A$ is obtained by pullback from some projective flat family $\cX^i /A_i$ for some $i \in I$
(since $\Hilb_{P}(\bP^N/\bC)$ is of finite type, where $P$ denotes the Hilbert polynomial of the fibres of $\cX/A$).
By Lemma~\ref{lem-Emb_and_F} we have that $\sD \subset \cX$ is a codimension $1$ closed subscheme which is flat over $A$.
It follows (since $\Hilb_{Q}(\cX^i/A_i)$ is of finite type, where $Q$ denotes the Hilbert polynomial of the fibres of $\sD$)
that there exists $j \in I$ such that $\Spec A_j \map \Spec A_i$, and $(\cX^j,\sD^j)/A_j \in \sH ilb_Q(\cX^i/A_i)(A_j)$, 
such that $(\cX,\sD)/A$ is obtained from $(\cX^j,\sD^j)/A_j$ by pullback to $A$.
By Lemma~\ref{lem-open_loci}(1), we may also assume that each fibre $\cX^j_s$ is CM, reduced and Gorenstein in codimension $1$,
and that $\sD^j$ and $K_{\cX^j/A_j}$ are relative Weil divisors.

We now analyse the push forward and base change conditions.
Let $\cF^j$ be a coherent sheaf on an open subset $i^j \colon \cU^j \inj \cX^j/A_j$, flat over $A_j$,
and $\cF$, $i \colon \cU \inj \cX/A$ the corresponding objects obtained by pullback to $A$.
We are interested in the cases 
\begin{enumerate}
\item $\cF^j= \omega_{\cU^j/A_j}^{[n]}$ where $\cU^j \subset \cX^j$ is the locus where
$\omega_{\cX^j/A_j}$ is Cartier.
\item $\cF^j= \cO_{\cU^j}(n\sD^j)$, where $\cU^j \subset \cX^j$ is the locus where
$\sD^j$ is Cartier.
\end{enumerate}
Assume that the push forward of $\cF$ commutes with base change, and that 
$i^P_{\star}\cF^P$ is coherent for each $P \in \Spec A$.
We work locally at $P \in \Spec A$, say $P \mapsto Q \in \Spec A_j$.
The natural map $i_{\star}\cF \map i^P_{\star}\cF^P$ is surjective, so pick a finite set of elements of $\cF$ which generate
$i^P_{\star}\cF^P$ over $\cO_{\cX_P}$ (recall $i^P_{\star}\cF^P$ is assumed to be coherent).
Since $\cF=\cF^j \otimes_{A_j} A$, these are defined over some $A_k$ where $\Spec A_k \map \Spec A_j$.
Write $\cF^k$, $i^k \colon \cU^k \inj \cX^k/A$ for the objects obtained by pullback from $A_j$ to $A_k$,
and say $P \mapsto R \in \Spec A_k$.
Then, by construction, we have that $i^k_{\star}\cF^k \map i^R_{\star} F^R$ is surjective, 
hence by Lemma~\ref{lem-pfbc_closed_pts}(3)
the push forward of $\cF^k$ commutes with base change in a neighbourhood of $R$, which we may assume is $\Spec A_k$.

As in Remark~\ref{rem-finite}, we only need to consider a finite number of sheaves $\cF$.
Thus there exists $k \in I$ such that $\Spec A_k \map \Spec A_j$ and, writing   
$(\cX^k,\sD^k)/A_k$ for the pullback of $(\cX^j,\sD^j)/A_j$, the sheaves 
$\omega_{\cX^k/A_k}^{[n]}$ and $\cO_{\cX^k}(n\sD^k)$ commute with base change for all $n \in \bZ$.
                                                          
Finally, by Lemma~\ref{lem-open_loci}(2) we may assume that every geometric fibre of $(\cX^k,\sD^k)/A_k$ is quasistable.
Then $(\cX^k,\sD^k)/A_k$ is an element of $F(A_k)$. This completes the proof of the lemma.
\end{proof}

\begin{lem} \label{lem-F-S1,2}
F satisfies conditions (S1)(a),(b) and (S2).
\end{lem}

\begin{proof}
Let
\begin{eqnarray*}
\begin{array}{ccc}
   &      & B \\
   &      & \da \\
A' & \map & A
\end{array}
\end{eqnarray*}
be a diagram in $\underline{C}$ as in the statement of condition (S1)(a).
Write $B'=A' \times_A B$, then $B'$ is an infinitesimal extension of $B$.
We need to show that $\bar{F}_a(B') \map \bar{F}_a(A') \times \bar{F}_a(B)$ is surjective
for $a \in F(A)$. So, let $(\cX_{A'},\sD_{A'})/A' \in F(A')$ and $(\cX_{B},\sD_{B})/B \in F(B)$ be families which extend some
$(\cX_A,\sD_A)/A \in F(A)$. Define $(\cX_{B'},\sD_{B'})$ as follows: let $\us (\cX_{B'},\sD_{B'})= \us (\cX_B,\sD_B)$
( where $\us$ denotes the underlying topological spaces), and $\cO_{\cX_{B'}}=\cO_{\cX_{A'}} \times_{\cO_{\cX_A}} \cO_{\cX_B}$,
$\cO_{\sD_{B'}}=\cO_{\sD_{A'}} \times_{\cO_{\sD_A}} \cO_{\sD_B}$. Then $(\cX_{B'},\sD_{B'})/B'$ is a flat family of pairs
extending $(\cX_{A'},\sD_{A'})/A'$ and $(\cX_B,\sD_B)/B$ (to prove flatness, use \cite{Sch} p. 216 Lemma~3.4). Thus, by
Lemma~\ref{lem-Emb_and_F}, we only need to verify that $\cX_{B'}/B'$ is a $\bQ$-Gorenstein family to obtain
$(\cX_{B'},\sD_{B'})/B' \in F(B')$ as required. To see this, take an index one cover $\cZ_B \map \cX_B$, this gives an index one 
cover $\cZ_A \map \cX_A$ on restriction to $A$ by Proposition~\ref{prop-QG_def}, extend this to an index one cover 
$\cZ_{A'} \map \cX_{A'}$. Define $\cZ_{B'}$ by $\cO_{\cZ_{B'}}=\cO_{\cZ_{A'}} \times_{\cO_{\cZ_A}} \cO_{\cZ_B}$
as above, then $\cZ_{B'} \map \cX_{B'}$ is an index one cover extending $\cZ_B \map \cX_B$ 
and thus $\cX_{B'}/B'$ is a $\bQ$-Gorenstein family by Proposition~\ref{prop-QG_def} as required.

Suppose given $(\cX_0,\sD_0)/A_0 \in F(A_0)$, an extension $A \map A_0$ in $\underline{C}$,
a finite $A_0$-module $M$, and $(\cX,\sD)/A \in F_{(\cX_0,\sD_0)/A_0}(A)$.
The condition (S1)(b) states that the natural map 
$$D_{(\cX,\sD)/A}(M) \map D_{(\cX_0,\sD_0)/A_0}(M)$$
is an isomorphism. By Theorem~\ref{thm-T^i_{QG}(X,D)}(2), each side is naturally identified with 
$T^1_{QG}((\cX_0,\sD_0)/A_0,\cO_{\cX_0} \otimes_{A_0} M)$, so the map is an isomorphism as required.
Finally, the finiteness condition (S2) for $D_{(\cX_0,\sD_0)/A_0}(M)$ is obvious from the construction of the module 
$T^1_{QG}((\cX_0,\sD_0)/A_0,\cO_{\cX_0} \otimes_{A_0} M)$ using the properness of $\cX_0/A_0$.
\end{proof}

\begin{lem} \label{lem-F-inv_lim}
For $\hat{A}$ a complete local ring in $\underline{C}$,
the map  
$$\bar{F}(\hat{A}) \map \stackrel{\lim}{\leftarrow} \bar{F}(\hat{A}/m^n)$$
is bijective.  
\end{lem}

\begin{proof}
Write $A_n= \hat{A} /m^n$.
An element of $\stackrel{\lim}{\leftarrow} \bar{F}(A_n)$ is a sequence of compatible families 
$(\cX_n,\sD_n)/ A_n \in F(A_n)$. This defines a pair of formal schemes $(\frak{X},\frak{D})/ \Spf\hat{A}$,
where $\Spf$ denotes the formal spectrum. We have a line bundle $\sL$ on $\frak{X}$ such that the restriction to the 
special fibre is ample, e.g. $\sL= \stackrel{\lim}{\leftarrow} \omega_{\cX_n/A_n}^{[-N]}$ where $N=\ind \cX_0$.
By Grothendieck's Existence Theorem (\cite{EGA}, III.5.4.5), $(\frak{X},\frak{D})$ is the completion of a proper pair
$(\cX,\sD) / \hat{A}$ along the fibre $\cX_0$. We claim that $(\cX,\sD) / \hat{A} \in F(\hat{A})$. 
First observe that $\omega_{\cX/\hat{A}}^{[i]}$ and $\cO_{\cX}(\sD)^{[i]}$ commute with base change for all $i \in \bZ$ by 
Lemma~\ref{lem-pfbc_inv_lim}(3). Then, by Lemma~\ref{lem-open_loci}, the set of points 
$$\{ P \in \Spec \hat{A} \mid \mbox{ The geometric fibre of } (\cX,\sD)/\hat{A} \mbox{ over } P \mbox{ is quasistable }\}$$
is open. But it contains the closed point by assumption, hence it is the whole of $\Spec \hat{A}$.
Thus $(\cX,\sD) / \hat{A} \in F(\hat{A})$ as claimed, so the map 
$\bar{F}(\hat{A}) \map \stackrel{\lim}{\leftarrow} \bar{F}(\hat{A}/m^n)$ is surjective.

Now suppose given $(\cX,\sD)/ \hat{A}$ and $(\cX',\sD')/ \hat{A} \in F(\hat{A})$
which give the same element of $\stackrel{\lim}{\leftarrow} \bar{F}(\hat{A}/m^n)$, i.e., we have compatible isomorphisms
$$\phi_n \colon (\cX,\sD) \otimes_{\hat{A}} A_n \map (\cX',\sD') \otimes_{\hat{A}} A_n$$
for each $n \in \bN$. Equivalently, using Lemma~\ref{lem-Isom}, we have compatible maps
$$\Spec A_n \map \underline{\Isom}_{\hat{A}}((\cX,\sD),(\cX',\sD')),$$
and thus a map 
$$\Spec \hat{A} \map \underline{\Isom}_{\hat{A}}((\cX,\sD),(\cX',\sD')),$$
so there is an isomorphism $\phi \colon (\cX,\sD) \map (\cX',\sD')$ over $\hat{A}$ extending the $\phi_n$.
Thus the map $\bar{F}(\hat{A}) \map \stackrel{\lim}{\leftarrow} \bar{F}(\hat{A}/m^n)$ is injective.
\end{proof}

\begin{lem} \label{lem-F-D_and_O}
$D$ and $\cO$ satisfy the conditions \ref{cond-D_and_O} for algebraic $A,B \in \underline{C}$.
\end{lem}

\begin{proof}
(1) Let $A \map B$ in $\underline{C}$ be etale (in fact we only require $A \map B$ flat).
Given a flat family of schemes $\cZ/A$ and a coherent sheaf $\cG$ on $\cZ$ 
we have natural isomorphisms
$$\cT^i(\cZ/A,\cG) \otimes_A B \cong \cT^i(\cZ \otimes_A B / B, \cG \otimes_A B)$$
by \cite{LS}, p. 50, 2.3.2, using $A \map B$ flat.  
Given a $\bQ$-Gorenstein family of slc surfaces $\cX/A$ and a coherent sheaf $\cF$ on $\cX$, 
let $\cZ \map \cX/A$ be a local index one cover (a $\mu_N$ quotient, say), then 
$$\cT^i_{QG}(\cX/A,\cF)= (\pi_{\star}\cT^i(\cZ/A,\pi^{\star}\cF))^{\mu_N}$$
by Remark~\ref{rem-T^i_{QG}}. Thus we obtain 
$$\cT^i_{QG}(\cX/A,\cF) \otimes_A B \cong \cT^i_{QG}(\cX \otimes_A B / B, \cF \otimes_A B)$$
by applying $(\pi_{\star} \cdot)^{\mu_N}$ to the natural isomorphism above with $\cG=\pi^{\star}\cF$.
By \cite{Har}, p. 255, Proposition~9.3 we have natural isomorphisms
$$H^j(\cX,\sH) \otimes_A B \cong H^j(\cX \otimes_A B, \sH \otimes_A B)$$
for any quasi-coherent sheaf $\sH$ since $A \map B$ is flat.
Thus
\begin{equation}\tag{a}
T^i_{QG}(\cX/A,\cF) \otimes_A B \cong T^i_{QG}(\cX \otimes_A B / B, \cF \otimes_A B)
\end{equation}
for each $i$, using the spectral sequence of Remark~\ref{rem-QG_ss} and flatness of $A \map B$.
Also, given  coherent sheaves $\cF$, $\cG$ on $\cX$, we have 
\begin{equation} \tag{b}
\Hom_{\cX}(\cF,\cG) \otimes_A B \cong \Hom_{\cX \otimes_A B}(\cF \otimes_A B, \cG \otimes_A B)
\end{equation}
since $A \map B$ is flat (\cite{Mat}, p. 52, Theorem~7.11). 

Now suppose given  $A \map A_0$ an extension in $\underline{C}$, $(\cX_0,\sD_0)/A_0 \in F(A_0)$, 
$(\cX,\sD)/A \in F_{(\cX_0,\sD_0)/A_0}(A)$ and $M$ a finite $A_0$-module.  
Since
$$\cO_{(\cX,\sD)/A}(M)=T^2_{QG}(\cX_0/A_0, \cO_{\cX_0} \otimes_{A_0} M),$$
we see that $\cO$ commutes with etale localisation by (a) above. 
Using the exact sequence of Lemma~\ref{lem-T^i_{QG}(X,D)} together with (a) and (b) 
above we find that $D_{(\cX_0,\sD_0)/A_0}(M) = T^1_{QG}((\cX_0,\sD_0)/A_0, \cO_{\cX_0} \otimes_{A_0} M)$ 
commutes with etale localisation.

(2) We need to show that, given $a_0 \in F(A_0)$, $M$ a finite $A_0$-module, and $m \subset A$ a maximal ideal, we have
$$D_{a_0}(M) \otimes \hat{A_0} \cong \stackrel{\lim}{\leftarrow} D_{a_0}(M / m^nM).$$
By (1), we may assume that $A_0 = \hat{A_0}$, then we need to show
$$D_{a_0}(M) \cong \stackrel{\lim}{\leftarrow} D_{a_0}(M /m^n M),$$
that is,
$$\bar{F}_{a_0}(A_0+M) \cong \stackrel{\lim}{\leftarrow} \bar{F}_{a_0}(A_0+M /m^n M).$$
This is a result in the style of Lemma~\ref{lem-F-inv_lim}, and is proved in the same way.

(3) We may assume that $A_0$ is an integral domain, since we are given that $A_0$ is reduced, and $D$ and $\cO$ commute with etale
localisation. 

We first show that, for $A$ integral, and $\cX/A$ a family of schemes over $A$, 
there is an open affine subset $\Spec \tilde{A} \subset \Spec A$ such that
if $A \map B$ is a morphism in $\underline{C}$ which factors through $A \map \tilde{A}$, then the natural maps
$$ T^i(\cX/A,\cF) \otimes_A B \map T^i(\cX \otimes_A B/B, \cF \otimes_A B)$$  
are isomorphisms. Recall the construction of the $T^i$:
$$\cT^i(\cX/A,\cF)= \sH^i(\cHom(\sL_{\cdot},\cF)),$$
where $\sL_{\cdot}$ is a cotangent complex for $\cX/A$, and $T^i$ is obtained via the spectral sequence
$$E^{pq}_2=H^p(\cT^q) \Rightarrow T^{p+q}.$$
If $\sL_{\cdot}$ is a cotangent complex for $\cX/A$, then given $A \map B$, $\sL \otimes_A B$ is a cotangent complex for 
$\cX \otimes_A B /B$ (\cite{LS}, p. 47, 2.2.1(c)). Recall that, for a coherent sheaf $\cF$ on $\cX/A$ (of finite type),
where $A$ integral, 
there is a non-empty open affine subset $\Spec \tilde{A}$ of $\Spec A$ such that $\cF \otimes_A \tilde{A}$ is flat over $\tilde{A}$
(\cite{Mu}, Lecture 8, p. 57). 
It follows that, given coherent sheaves $\cF,\cG$ on $\cX/A$, there exists an open affine 
subset $\Spec \tilde{A} \subset \Spec A$ such that the natural map
$$\cHom(\cF,\cG) \otimes_A B \map \cHom(\cF \otimes_A B, \cG \otimes_A B)$$
is an isomorphism whenever $A \map B$ factors through $A \map \tilde{A}$ 
(compare the proof of \cite{Mat}, p. 52, Theorem~7.11). Thus, in our case, there is an open affine subset
$\Spec \tilde{A} \subset \Spec A$ such that the natural maps
$$\cHom(\sL_j,\cF) \otimes_A B  \map \cHom(\sL_j \otimes_A B, \cF \otimes_A B)$$
are isomorphisms for all $j$ and all $A \map \tilde{A} \map B$.
Similiarly, shrinking $\Spec \tilde{A}$ if necessary, we may assume the natural maps
$$\sH^q(\cHom(\sL_{\cdot},\cF)) \otimes_A B \map \sH^q(\cHom(\sL_{\cdot},\cF) \otimes_A B)$$
are isomorphisms for all $q$ and all $A \map \tilde{A} \map B$. So combining, we have isomorphisms
$$\cT^q(\cX/A,\cF) \otimes_A B \map \cT^q(\cX \otimes_A B /B, \cF \otimes_A B)$$
for all $q$ and all $A  \map \tilde{A} \map B$.
Next, we may assume that the natural maps
$$H^p(\cT^q) \otimes_A B \map H^p(\cT^q \otimes_A B)$$
are isomorphisms for all $p$, $q$ and $A \map \tilde{A} \map B$, 
by the theory of cohomology and base change for projective morphisms
(see e.g. \cite{Har}, III.12, in particular p. 288, Corollary~12.9), thus
$$H^p(\cT^q(\cX/A,\cF) ) \otimes_A B \cong H^p(\cT^q(\cX \otimes_A B /B, \cF \otimes_A B)).$$
Finally, we may assume that we obtain
$$T^{p+q}(\cX/A,\cF) \otimes_A B  \cong T^{p+q}(\cX \otimes_A B / B, \cF \otimes_A B)$$
using the local-to-global spectral sequence, by imposing certain flatness requirements.

The same argument shows that the  $T^i_{QG}(\cX/A,\cF)$ commute with base changes $A \map \tilde{A} \map B$ for some $\tilde{A}$,
where $\cX/A$ is a family of slc surfaces.
We use the local identification 
$$\cT^i_{QG}(\cX/A,\cF)=(\pi_{\star}\cT^i(\cZ/A,\pi^{\star}\cF))^{\mu_N}$$
where $\pi \colon \cZ \map \cX$ is a local index one cover which is a $\mu_N$ quotient.

We use the exact sequence of Lemma~\ref{lem-T^i_{QG}(X,D)} to show that $T^1_{QG}((\cX,\sD)/A,\cF)$
commutes with base changes $A \map \tilde{A} \map B$ for some $\tilde{A}$.
We may assume that the modules $T^0(\cX/A,\cF)$, $\Hom_{\cO_X}(\cI_{\sD},\cO_{\sD} \otimes_A M)$ and 
$T^1_{QG}(\cX/A,\cF)$ commute with base changes $A \map B$ factoring through some $A \map \tilde{A}$ by the above,
and if we assume some flatness conditions we obtain that $T^1_{QG}((\cX,\sD)/A,\cF)$ commutes with base changes 
$A \map \tilde{A} \map B$ using the 5-lemma. 

Thus, given $A_0$ integral, $M$ a finite $A_0$-module, and $A \map A_0$ an extension, 
$D_{(\cX_0,\sD_0)/A_0}(M)=T^1_{QG}((\cX_0,\sD_0)/A_0,\cO_{\cX_0}\otimes_{A_0} M)$ and 
$\cO_{(\cX,\sD)/A}(M)= T^2_{QG}(\cX_0/A_0,\cO_{\cX_0} \otimes_{A_0} M)$ commute with base changes
$A_0 \map \tilde{A_0} \map B_0$ for some open affine subset $\Spec \tilde{A_0} \subset \Spec A_0$.
This completes the proof.
\end{proof}

\begin{thm} \label{thm-F_stack}
$F$ is a Deligne-Mumford algebraic stack, locally of finite type over $\bC$.
\end{thm}
\begin{proof}
First, it is clear that $F$ defines a stack. The only non-trivial point is that etale descent data for $F$ are effective ---
this follows since for each $S \in \underline{Sch}$ and $(\cX,\sD)/S \in F(S)$ 
we have a canonical polarisation $\omega_{\cX/S}^{[-N]}$ for some $N \in \bN$.

By Theorem~\ref{thm-F_patches}, we can construct a scheme $V$ locally of finite type over $\bC$ and $v \in F(V)$
such that the morphism of stacks $\underline{V} \map F$ is smooth and surjective.
The construction is as follows: For every $(X,D) \in F(\bC)$, let $R$ be an algebraic ring with a closed point $0 \in \Spec R$ and 
$(\cX,\sD)/R \in F(R)$, formally smooth over $F$, such that $(\cX_0,\sD_0)=(X,D)$. 
Let $V$ be the disjoint union of the schemes $\Spec R$, and $v \in F(V)$ the union of the families $(\cX,\sD)/R$.

We have that the diagonal $F \map F \times F$ is representable, quasicompact and separated, i.e., for every $S$ and 
every pair of elements 
$(\cX,\sD)/S$, $(\cX',\sD')/S \in F(S)$, the functor 
$$\Isom_{S}((\cX,\sD),(\cX',\sD')) \colon (\emph{Schemes}/S) \map (\emph{Sets})$$
$$ T \mapsto \{ \mbox{ Isomorphisms } \phi \colon (\cX,\sD) \times_S T \map (\cX',\sD') \times_S T \}$$
is represented by a scheme which is quasicompact and separated over $S$. 
For $\Isom_S((\cX,\sD),(\cX',\sD'))$ is represented by a 
scheme $\underline{\Isom}_S((\cX,\sD), (\cX',\sD'))/S$, quasiprojective over $S$, by Lemma~\ref{lem-Isom}. 
Thus $F$ is an algebraic stack, locally of finite type over $\bC$.

Finally, we claim that $F$ is a Deligne-Mumford algebraic stack, i.e., there exists a scheme $V$
locally of finite type over $\bC$ and $v \in F(V)$ such that the morphism of stacks 
$\underline{V} \map F$ is $\emph{etale}$ and surjective.
By \cite{DM}, p. 104, Theorem~4.21, it is enough to show that the diagonal map $F \map F \times F$ is unramified, i.e.,
for every $S$ and $(\cX,\sD)/S$, $(\cX',\sD')/S \in F(S)$, the scheme $\underline{\Isom}_S((\cX,\sD), (\cX',\sD'))/S$ is 
unramified over $S$.
We may assume that $S= \Spec k$, $k$ algebraically closed.
We know that $\Aut((X,D)/k)$ is finite for $(X,D)/k \in F(k)$ by \cite{Ii}, since $K_X+(\frac{3}{d}+\epsilon)D$ is slc and ample. 
It follows that $(X,D)/k$ has no infinitesimal automorphisms
(since $\ch k =0$) and hence $F \map F \times F$ is unramified as required.  
\end{proof}

\begin{thm} \label{thm-stack}
$\cM_d$ is a Deligne-Mumford algebraic stack, of finite type over $\bC$.
\end{thm}
\begin{proof}
It is immediate that $\cM_d \subset F_d$ is a Deligne-Mumford algebraic stack, locally of finite type over $\bC$ by 
Theorem~\ref{thm-F_stack} and Definition~\ref{defn-smoothable_deformation}.
It remains to show that it is of finite type. By Theorem~\ref{thm-ind} there exists $N(d) \in \bN$ such that
$N(d)K_X$ is Cartier for every stable pair $(X,D)/k$ of degree $d$. 
By smoothability, the Hilbert polynomials of $X$ and $D$ with respect to the polarisation $-N(d)K_X$ are
the same as the Hilbert polynomials of $\bP^2$ and $C \subset \bP^2$ a curve of degree $d$ with respect to the
polarisation $-N(d)K_{\bP^2}$. 
Also, by \cite{Ko5}, Theorem~2.1.2, there exists $M(d) \in \bN$ such that $-M(d)K_X$ is very ample and has no higher cohomology 
for all $(X,D)$.
Thus there is a Hilbert scheme $H$ of pairs $(Y,B) \inj \bP^K$ 
such that every stable pair of degree $d$ occurs in the universal family. 
Now, let $S$ be a local affine patch of $\cM_d$ with universal family $(\cX,\sD)/S$.
Pick an projective embedding $(\cX,\sD)/S \inj \bP^K_S$ defined by $-M(d)K_{\cX/S}$.
We obtain a map $S \map H$, and moreover taking all possible such embeddings we obtain a map 
$$\phi \colon \underline{\Aut}(\bP^K) \times S =PGL(K+1) \times S \map H.$$ 
We perform this construction everywhere locally on $\cM_d$.
The (set-theoretic) image of the maps $\phi$ is the set $I$ of points $[(Y,B) \inj \bP^K]$ of $H$ 
such that $(Y,B)$ is a stable pair of degree $d$ embedded by $-M(d)K_Y$. We claim that $I$ is a locally closed subset of $H$  
and moreover that the images of the maps $\phi$ give an open cover of $I$. Assuming this for the moment, we deduce that 
there exists a finite open subcover of $I$, and hence a finite affine open cover of $\cM_d$. 
Thus $\cM_d$ is of finite type as required. 

We now prove our claim above.
We first apply Lemma~\ref{lem-open_loci}(1) to deduce that there exists an open locus $I^1 \subset H$ containing $I$
such that for all $[(Y,B) \inj \bP^K] \in I^1$, $Y$ is CM, reduced and Gorenstein in codimension $1$,
and $B$ is Cartier in codimension $1$. Moreover, writing $(\cY^1,\cB^1)/I^1 \inj \bP^K_{I^1}$ for the universal family,
$K_{\cY^1/I^1}$ and $\cB^1$ are relative Weil divisors on $\cY^1/I^1$. Next, we apply Theorem~\ref{thm-pfbc_stratification}, 
to deduce that there exists a locally closed subset $I^2 \subset I^1$ containing $I$ such that, writing
$(\cY^2,\cB^2) \inj \bP^K_{I^2}$ for the universal family, 
$\omega_{\cY^2/I^2}^{[i]}$ and $\cO_{\cY^2}(i\cB)$ commute with base changes $T \map I^2$ with $T$ reduced,  for all $i \in \bZ$.
There now exists an open subset $I^3 \subset I^2$ containing $I$ such that $(Y,B)$ is quasistable for all 
$[(Y,B) \inj \bP^K] \in I^3$ by Lemma~\ref{lem-open_loci}(2). 
Let $I^4$ be the closure in $I^3$ of the locus where $Y \cong \bP^2$,
then $I^4$ contains $I$ and $(Y,B)$ is a stable pair of degree $d$ for all $[(Y,B) \inj \bP^K] \in I^4$. Finally, $I \subset I^4$
is the locus where $(Y,B) \inj \bP^K$ is defined by $-M(d)K_{Y}$ --- this is an locally closed condition. Thus $I$ is a locally 
closed subset of $H$. Now let $S$ be a local patch of $\cM_d$ with universal family $(\cX,\sD)/S$ and 
$\phi \colon \underline{\Aut}(\bP^K) \times S \map H$ a map as constructed above, we claim that the image of $\phi$ is 
an open subset of $I$. Let $T$ be the spectrum of a DVR with closed point $0$ and generic point $\eta$, and 
$h \colon T \map I$ a morphism (where $I \subset H$ is given its reduced structure) such that $h(0)$ lies in the image of $\phi$.
It is enough to show that $h(\eta)$ also lies in the image of $\phi$. Writing $(\cY_T,\cB_T)/T$ for the pullback of the universal
family of $H$ to $T$, we have $(\cY_T,\cB_T)/T \in \cM_d(T)$ (note that the base change conditions are satisfied by construction).
We deduce our result by the versality of $(\cX,\sD)/S$.
\end{proof}

\section{Classification of the type B surfaces} \label{B}

We provide a characterisation of the slc del Pezzo surfaces of type B which have a $\bQ$-Gorenstein smoothing to $\bP^2$.

\begin{thm} \label{thm-B}
Suppose $X$ is a surface of type B. Then $X$ admits a $\bQ$-Gorenstein smoothing to $\bP^2$ iff
\begin{enumerate}
\item $X$ has singularities of types $\frac{1}{n^2}(1,na-1)$, $(a,n)=1$,
 and $(xy=0) \subset \bA^3_{x,y,z} / \frac{1}{r}(1,-1,a)$, $(a,r)=1$.
 Moreover there are at most two non-normal singularities with index $r>1$.
\item $K_X^2=9$.
\item $\rho(X_1)=\rho(X_2)=1$ or $\{ \rho(X_1),\rho(X_2) \}=\{ 1,2 \}$.
\end{enumerate}
\end{thm} 

We use the $\bQ$-Gorenstein deformation theory developed in Section~\ref{QG}.
Write $A_1=k[t]/(t^2)$.
We use the short hand $\cT^i_X=\cT^i(X/k,\cO_X \otimes_k A_1)$, and similiarly 
$T^i_X$, $\cT^i_{QG,X}$ and $T^i_{QG,X}$. 
The following theorem describes how to calculate the $T^i_X$. 

\begin{thm}
Let $X$ be a scheme over $k$. We have
\begin{enumerate}
\item $\cT^0_X = \cHom(\Omega_X,\cO_X)=T_X$, the tangent sheaf of $X$.
\item $\cT^1_X$ can be calculated locally from an embedding $X \inj \bA$ of $X$ in an affine space via the exact sequence
	$T_{\bA}|_X \map \cN_{X/\bA} \map \cT^1_X \map 0$. In particular $\cT^1_X =0$ if $X$ is smooth.
\item $\cT^2_X =0$ if $X$ is a local complete intersection.
\end{enumerate}
\end{thm}

\begin{proof}
All this follows from Theorem~\ref{thm-calc_cT} apart from the exact sequence of (2) --- see \cite{LS}, p. 51, 2.3.6.
\end{proof}

We also recall the following results from Section~\ref{QG} --- if $\pi \colon Z \map X$ is a local index one cover of 
$X$, a $\mu_N$ quotient say, then 
$$\cT^i_{QG,X}=(\pi_{\star}\cT^i_Z)^{\mu_N}$$
and moreover $\cT^0_{QG,X}=\cT^0_X$.

\begin{lem} \label{lem-T^1}
Let $X$ be a surface with two normal irreducible components $X_1$, $X_2$ meeting in a smooth double curve $\Delta$.
Suppose $X$ has only singularities of the form $(xy=0) \subset \frac{1}{r}(1,-1,a)$ at $\Delta$. Then, in a neighbourhood of
$\Delta$, $$\cT^1_{QG,X} \cong \cO_{\Delta}(\Delta_1|_{\Delta} + \Delta_2 |_{\Delta}).$$
Here $\Delta_i$ is the double curve on $X_i$, and we calculate $\Delta_i |_{\Delta}$ by moving $\Delta_i$ on $X_i$, and 
restricting to $\Delta$ --- we obtain a $\bQ$-divisor on $\Delta$ which is well defined modulo linear equivalence.
We have that $\Delta_1|_{\Delta} + \Delta_2|_{\Delta}$ is a $\bZ$-divisor on $\Delta$.
In particular, $\cT^1_{QG,X}$ is a line bundle on $\Delta$ of degree $\Delta_1^2+\Delta_2^2$.
\end{lem} 
\begin{proof}
Let $P \in \Delta$, we first work locally analytically at $P$. 
Suppose that $X$ is normal crossing at $P$, i.e., 
$$X \cong (xy=0) \subset \bA^3_{x,y,z}.$$
Then, using the exact sequence $T_{\bA^3}|_X \map \cN_{X/\bA^3} \map \cT^1_X \map 0$, we deduce that $\cT^1_X$ is a line 
bundle on $\Delta$. Moreover we obtain a natural isomorphism 
$\cT^1_X \cong \cO_{\Delta_1}(\Delta_1) \otimes_{\cO_{\Delta}} \cO_{\Delta_2}(\Delta_2)$ --- for we have
$$\cT^1_X \cong \cN_{X/\bA^3}|_{\Delta} \cong \cO_{\Delta}(X_2) \otimes \cO_{\Delta}(X_1)
\cong \cO_{\Delta_1}(\Delta_1) \otimes_{\cO_{\Delta}} \cO_{\Delta_2}(\Delta_2).$$
Suppose now that $X$ has a singularity $(xy=0) \subset \frac{1}{r}(1,-1,a)$ at $P$. 
Let $\pi \colon Z \map X$ be the canonical cover, then $Z$ is normal crossing.
Write $Z_1$, $Z_2$ for the inverse images of $X_1$, $X_2$, and $\Delta'$, $\Delta_1'$, $\Delta_2'$ 
for the inverse images of $\Delta$, $\Delta_1$, $\Delta_2$ on $Z$, $Z_1$, $Z_2$.  
We have a natural isomorphism $\cT^1_Z \cong \cO_{\Delta_1'}(\Delta_1') \otimes_{\cO_{\Delta}} \cO_{\Delta_2'}(\Delta_2')$, 
and $\cT^1_{QG,X}= (\pi_{\star}\cT^1_Z)^{\mu_r}$.

Patching together this local information, we obtain that globally $\cT^1_{QG,X} \cong \cO_{\Delta}(\Delta_1|_{\Delta} +\Delta_2|_{\Delta})$, where $\Delta_i |_{\Delta}$
is defined as detailed above.
\end{proof}

\begin{lem} \label{lem-Delta^2}
Suppose $X$ is a surface of type B which satisfies the conditions (1),(2),(3) of Theorem~\ref{thm-B}.
Then we have $\Delta_1^2+\Delta_2^2= 3 - (\rho(X_1)+\rho(X_2))$.
\end{lem}
\begin{proof}
Let $\tilde{X_i} \map X_i$ for $i=1$ and $2$ be the minimal resolutions of the components of $X$.
Then 
$$K_{\tilde{X_i}}^2+\rho(\tilde{X_i})=10$$
for $i=1$ and $2$ by Noether's formula, and 
$$K_{\tilde{X_1}}^2+K_{\tilde{X_2}}^2+\rho(\tilde{X_1})+\rho(\tilde{X_2})=
K_{X_1}^2+K_{X_2}^2+\rho(X_1)+\rho(X_2)+ 4 \sum (1-\frac{1}{r_j}),$$
where the $r_j$ are the indices of the non-Gorenstein singularities of $X$ at $\Delta$ (see the proof of 
Theorem~\ref{thm-B*}). Thus
$$K_{X_1}^2+K_{X_2}^2 = 20 -(\rho(X_1)+\rho(X_2))-4 \sum (1-\frac{1}{r_j}).$$
We also have 
$$(K_{X_1}+\Delta_1)^2+(K_{X_2}+\Delta_2)^2 = 9$$
and
$$K_{X_i}\Delta_i+\Delta_i^2=-2+ \sum (1-\frac{1}{r_j})$$
for $i=1$ and $2$. Solving for $\Delta_1^2+\Delta_2^2$ we obtain $\Delta_1^2+\Delta_2^2= 3 - (\rho(X_1)+\rho(X_2))$.
\end{proof}

\begin{lem} \label{lem-H^2(T)}
Suppose $X$ is a surface of type B which satisfies the conditions (1),(2),(3) of Theorem~\ref{thm-B}. Then $H^2(T_X)=0$.
\end{lem}

\begin{proof}
We have an exact sequence
$$ 0 \map \cO_{X_1}(-\Delta_1) \oplus \cO_{X_2}(-\Delta_2) \map \cO_X \map \cO_{\Delta} \map 0.$$
Apply the functor $\cHom_{\cO_X}(\Omega_X, \cdot)$ --- we obtain an exact sequence
$$ 0 \map T_{X_1}(-\Delta_1) \oplus T_{X_2}(-\Delta_2) \map T_X \map \cHom_{\cO_{\Delta}}(\Omega_X|_{\Delta}, \cO_{\Delta}).$$
Thus we have an inclusion $T_{X_1}(-\Delta_1) \oplus T_{X_2}(-\Delta_2) \inj T_X$ with cokernel supported on $\Delta$.
It follows that $H^2(T_{X_1}(-\Delta_1)) \oplus H^2(T_{X_2}(-\Delta_2)) \map H^2(T_X)$ is surjective.
So it is enough to show that $H^2(T_{X_i}(-\Delta_i))=0$ for $i=1,2$.

Write $(Y,C)$ for a component $(X_i,\Delta_i)$.
By Serre duality, 
$$H^2(T_{Y}(-C)) \cong \Hom(T_{Y}(-C),\cO_{Y}(K_{Y}))^{\vee}
= \Hom(T_{Y},\cO_{Y}(K_{Y}+C))^{\vee}.$$
We claim that $\cO_Y(-K_{Y}-C)$ has a nonzero global section. Assuming this, we have
$$\Hom(T_{Y},\cO_{Y}(K_{Y}+C)) \inj \Hom(T_{Y},\cO_{Y})=H^0(\Omega_{Y}^{\vee\vee}).$$
Now, letting $\pi : \tilde{Y} \map Y$ be the minimal resolution, we have 
$\Omega_{Y}^{\vee\vee}=\pi_{\star}\Omega_{\tilde{Y}}$ since $Y$ has only quotient singularities 
(\cite{St1}, Lemma~1.11). Thus $h^0(\Omega_{Y}^{\vee\vee})=h^0(\Omega_{\tilde{Y}})=h^1(\cO_{\tilde{Y}})=0$, since 
$\tilde{Y}$ is rational by Theorem~\ref{thm-min-res}. So $H^2(T_{Y}(-C))=0$ as required.

It remains to show that $\cO_Y(-K_{Y}-C)$ has a nonzero global section.
We have an exact sequence
$$0 \map \cO_Y(-K_Y-2C) \map \cO_Y(-K_Y-C) \map \cO_C(-K_Y-C) \map 0.$$
Now $h^1(\cO_Y(-K_Y-2C))=h^1(\cO_Y(2K_Y+2C))=0$ by Serre duality and Kodaira vanishing (recall that $Y$ is log terminal
and $-(K_Y+C)$ is ample). So we are done if $\cO_C(-K_Y-C)$ has a nonzero global section.
A local calculation shows that $\cO_C(-K_Y-C) \cong \cO_C(-K_C-S)$, where $S$ is the sum of the singular points of $Y$ 
lying on $C$. Now $C$ is isomorphic to $\bP^1$, and there are at most $2$ singular points of $Y$ on $C$ by assumption, 
thus $\deg(-K_C-S) \ge 0$ and $\cO_C(-K_Y-C)$ has a nonzero global section as required.
\end{proof}

\begin{lem} \label{lem-le2}
Let $X$ be a surface of type B which admits a $\bQ$-Gorenstein smoothing to $\bP^2$.
Then there are at most two non-normal points of $X$ of index greater than $1$.
\end{lem}

\begin{proof}  
For a smoothable surface $X$ of type B, 
the non-Gorenstein singularities at $\Delta$ are of type $(xy=0) \subset \frac{1}{r}(1,-1,a)$. 
Moreover, if we have $k$ such singularities with indices $r_1, \ldots, r_k$ then $\sum (1-\frac{1}{r_i}) < 2$, 
in particular $k \le 3$. We show that if $X$ smoothes to $\bP^2$, then $k \le 2$. We assume $k=3$ and obtain a contradiction.
Let $\cX/T$ be a $\bQ$-Gorenstein smoothing of $X$ over the spectrum $T$ of a DVR with generic point $\eta$, such that 
$\cX_{\eta} \cong \bP^2_{\eta}$. Write $X=X_1 \cup X_2$, then there are two cases :
\begin{enumerate}
\item $\rho(X_1)=1$ and $\rho(X_2)=2$, $\rho(\cX/T)=2$, $\cX$ is $\bQ$-factorial.
\item  $\rho(X_1)=\rho(X_2)=1$, $\rho(\cX/T)=1$, $\cX$ is not $\bQ$-factorial.
\end{enumerate}
Here we use Corollary~\ref{cor-rho}, also $\Cl(\cX) \cong \bZ^{\oplus 2}$ by Lemma~\ref{lem-Cl}(1) and 
we can calculate $\rho(\cX/T)$ using Lemma~\ref{lem-MV}(2).

Suppose first $X$ is of type (1). We claim that we can contract the divisor $X_1 \subset \cX$ to obtain a relative minimal model
$\bar{\cX}/T$. We have $\cX$ canonical by Lemma~\ref{lem-iofa}, so we can use the relative MMP theory.
It's enough to show that the double curve $\Delta$ of $X$ generates an extremal ray in $\overline{NE}(\cX/T)$
(note that $-K_{\cX}$ is relatively ample so certainly $K_{\cX} \Delta <0$).
Let $\cX$ have singularity $((xy=t^n) \subset \bA^4_{x,y,z,t}) \cong \frac{1}{n}(1,-1) \times \bA^1_z$ generically along $\Delta$.
Then 
$$\Delta^2_{X_1}+\Delta^2_{X_2}=n X_1 X_2 (X_2+X_1) =0.$$
Now $\Delta^2_{X_1} > 0$ since $\rho(X_1)=1$, thus $\Delta^2_{X_2} < 0$ and $\Delta$ generates an extremal ray on $X_2$.
By the exact sequence
$$0 \map \bZ \Delta \map N_1(X_1) \oplus N_1(X_2) \map N_1(\cX/T) \map 0$$
(c.f. Lemma~\ref{lem-MV}), we see that $\Delta$ generates an extremal ray on $\cX$.
Hence we can contract $X_1$ to obtain a $\bQ$-factorial family of surfaces $\bar{\cX}/T$ with generic fibre $\bP^2_{\eta}$.
The special fibre $\bar{X}$ is obtained from $X_2$ by contracting the double curve $\Delta$. Thus $\bar{X}$ has a 
log-terminal singularity whose minimal resolution has exceptional locus a tree of rational curves with one fork.
But this singularity is not smoothable, a contradiction.

Now suppose that $X$ is of type (2). Then $\cX$ is not $\bQ$-factorial. Thus we have a point $P \in \Delta$ such that,
locally analytically at $P$,
$$\cX \cong (xy+t^ng(z^r,t)=0) \subset \bA^4_{x,y,z,t} / \frac{1}{r}(1,-1,a,0),$$
where $t \nd g(z^r,t)$, $n \in \bN$. 
Working locally at $P \in \cX$, we construct a small partial resolution $\pi \colon \tilde{\cX} \map \cX$.
Assume first that $r=1$. Let $\pi \colon \tilde{\cX} \map \cX$ be the blowup of $(x=g=0) \subset \cX$.
Writing $u=\frac{g}{x}$ and $v=\frac{x}{g}$, we have two affine pieces of $\tilde{\cX}$ as follows:
$$(vy + t^n=0) \subset \bA^4_{v,y,z,t}$$
$$((y + t^n u =0, xu = g) \subset \bA^5_{x,y,z,t,u}) \cong ((xu=g) \subset \bA^4_{x,u,z,t})$$
Thus $\tilde{X}$ is normal crossing at the strict transform $\tilde{\Delta}$ of $\Delta$ and $\tilde{\cX}$
has a $\frac{1}{n}(1,-1) \times \bA^1_z$ singularity along $\tilde{\Delta}$. The only possible singularity of 
$\tilde{X} \subset \tilde{\cX}$ away from $\Delta$ is a $cA$ singularity at $(x,u,z,t)=0$.
Write $X_1=(x=0)$, $X_2=(y=0)$, and $\tilde{X_1}$, $\tilde{X_2}$ for the strict transforms.
Then $\tilde{X_1} \map X_1$ is an isomorphism, and $\tilde{X_2} \map X_2$ contracts a $\bP^1$ to the smooth point $P \in X_1$.
In the case $r>1$ we obtain $\pi \colon \tilde{\cX} \map \cX$ as the $\mu_r$ quotient of the construction above.

Now work globally on $\tilde{\cX}$. 
Note that $\tilde{X}$ is slt, $\tilde{\cX}$ is canonical, and $K_{\tilde{\cX}}=\pi^{\star}K_{\cX}$.
We first claim that $\tilde{\cX}/T$ is projective. For it is easy to see that the special fibre $\tilde{X}$ is projective, 
and the restriction map $\Pic \tilde{\cX} \map \Pic \tilde{X}$ is an isomorphism (using $-K_{\tilde{X}}$ nef and big, c.f.
proof of Lemma~\ref{lem-MV}(1)). Moreover $\tilde{\cX}$ is $\bQ$ factorial (c.f. Corollary~\ref{cor-rho}). 
We now obtain a contradiction as for case (1) above.  
\end{proof}

\begin{proof}[Proof of Theorem~\ref{thm-B}]
First suppose that $X$ is a surface of type B that admits a $\bQ$-Gorenstein smoothing to $\bP^2$.
The condition (1) on the singularities of $X$ is satisified by Theorem~\ref{thm-smoothability}, Lemma~\ref{lem-T_1} and 
Lemma~\ref{lem-le2}.
Let $\cX/T$ be a $\bQ$-Gorenstein smoothing of $X/k$, where $T$ is the spectrum of a DVR with generic point $\eta$,
such that $\cX_{\eta} \cong \bP^2_{\eta}$.
We have $K_X^2=K_{\cX_{\eta}}^2=9$ since $K_{\cX}$ is $\bQ$-Cartier, so (2) holds.
Finally, (3) is satisfied by Corollary \ref{cor-rho}. Thus the conditions are necessary.

Now let $X$ be a surface of type B which satisfies the conditions (1),(2) and (3). 
We first construct a first order $\bQ$-Gorenstein deformation $\cX^1 / A_1$ of $X/k$, where $A_1=k[t]/(t^2)$.
We then show that this extends to a $\bQ$-Gorenstein smoothing of $X$ over a DVR.
Recall that there is a natural isomorphism $\Def^{QG}_{X/k}(A_1) \cong T^1_X$. 
Now, by the local-to-global spectral sequence, we have an exact sequence
$$0 \map H^1(\cT^0_X) \map T^1_{QG,X} \map H^0(\cT^1_{QG,X}) \map H^2(\cT^0_X).$$
By Lemma~\ref{lem-H^2(T)} we have $H^2(\cT^0_X)=H^2(T_X)=0$, so 
$$0 \map H^1(\cT^0_X) \map T^1_{QG,X} \map H^0(\cT^1_{QG,X}) \map 0$$
is exact. 
We next specify an element of $H^0(\cT^1_{QG,X})$, and describe the local first order deformations of $X$ it determines.
We can then lift this to an element of $T^1_{QG,X}$ defining a global $\bQ$-Gorenstein first order deformation
$\cX^1/A_1$ of $X/k$. 

The sheaf $\cT^1_{QG,X}$ is supported on the singular locus of $X$.
So, to define a global section of $\cT^1_{QG,X}$, we have to define a section in a neighbourhood 
of each connected component of the singular locus of $X$.
 An isolated singularity $P$ of $X$ is of the form 
$\frac{1}{n^2}(1,na-1)$ by assumption. Write $Z \map X$ for the local index one cover of $X$, then   
\begin{eqnarray*}
\begin{array}{ccccc}
Z   & \cong &(xy+z^n=0) & \subset & \bA^3_{x,y,z} \\ 
\da &       & \da	& 	  & \da \\
X   & \cong &(xy+z^n=0) & \subset & \bA^3_{x,y,z} / \frac{1}{n}(1,-1,a).
\end{array}
\end{eqnarray*}
We define a local first order $\bQ$-Gorenstein deformation of $X$ as follows:
$$(xy+z^n+t=0) \subset  \bA^3_{x,y,z} / \frac{1}{n}(1,-1,a) \times \Spec A_1.$$
This corresponds to a local section of ${\cT^1_X}'$.
Now consider $\Delta$, the locus of non-isolated singularities of $X$. In a neighbourhood of $\Delta$,
$\cT^1_{QG,X}$ is a line bundle on $\Delta$
of degree $\Delta_1^2 + \Delta_2^2$ by Lemma~\ref{lem-T_1}. Moreover $\Delta \cong \bP^1$ and 
$\Delta_1^2+\Delta_2^2= 3 - (\rho(X_1)+\rho(X_2)) \ge 0$ by Lemma~\ref{lem-Delta^2}. So there exists a  
section $s$ of $\cT^1_{QG,X}$ at $\Delta$,
with reduced divisor of zeroes missing the non-Gorenstein points of $X$, i.e.,
$$0 \map \cO_{\Delta} \stackrel{s}{\map} \cT^1_{QG,X} \map \oplus_{i=1}^k k(P_i) \map 0$$
where $P_1,\ldots,P_k$ are distinct points of index one on $X$. Then $s$ defines local first order deformations of $X$ at the 
points $P$ of $\Delta$ of the following forms. If $P \neq P_1, \ldots, P_k$ and 
$X \cong (xy=0) \subset \bA^3_{x,y,z} / \frac{1}{r}(1,-1,a)$ at $P$, some $r \ge 1$, then we have
$$(xy+ tu=0) \subset \bA^3_{x,y,z} / \frac{1}{r}(1,-1,a) \times \Spec A_1,$$
for some unit $u \in k[[z^r]]$.
If $P=P_i$ for some $i$ then  $X \cong (xy=0) \subset \bA^3_{x,y,z}$ at $P$ and we have
$$(xy+ tzu=0) \subset \bA^3_{x,y,z} \times \Spec A_1,$$
for some unit $u \in k[[z]]$.
 
Observe that $X$ has unobstructed $\bQ$-Gorenstein deformations. 
For the local index one covers of $X$ are local complete intersections, so $\cT^2_{QG,X}=0$ 
(thus $X$ has unobstructed $\bQ$-Gorenstein deformations locally).
Also, we have $H^1(\cT^1_{QG,X})=0$, since $\cT^1_{QG,X}$ is an invertible sheaf on $\Delta \cong \bP^1$ of non-negative degree, and 
$H^2(T_X)=0$ by Lemma~\ref{lem-H^2(T)}. Thus $T^2_{QG,X}=0$ by 
the local-to-global spectral sequence, so $X$ has unobstructed $\bQ$-Gorenstein deformations as claimed. 
Hence we can extend the deformation $\cX^1/A_1$ to a sequence of compatible 
$\bQ$-Gorenstein deformations $\cX^n/A_n$, where $A_n = k[t]/(t^{n+1})$.
They determine a formal scheme $\hat{\cX}/ \Spf \hat{A}$, where $A=k[t]_{(t)}$, $\hat{A}$ is the completion of $A$, and $\Spf$
denotes the formal spectrum. By the Grothendieck Existence Theorem, $\hat{\cX}/ \Spf \hat{A}$ is the completion of a projective 
$\bQ$-Gorenstein family $\cX / \Spec \hat{A}$ (c.f. Proof of Lemma~\ref{lem-F-inv_lim}).
We claim that this is a smoothing of $X$. For, from the explicit descriptions of the local 
first order deformations above, we see that $\cX$ has only singularities of types $\frac{1}{n}(1,-1,a)$ and 
$(xy+zt=0) \subset \bA^4_{x,y,z,t}$, hence $\cX$ is smooth over the generic point $\eta$ of $\Spec \hat{A}$ as required.
The geometric generic fibre $\cX_{\bar{\eta}}$ is a smooth del Pezzo surface with $K^2=9$, hence 
$\cX_{\bar{\eta}} \cong \bP^2_{\bar{\eta}}$.
Thus $\cX' / \Spec \hat{A}$ is a $\bQ$-Gorenstein smoothing of $X$ to $\bP^2$ as required.
\end{proof}

We include below an important result we have proved along the way. 

\begin{thm} \label{thm-slt_unobs}
Let $X$ be an slt del Pezzo surface that admits a $\bQ$-Gorenstein smoothing to $\bP^2$.
Then $X$ has unobstructed $\bQ$-Gorenstein deformations.
\end{thm} 

\begin{proof}
$X$ is either a Manetti surface or of type B. We proved the result in the type B case in the proof of Theorem~\ref{thm-B}.
If $X$ is a Manetti surface then the local index one covers of $X$ are local complete intersections so $\cT^2_{QG,X}=0$, 
$H^1(\cT^1_{QG,X})=0$ since $\cT^1_{QG,X}$ has $0$ dimensional support, 
and $H^2(\cT^0_X)=0$ by \cite{Ma}, p. 113, Proof of Theorem~21.
Thus $T^2_{QG,X}=0$, so $X$ has unobstructed $\bQ$-Gorenstein deformations as required.
\end{proof}

\section{Simplifications in the case $3 \nd d$} \label{3_divis}
We state and prove two major simplifications we obtain in the case $3 \nd d$.
First, if $(X,D)$ is a stable pair of degree $d$ then $X$ is slt. 
Second, the stack $\cM_d$ is smooth.

\begin{thm} \label{thm-simp}
Let $(X,D)$ be a stable pair of degree $d$.
Suppose $3 \nd d$. Then $X$ is slt. So $X$ is either a Manetti surface 
(that is, normal, with log terminal singularities) or of type B.
In particular $X$ has $1$ or $2$ components.
\end{thm}

\begin{proof}
Write $(Y,C)$ for a component of $(X^{\nu},\Delta^{\nu})$, and let $D_Y$ denote $D \mid_Y$. We need to show that $(Y,C)$ is 
log terminal.
First we show that there are no singularities of types $(\frac{1}{r}(1,a), 2 \Delta)$ or $(D,\Delta)$.
For suppose $\Gamma$ is a component of $C$ which passes through such a point.
Then there is at most one other singularity of $(Y,C)$ on $\Gamma$, of type $(\frac{1}{s}(1,b),\Delta)$.
We calculate  $(K_Y+C)\Gamma = -2 +1+(1-\frac{1}{s})=-\frac{1}{s}$ (we allow $s=1$ --- this is the case where there are no other
singularities on $C$).
Now $K_Y+C+\frac{3}{d}D_Y \equiv 0$, so $D_Y \Gamma = \frac{d}{3} \frac{1}{s}$.
But $D_Y$ misses the strictly log canonical point of $(Y,C)$, since $K_X+(\frac{3}{d}+\epsilon)D$ is slc,
thus $sD_Y$ is Cartier near $\Gamma$. So $D_Y \Gamma \in \frac{1}{s} \bZ$, a contradiction.

It only remains to show that no normal strictly log canonical singularities can occur.
Otherwise, $X$ is a normal surface, so is an elliptic cone by Theorem~\ref{thm-normal-case}. Let $l$ be a ruling,
$\pi : \tilde{X} \map X$ the minimal resolution of $X$ and $E$ the exceptional curve. Thus
$\tilde{X}$ is a ruled surface over an elliptic curve, $E$ is the negative section, and the strict transform $l'$ of $l$
is a ruling. We calculate $K_X l = \pi^{\star}K_X l'=(K_{\tilde{X}} +E)l'=-2+1=-1$.
Now $K_X+\frac{3}{d}D \equiv 0$ implies $Dl=\frac{d}{3}$. 
But $D$ misses the singularity of $X$ since $K_X+(\frac{3}{d}+\epsilon)D$
is log canonical, so $D$ is Cartier. Thus $Dl \in \bZ$, a contradiction.
\end{proof}

\begin{rem}
In particular we have a classification of the surfaces $X$ occurring by Theorem~\ref{thm-Manetti} and Theorem~\ref{thm-B}.
\end{rem}

\begin{rem}
If $3 \mid d$ the situation is much more complicated. 
For example in the case $d=6$ the maximum number of components of a surface
$X$ is $18$. See Example~\ref{ex-typeC}.
\end{rem}

\begin{thm}
$\cM_d$ is smooth for $3 \nd d$.
\end{thm}
\begin{proof}
First, by Theorem~\ref{thm-slt_unobs} and Theorem~\ref{thm-simp} above, 
for $(X,D)/k \in \cM_d(k)$, $X/k$ has unobstructed $\bQ$-Gorenstein deformations.
Second, by Theorem~\ref{thm-F_to_Def^QG_smooth}, the map of functors $\bar{\cM}_{d,(X,D)/k} \map \Def^{QG}_{X/k}$
is smooth. Note that $\bar{\cM}_{d,(X,D)/k}=\bar{F}_{d,(X,D)/k}$ in this case, i.e., every $\bQ$-Gorenstein deformation
of $(X,D)/k$ is smoothable. Thus $\cM_d$ is smooth as required.
\end{proof}

\begin{rem}
$\cM_d$ is not smooth if $3 \mid d$ --- see Example~\ref{ex-elliptic_cone_defns}.
\end{rem}

\section{Classification of stable pairs for $d=4,5$ and $6$} \label{small_d}

We set out below the complete classification of stable pairs of degrees 4 and 5.
In degree 6 we give a complete list of candidates for the surfaces that occur, 
however we have yet to establish which surfaces of types C and D
are smoothable. I explain how the classification was obtained in Section~\ref{classn}.

\begin{notn}
\begin{enumerate}
\item Given a weighted projective space $\bP$, let $H$ denote a section of $\cO_{\bP}(1)$. We also write $\sim kH$
to denote a general section of $\cO_{\bP}(k)$.\\
\item Given a rational ruled surface $\bF_n$, let $A$ and $B$ denote a fibre and the negative section respectively.
\end{enumerate}
\end{notn}

When we refer to the singularities of $(X,D)$ below, we mean the points of $X$ where $(X,D)$ is not normal crossing, i.e.,
where we do not have $X$ normal crossing and $D^{\nu}+\Delta^{\nu} \subset X^{\nu}$ normal crossing. 

\subsection{$d=4$} Surfaces $X$:
\begin{eqnarray*}
\renewcommand{\arraystretch}{1.5}
\begin{array}{|l|l|l|} \hline
\mbox{Surface}      	& \mbox{Double curve}	& \mbox{Singularities} \\ \hline
\bP^2				&			&			\\
\bP(1,1,4)			&			& \frac{1}{4}(1,1)	\\
\bP(1,1,2) \cup \bP(1,1,2)	& H,H 			& (xy=0) \subset \frac{1}{2}(1,1,1) \\ \hline
\end{array}
\end{eqnarray*}

Singularities of $(X,D)$:
\begin{eqnarray*}
\renewcommand{\arraystretch}{1.5}
\begin{array}{|l|l|} \hline
X      	& D	\\ \hline
\bA^2_{x,y}				& (y^2+x^3=0) \\
\frac{1}{4}(1,1)		&   0 \\
(xy=0) \subset \frac{1}{2}(1,1,1) & 0\\ \hline
\end{array}
\end{eqnarray*}

\subsection{$d=5$} Surfaces $X$:
\begin{eqnarray*}
\renewcommand{\arraystretch}{1.5}
\begin{array}{|l|l|l|} \hline
\mbox{Surface}   	   	& \mbox{Double curve}	& \mbox{Singularities} \\ \hline
\bP^2					&			&			\\
\bP(1,1,4)				&			& \frac{1}{4}(1,1)	\\
X_{26} \subset \bP(1,2,13,25)		& 			& \frac{1}{25}(1,4)	\\
\bP(1,4,25)				&			& \frac{1}{4}(1,1), \frac{1}{25}(1,4) \\
\bP(1,1,2) \cup \bP(1,1,2)		& H,H 			& (xy=0) \subset \frac{1}{2}(1,1,1) \\
\bP(1,1,5) \cup (X_6 \subset \bP(1,2,3,5))& H, \sim 2H		& (xy=0) \subset \frac{1}{5}(1,-1,1) \\
\bP(1,1,5) \cup \bP(1,4,5)		& H, \sim 4H		& \frac{1}{4}(1,1), (xy=0) \subset \frac{1}{5}(1,-1,1) \\ \hline
\end{array}
\end{eqnarray*}

Singularities of $(X,D)$:
\begin{eqnarray*}
\renewcommand{\arraystretch}{1.5}
\begin{array}{|l|l|} \hline
X      		& D	\\ \hline
\bA^2_{x,y}	& (y^2+x^n=0) \mbox{ for } 3 \le n \le 9\\
\bA^2_{x,y}	& (x(y^2+x^n)=0) \mbox{ for } n=2,3 \\
\bA^2_{x,y} / \frac{1}{4}(1,1) & (y^2+x^n=0) \mbox{ for } n=2,6 \\
(xy=0) \subset \bA^3_{x,y,z} / \frac{1}{2}(1,1,1) & (z=0) \\ 
\frac{1}{25}(1,4) & 0 \\
(xy=0) \subset \frac{1}{5}(1,-1,1) & 0 \\ \hline
\end{array}
\end{eqnarray*}

\pagebreak

\subsection{$d=6$}
We have the following cases for the surfaces $X$:
\subsubsection{Types A and B}

\begin{eqnarray*}
\renewcommand{\arraystretch}{1.5}
\begin{array}{|l|l|l|} \hline
\mbox{Surface}  		& \mbox{Double curve} 	& \mbox{Singularities}	\\ \hline
\bP^2					& 			&			\\		
\bP(1,1,4)				&			& \frac{1}{4}(1,1)	\\
\mbox{Elliptic cone, degree $9$}	&			& \mbox{simple elliptic singularity} \\
\bP^2 \cup \bF_1				& H,B 			&			\\		
\bP^2 \cup \bF_4				& \sim 2H, B		&			\\
\bP(1,1,4) \cup \bF_4			& \sim 4H, B		& \frac{1}{4}(1,1)	\\
\bP(1,1,2) \cup \bP(1,1,2)		& H,H			& (xy=0) \subset \frac{1}{2}(1,1,1) \\ \hline
\end{array}
\end{eqnarray*}

\subsubsection{Types C and D} \label{sss-typeC&D}
In this case we only have a partial solution -- we give a complete list of candidates, but we have yet to establish which are smoothable. 

We give below a list of possible components $(Y,C)$ of the surfaces $X$, 
i.e., $Y$ is a component of the normalisation of $X$ and $C$ is the inverse image of the double curve. 
We then glue these together, following the instructions in Theorem~\ref{thm-glueing} and Theorem~\ref{thm-smoothability}, 
to recover $X$. We have the following constraints:
\begin{enumerate}
\item  $K_X^2=9$. Equivalently, $\sum (K_Y+C)^2 =9$, where $(Y,C)$ runs over the components of $X$.
\item  $\sum \rho(Y) \le 2V$ in case C, $\sum \rho(Y) \le 2V-1$ in case D, where $V$ is the number of components of $X$ --- see 
Corollary~\ref{cor-rho}.
\end{enumerate}  

Recall that the components $(Y,C)$ are of two types III and IV. We tabulate these separately.

\begin{notn} For a singularity of type $(D,\Delta)$, the exceptional locus of the minimal resolution consists of a chain 
$F_1, \ldots, F_k$ of curves together with two -2 curves meeting $F_k$. The strict transform of $\Delta$ meets $F_1$.
We use the sequence $-F_1^2,\ldots,-F_k^2$ of self-intersections to describe the $(D,\Delta)$ singularities occurring below.
\end{notn}

Type III:
\begin{eqnarray*} 
\renewcommand{\arraystretch}{1.5}
\begin{array}{|l|l|l|l|c|ll|} \hline
\mbox{Case} & \mbox{Parameter}	&\rho(Y)& (K_Y+C)^2	& \multicolumn{3}{c|}{\mbox{Singularities}} \\
	&			&	&		&\mbox{Number of}& \multicolumn{2}{l|}{(\frac{1}{r}(1,a),2\Delta) \mbox{ sing.}} \\
	&			&	&		&\mbox{$(\frac{1}{2}(1,1),\Delta)$'s}& r 	& a\\ \hline		
1. 	& n \ge 1		& 1	& n		& 		& n 	& 1 \\
2. 	& n \ge 1		& 1	& n-\frac{1}{2}	& 1		& 2n-1 	& 2 \\
3.	& n \ge 2		& 1 	& n-1		& 2		&4(n-1) & 2n-1 \\
4.	& n \ge 0		& 2	& n+2		& 		&	&\\
5.	& n \ge 1		& 2 	& n+\frac{3}{2}	& 1		& 2 	& 1 \\
6.	& n \ge 2		& 2	& n -\frac{1}{2}& 1		& n 	& 1 \\
7.	& n \ge	2		& 2 	& n-1 		& 2	        &2n-1 	& 2\\
8.	& 2			& 2 	& \frac{7}{2} 	& 1		&	&\\
9.	& 2			& 2 	& 4		& 2		&	&\\
10.	& 2			& 2 	& 3		& 2 		& 3 	& 1 \\
11.	& n \ge 1		& 3 	& n		& 2		& 2 	& 1 \\
12.	& n \ge 2		& 3 	& n+\frac{1}{2}	& 1		&	&\\
13.	& n \ge 2		& 3	& n+\frac{3}{2}	& 1		&	&\\
14.	& n \ge 5		& 3     & n-1		& 2		& n 	& 1\\ \hline
\end{array}
\end{eqnarray*}
\smallskip
\newpage

Type IV: 
\begin{eqnarray*}
\renewcommand{\arraystretch}{1.5}
\begin{array}{|l|l|l|l|c|l|} \hline
\mbox{Case} & \mbox{Parameter} &\rho(Y)& (K_Y+C)^2	& \multicolumn{2}{c|}{\mbox{Singularities}} \\
	&		&	&		&\mbox{Number of}	& (D,\Delta) \mbox{ sing.} \\
	&		&	&		&\mbox{$(\frac{1}{2}(1,1),\Delta)$'s}& -F_1^2,\ldots,-F_k^2 \\ \hline
1.	& n \ge 2	& 1	& n-1		& 0		& n,2,2 \\
2.	& n \ge 2	& 1	& n-\frac{3}{2}	& 1		& 2,n,2,2 \\
3.	& n \ge 2	& 2	& n+1		& 0		& 2,2 \\
4.	& n \ge 2	& 2	& n-\frac{3}{2}	& 1		& n,2,2 \\
5.	& n \ge 2	& 2	& n-1		& 0		& n,2 \\
6.	& n \ge 3	& 2	& n-\frac{3}{2}	& 1		& 2,n,2 \\
7.	& 2		& 2	& \frac{7}{2}	& 1		& 2 \\
8.	& 2		& 2	& \frac{5}{2}	& 1		& 3,2 \\ \hline
\end{array}
\end{eqnarray*}

\begin{rem} In almost all cases the parameter is just the \emph{weight} of the surface $Y$,
namely the maximal $w$ such that there exists a birational morphism $\tilde{Y} \map \bF_w$,
where $\tilde{Y}$ is the minimal resolution of $Y$. This is only false when $n=1$ in type III 
cases  1,2,10 and 12
-- then  the weight is undefined in case 1 (because $Y=\bP^2$) and equals $2$ in the other cases.
\end{rem}

We describe the components of type III concretely below, using toric language.
Note that we can also give an explicit description of the type IV cases --- we describe $Y$ by expressing its minimal resolution 
$\tilde{Y}$  as a blowup of $\bF_n$. This is omitted here.

\begin{notn}
Let $\Bl_{(m,n)}S \map S$ denote the weighted blowup of a smooth point of a surface $S$ 
with weights $(m,n)$ with respect to some
local analytic coordinates. Unless otherwise stated, we assume that the point and choice of coordinates is general. 
We write $E$ for the exceptional curve and use primes to 
denote strict transforms of curves. If we refer to $C'$, it is assumed that $C$ passes through the centre $P$ of the blowup, e.g.
if $S=\bP$, $H'$ denotes the strict transform of a general section of $\cO_{\bP}(1)$ through $P$.

Let $\bF_{n-\frac{1}{2}}$ $(n \ge 1)$ 
denote the surface obtained from $\bF_n$ in the following way: First perform a sequence of two blowups 
away from the negative section to obtain a degenerate fibre which is a chain of curves of self intersections $-2,-1,-2$. Then 
contract the two -2 curves. Thus $\bF_{n-\frac{1}{2}}$ has one double fibre
with two $\frac{1}{2}(1,1)$ singularities on it, and a negative section with square $-(n-\frac{1}{2})$. We let $A$ and $B$ denote the
negative section and fibre as usual. We write $\frac{1}{2}A$ for the double fibre with its reduced structure.
\end{notn}

Type III (toric descriptions):
\begin{eqnarray*}
\renewcommand{\arraystretch}{1.5}
\begin{array}{|l|l|l|} \hline
\mbox{Case}	&\mbox{Surface component} 	& \mbox{Double curve} \\ \hline
1. 	&\bP(1,1,n)		  	&  H_1+H_2\\
2. 	&\bP(1,2,2n-1)			&  H+ \sim 2H\\
3.	&\bP(1,1,2(n-1)) / \mu_2		&  H_1 / \mu_2 +H_2 / \mu_2\\	
4.	&\bF_n	 			& A+B\\
5.      &\bF_{n-\frac{1}{2}}		& \frac{1}{2}A+B\\
6.	&\Bl_{(1,2)}\bP(1,1,n)		& H'_1+H_2\\
7.	&\Bl_{(1,2)}\bP(1,2,2n-1)	& H+(\sim 2H)'\\	
8.	&\Bl \bP(1,1,2)			& H'+E\\
9.	&\Bl_{(1,2)}\bP(1,1,2), \wt(H)=2	& H'+E\\
10.	&\Bl_{(2,3)}\bP(1,1,2)		& H'+E\\
11,n>1.	&\Bl_{(1,2)}\bF_{n-\frac{3}{2}}, \wt(B)=1 & \frac{1}{2}A+B'\\		
11,n=1.	&\Bl_{(1,2)}\bF_{\frac{1}{2}}    & \frac{1}{2}A + (\sim B+\frac{1}{2}A)'\\   
12.	&\Bl_{(1,2)}\bF_{n-1}, \wt(B)=1	& A+B'\\
13.	&\Bl_{(1,2)}\bF_n		& A'+B\\			
14.	&\Bl_{(1,2)}^2\bP(1,1,n)		& H'_1+H'_2\\ \hline		
\end{array} 
\end{eqnarray*}

\section{Derivation of the classification for $d=4,5$ and $6$}
\label{classn}

We aim to classify the surfaces $X$ for $d=4,5,6$. We provide a brief overview of our method below. 
We first obtain a list of candidates satisfying certain combinatorial conditions.
We then check which of these candidates actually occur --- the main point is to establish the smoothability of $X$.
This last step is incomplete in the case $d=6$.	

\begin{notn} \label{notn-tildeY}
Let $(X,D)$ be a stable pair of degree $d$ and $(Y,C)$ a component of $(X^{\nu},\Delta^{\nu})$.
Let $\pi : \tilde{Y} \map Y$ be the minimal resolution of $Y$ and define a $\bQ$-divisor $\tilde{C}$ 
on $\tilde{Y}$ by $K_{\tilde{Y}}+\tilde{C}= \pi^{\star}(K_Y+C)$, $\pi_{\star}\tilde{C} =C$
(as in Notation~\ref{notn-cpt}).
\end{notn}

We first give a list of candidates for the pairs $(Y,C)$.
We then glue these together as instructed in Theorem~\ref{thm-glueing} to obtain our list of candidates for $X$.  
We do not work directly with the pairs $(Y,C)$ --- instead we consider the pairs $(\tilde{Y},\tilde{C})$.
Note that it is immediate to recover $(Y,C)$ from $(\tilde{Y},\tilde{C})$ --- we just contract the components of $\tilde{C}$ 
where $-(K_{\tilde{Y}}+\tilde{C})$ is zero. For $-(K_Y+C)$ is ample and $-(K_{\tilde{Y}}+\tilde{C})=-\pi^{\star}(K_Y+C)$ by
definition, moreover $\Exc(\pi) \subset \Supp \tilde{C}$ by Lemma~\ref{lem-exc_locus} below. 
Note that we may assume $\tilde{Y}$ is rational, since otherwise $Y$ is an elliptic cone and $X=Y$.

We begin our classification of the $(Y,C)$ by giving a list of possible singularities.
We deduce the possible configurations 
of rational curves making up the $\bQ$-divisor $\tilde{C}$, together with multiplicities.
We obtain a finite list of possible slt singularities on $(Y,C)$ by using the conditions 
$dK_X+3D \sim 0$ and $K_X+(\frac{3}{d}+\epsilon)D$ slc to bound the index of $X$ 
(c.f. the proof of Theorem~\ref{thm-ind-explicit}). 
We are content to note that the strictly log canonical singularities of $(Y,C)$ are of two types
$(\frac{1}{r}(1,a),2\Delta)$ and $(D,\Delta)$. 
The corresponding configurations of components of
$\tilde{C}$ are easily understood --- working locally over the singular point $P \in Y$ we have:
\begin{enumerate}
\item $(\frac{1}{r}(1,a),2\Delta)$ : $\tilde{C}=C'_1+E_1+ \cdots+E_k+C'_2$, where $E_1,\ldots,E_k$ is a chain of
rational curves, and $C'_1$ and $C'_2$ intersect $E_1$ and $E_k$ respectively. 
\item $(D,\Delta)$ : $\tilde{C}=C'+F_1+\cdots+F_l+\frac{1}{2}G_1+\frac{1}{2}G_2$, where $F_1,\ldots,F_l$ is a chain of rational 
curves, $G_1$ and $G_2$ are $-2$-curves meeting $F_l$, and $C'$ intersects $F_1$.  
\end{enumerate}
Note that, if $3 \nd d$, $(Y,C)$ is log terminal by Theorem~\ref{thm-simp}, so these cases do not occur.

\begin{notn} \label{notn-bir_ruling}
Assume $Y$ is rational and $Y \not \cong \bP^2$, then there exists a birational morphism $\tilde{Y} \map \bF_w$, where $w \ge 0$
--- fix one such morphism $\phi$, with $w$ maximal. Then $\phi$ is an isomorphism over the negative section $B$ of $\bF_w$.
Let $p$ denote the composite $\tilde{Y} \map \bF_w \map \bP^1$.
\end{notn}

We next analyse all the ways that we can fit the possible $\tilde{C}$ configurations into a surface
$\tilde{Y}$ with a birational ruling $\tilde{Y} \map \bF_w \map \bP^1$ as above. 
We require that $-(K_{\tilde{Y}}+\tilde{C})$ is nef and big, and positive outside $\tilde{C}$
--- this imposes strong restrictions on curves in degenerate fibres which are not contained in 
$\Supp \tilde{C}$ (see Lemma~\ref{lem-Gamma}). 
Also, we have $B' \subset \Supp \tilde{C}$ if $w \ge 2$, where $B'$ denotes the strict transform of 
$B \subset \bF_w$ under $\phi \colon \tilde{Y} \map \bF_w$. 
We obtain the candidates for $d=4$ and $5$ in this way in the Proof of Theorem~\ref{thm-d=4,5}. 
We omit the derivation of the list of candidates for $d=6$. It is rather different in style, since $(Y,C)$ may have 
strictly log canonical singularities. In particular we do \emph{not} obtain a finite list of possible singularities 
as our first step. We begin by giving a list of possible degenerate fibres of $(\tilde{Y},\tilde{C})$. 
Then, with some work, we obtain a complete list of candidates. In each case $\Supp \tilde{C}$ is a collection of components of 
fibres together with $B'$ and at most one other horizontal component.

\begin{lem} \label{lem-exc_locus}
With the notation as above,  the $\bQ$-divisor  $\tilde{C}$ is effective and $\Supp(\tilde{C}) = C' \cup \Exc(\pi)$.
\end{lem}

\begin{proof}
We have $K_{\tilde{Y}}+C' = \pi^{\star}(K_Y+C) + \sum a_iE_i$ with $a_i \le 0$, and $\tilde{C} = C' + \sum (-a_i)E_i$ by 
definition (compare the proof of Proposition~\ref{prop-setup}).
We need to show $a_i < 0$ for all $i$. If $a_j=0$ for some $j$, then $a_k=0$ for all $E_k$ in the same connected component of
$\Exc(\pi)$ using $K_{\tilde{Y}}+\tilde{C}$ $\pi$-nef.
It follows that this connected component contracts to a canonical singularity of $(Y,C)$. 
This must be a Du Val singularity of $Y$ with $C=0$ locally. But then we have a normal log-terminal singularity on $X$ which is 
not of type $T_1$, a contradiction.
\end{proof}

\begin{lem} \label{lem-index_small_d}
Let $(X,D)$ be a stable pair of degree $d=4,5$ or $6$.
For $P \in X$  we have 
\begin{enumerate}
\item d=4 : $\ind_P K_X = 1,2$ or $4$.
\item d=5 : $\ind_P K_X = 1,2$ or $5$.
\item d=6 : $\ind_P K_X= 1$ or $2$.
\end{enumerate} 
\end{lem}

\begin{proof}
See the Proof of Theorem~\ref{thm-ind-explicit} --- a finer analysis in the cases $d=4,5$ and $6$ gives our result.
\end{proof}

\begin{lem} \label{lem-C_cpts}
Let $d=4$ or $5$. Then the connected components of $\tilde{C}$ are chains of smooth rational curves.
We have the following possibilities for the multiplicities and self-intersections of the components:
\begin{enumerate}
\item {$d=4$} 
\begin{eqnarray*}
\renewcommand{\arraystretch}{1.5}
\begin{array}{llll}
Case	& \mbox{Multiplicities} 		& \mbox{Self-Intersections} 	& \mbox{Image in $Y$}\\
(1)	& \frac{1}{2}    			& -4 				& \frac{1}{4}(1,1) \\
(2)	& 1,\frac{1}{2}  			& ?,-2 				& \bP^1 \ni \frac{1}{2}(1,1) \\
(3)	& \frac{3}{4},\frac{1}{2},\frac{1}{4}	& -6,-2,-2 			& \frac{1}{16}(1,3) \\
(4)	& \frac{1}{2},1,\frac{3}{4}		& -2,?,-4			& \bP^1 \ni \frac{1}{2}(1,1),\frac{1}{4}(1,1) \\
(5) 	& \frac{1}{2},1,\frac{3}{4},\frac{1}{2},\frac{1}{4} & -2,?,-2,-2,-2     & \bP^1 \ni \frac{1}{2}(1,1),\frac{1}{4}(1,3) 
\end{array}
\end{eqnarray*}\\

\item {$d=5$} 
\begin{eqnarray*}
\renewcommand{\arraystretch}{1.5}
\begin{array}{llll}
Case	& \mbox{Multiplicities} 		& \mbox{Self-Intersections} 	& \mbox{Image in $Y$}\\
(1)	& \frac{1}{2}    			& -4 				& \frac{1}{4}(1,1) \\
(2)	& 1,\frac{1}{2}  			& ?,-2 				& \bP^1 \ni \frac{1}{2}(1,1) \\
(3)	& \frac{4}{5},\frac{3}{5},\frac{2}{5},\frac{1}{5}	& -7,-2,-2,-2 	& \frac{1}{25}(1,4) \\
(4)	& \frac{3}{5},\frac{4}{5},\frac{2}{5}	& -3,-5,-2			& \frac{1}{25}(1,9) \\
(5) 	& 1,\frac{4}{5} 			& ?,-5    			& \bP^1 \ni \frac{1}{5}(1,1) \\
(6)	& 1,\frac{4}{5},\frac{3}{5}		& ?,-2,-3			& \bP^1 \ni \frac{1}{5}(1,2) \\
(7)	& 1,\frac{4}{5},\frac{2}{5}		& ?,-3,-2			& \bP^1 \ni \frac{1}{5}(1,3) \\
(8)     & 1,\frac{4}{5},\frac{3}{5},\frac{2}{5},\frac{1}{5} & ?,-2,-2,-2,-2	& \bP^1 \ni \frac{1}{5}(1,4)  
\end{array}
\end{eqnarray*}\\
\end{enumerate}
Here we write $\bP^1 \ni \frac{1}{r}(1,a)$ to denote a smooth rational curve $C$ on $Y$ passing through a singularity of 
type $\frac{1}{r}(1,a)$, such that locally analytically $(Y,C) \cong (\bA^2_{x,y} / \frac{1}{r}(1,a),(x=0))$.
\end{lem}

\begin{proof}
X is slt for $3 \nd d$ by Theorem~\ref{thm-simp}, thus $(Y,C)$ is log terminal. 
Recall that the log terminal singularities of $(Y,C)$ are of types $(\frac{1}{n^2}(1,na-1),0)$ and $(\frac{1}{r}(1,a),\Delta)$. 
Lemma~\ref{lem-index_small_d} states which values of the indices $n$ and $r$ are possible.
It only remains to determine what combinations of singularities can lie on $C \subset Y$ for $C \neq \emptyset$.
Suppose $C \neq \emptyset$, and let $(Y,C)$ have singularities of indices $r_1,\ldots,r_k$ at $C$.
Then $C \cong \bP^1$ and $-(K_Y+C)C=2 - \sum (1-\frac{1}{r_i}) >0$ (compare Theorem~\ref{thm-cpts}(II)).
Moreover, $-(K_Y+C)C$ is 3-divisible by Lemma~\ref{lem-Gamma3} below. Using the restrictions on the indices in 
Lemma~\ref{lem-index_small_d}, we obtain the solutions $k=1$, $r_1=2$, or $k=2$, $r_1=2$, $r_2=4$ for $d=4$,
and $k=1$, $r_1=2$, or $k=1$, $r_1=5$ for $d=5$.
Combining our results we obtain the lists of possible connected components of $\tilde{C}$ above.
\end{proof}

\begin{lem} \label{lem-Gamma3} Suppose $3 \nd d$. 
Let $\Gamma \subset Y$ be a curve. Then $-(K_Y+C)\Gamma$ is 
divisible by 3 (i.e. writing $-(K_Y+C)\Gamma=\frac{a}{b}$, $a,b \in \bN$, $(a,b)=1$, we have $3 \mid a$).
Equivalently, let $\Gamma' \subset \tilde{Y}$ be a curve that is not $\pi$-exceptional, 
then $-(K_{\tilde{Y}}+\tilde{C})\Gamma'$ is divisible by 3.
\end{lem}

\begin{proof}
We use the relation $d(K_Y+C)+3D_Y \sim 0$.
Thus $-(K_Y+C)\Gamma = \frac{3}{d} D_Y \Gamma$, so it's enough to show that
the index of $D_Y$ is not 3-divisible. See the proof of Lemma~\ref{lem-stable_and_quasistable}.
\end{proof}

\begin{lem} \label{lem-Gamma} Suppose $\Gamma$ is a curve contained in a degenerate fibre of $p$ which is not a 
component of $\tilde{C}$.
Then $\Gamma$ is a -1 curve and $\tilde{C} \Gamma <1$.
\end{lem}

\begin{proof} We know $\Gamma$ is not $\pi$-exceptional by Lemma~\ref{lem-exc_locus}, hence $(K_{\tilde{Y}}+\tilde{C})\Gamma <0$.
But $(K_{\tilde{Y}}+\tilde{C})\Gamma = -2 -\Gamma^2 +\tilde{C}\Gamma \ge -2 -\Gamma^2$.
Hence $\Gamma^2=-1$ and  $\tilde{C} \Gamma <1$ as claimed.
\end{proof}

\begin{notn}
We refer to such a curve $\Gamma$ as a $\gamma$-curve
\end{notn}

\begin{cor} \label{cor-Gamma}
Suppose $\Gamma$ is a $\gamma$-curve.
We list the possible intersections of $\Gamma$ with $\tilde{C}$.
\begin{enumerate}
\item $d=4$ : $\Gamma$ intersects one component of $\tilde{C}$ of multiplicity $\frac{1}{4}$.
\item $d=5$ : $\Gamma$ intersects one component of $\tilde{C}$ of multiplicity $\frac{2}{5}$,
		or two components of multiplicity $\frac{1}{5}$, 
		or two components of multiplicities $\frac{1}{5}$ and $\frac{1}{2}$.
\item $d=6$ : $\Gamma$ intersects one component of $\tilde{C}$ of multiplicity $\frac{1}{2}$.
\end{enumerate}
All intersections are transverse and at a single point, except possibly the $\frac{1}{5},\frac{1}{5}$ case for $d=5$ ---
here the two components of $\tilde{C}$ may coincide, 
and then further the two points of intersection may coincide to yield a point of contact
of order $2$.
\end{cor}

\begin{proof}
This is immediate from the index calculations of Lemma~\ref{lem-index_small_d} together with Lemma~\ref{lem-Gamma3} 
and  Lemma~\ref{lem-Gamma}.
\end{proof}

\begin{thm} \label{thm-d=4,5}
The stable pairs of degrees $4$ and $5$ are as described in Section~\ref{small_d}.
\end{thm}
\begin{proof}
Let $(X,D)$ be a stable pair of degree $d$, where $d=4$ or $5$.
Let $(Y,C)$ be a component of $(X^{\nu},\Delta^{\nu})$.
Then $Y$ is rational since $3 \nd d$.
If $Y \cong \bP^2$ then $C=0$, since  $-(K_Y+C)$ is ample and 3-divisible (again using $3 \nd d$). 
So assume $(Y,C) \not \cong (\bP^2,0)$, and choose $\phi : \tilde{Y} \map \bF_w$ 
as in Notation~\ref{notn-bir_ruling}. We classify the pairs $(\tilde{Y},\tilde{C})$, and hence the pairs $(Y,C)$, 
and finally glue these together to form the surfaces $X$.

Each degenerate fibre of the birational ruling $p: \tilde{Y} \map \bP^1$ consists of some irreducible components of $\tilde{C}$
and some $\gamma$-curves. We have a list of possible connected components of $\tilde{C}$ in Lemma~\ref{lem-C_cpts}, and 
Corollary~\ref{cor-Gamma} describes how $\gamma$-curves and irreducible components of $\tilde{C}$ may intersect. 
Note also that two $\gamma$-curves $\Gamma_1,\Gamma_2$ in a degenerate fibre $f_0$ do not intersect unless 
$f_0=\Gamma_1 \cup \Gamma_2$.
We classify the $(\tilde{Y},\tilde{C})$ by calculating all the ways we can piece together connected components of $\tilde{C}$
and $\gamma$-curves to form a net of curves which is a union of fibres and horizontal curves on some surface $\tilde{Y}$
with a birational ruling $p: \tilde{Y} \map \bP^1$.

So, suppose given some $(\tilde{Y},\tilde{C})$, and assume $\tilde{C} \neq 0$ (otherwise $X=Y$ is smooth, so $X=\bP^2$).
Let $F$ be a connected component of $\tilde{C}$. We consider the cases of Lemma~\ref{lem-C_cpts}. 

We first treat the case $d=4$. In case (4) $F$ cannot intersect a $\gamma$-curve. 
Thus if $F$ intersects a degenerate fibre, it must contain the whole degenerate fibre.
But by inspection $F$ cannot contain a degenerate fibre, so $F$ does not intersect any degenerate fibre.
It follows that there are no degenerate fibres. But we know there are at least two curves of negative self-intersection, 
a contradiction. So case (4) does not occur.
It now follows that case (5) cannot occur. For a $(Y,C)$ with singularities at $C$ corresponding to case (5) 
must be glued to a $(Y',C')$ with singularities at $C'$ corresponding to  case (4) in order to form a smoothable surface $X$.

In case (3), let $E_1,E_2,E_3$ be the components of $F$, ordered as in Lemma~\ref{lem-C_cpts}.
Only $E_3$ can intersect a $\gamma$-curve and $F$ cannot contain a degenerate fibre. Suppose that a component of $F$ is contained 
in a degenerate fibre $f_0$. Then $f_0$ consists of $E_3$ together with some $\gamma$-curves meeting
$E_3$, and possibly $E_1$ or $E_2$. The only possibility is $f_0=E_3 \cup \Gamma_1 \cup \Gamma_2$ where $\Gamma_1,\Gamma_2$ are 
$\gamma$-curves. Then $E_1$ and $E_2$ must be horizontal, but $E_1$ cannot intersect $f_0$, a contradiction.
So every component of $F$ is horizontal. In particular there can be no degenerate fibres (because $E_1$ cannot
intersect a $\gamma$-curve), a contradiction. Thus case (3) does not occur.

In case (1) $F$ cannot intersect a $\gamma$-curve. So there are no degenerate fibres (compare (4) above).
Thus $(\tilde{Y},\tilde{C})=(\bF_4,\frac{1}{2}B)$,  so $(Y,C)=(\bP(1,1,4),0)$.
In case (2) again $F$ cannot intersect a $\gamma$-curve. Thus $\tilde{Y}$ is ruled, $\tilde{Y}=\bF_2$, and so $Y=\bP(1,1,2)$.
Finally $(Y,C)=(\bP(1,1,2),H)$, since $-(K_Y+C)$ is ample and 3-divisible. 
This completes the classification of the $(Y,C)$ for $d=4$. 

Now let $d=5$. In case (6) $F$ cannot intersect a $\gamma$-curve. We obtain a 
contradiction as for $d=4$,(4). So (6) does not occur. It follows that case (7) cannot occur, as for $d=4$,(5). 
For a $(Y,C)$ with singularities at $C$ corresponding to case (7) must be glued to a $(Y',C')$ with singularities at 
$C'$ corresponding to case (6) to form a smoothable surface $X$.

In case (4), let $E_1,E_2,E_3$ be the components of $F$, ordered as above. Only $E_3$ can intersect a $\gamma$-curve, and 
a $\gamma$-curve intersecting $E_3$ does not intersect any other irreducible components of $\tilde{C}$. We proceed as in 
$d=4$,(3) to obtain a contradiction, thus (4) does not occur.

In case (5), $F$ cannot intersect a $\gamma$-curve. Thus $\tilde{Y}$ is ruled, $\tilde{Y}=\bF_5$ and $Y=\bP(1,1,5)$.
The curve $C'$ cannot be horizontal using $-(K_Y+C)$ ample and 3-divisible, thus $C'$ is a fibre and $(Y,C)=(\bP(1,1,5),H)$.

In case (3), write $E_1,E_2,E_3,E_4$ for the components of $F$, ordered as above.
Then only $E_3$ and $E_4$ can intersect a $\gamma$-curve. Now $F$ cannot contain a fibre, thus we see that $F$ 
has components in at most one fibre, and then this fibre contains $E_3$ or $E_4$. 
If every component of $F$ were horizontal, then $E_1$ would intersect a $\gamma$-curve in a degenerate fibre, a contradiction.
Thus $F$ has components in a unique degenerate fibre $f_0$. 

We classify the possible fibres $f_0$ above. We have $E_2 \subset f_0$ because $E_1$ intersects $f_0$, it follows
that $E_3 \subset f_0$. If $E_4 \not \subset f_0$, then $f_0$ consists of $E_2,E_3$, some $\gamma$-curves (meeting
$E_3$) and possibly $E_1$ --- but such a configuration is never a fibre, a contradiction. 
Hence $E_2,E_3,E_4 \subset f_0$. If there is a $\gamma$-curve $\Gamma$ meeting $E_3$,
then $f_0 = E_2 \cup E_3 \cup E_4 \cup \Gamma$. Otherwise, we have some $\gamma$-curves meeting $E_4$.
If $\Gamma_1$, $\Gamma_2$ are two such curves then $f_0=E_4 \cup \Gamma_1 \cup \Gamma_2$, a contradiction.  
Thus there is exactly one $\gamma$-curve, $\Gamma$ say, meeting $E_4$. Then $\Gamma$ intersects another irreducible component $E$ 
of $\tilde{C}$ which has multiplicity $\frac{1}{5}$ or $\frac{1}{2}$ (by Corollary~\ref{cor-Gamma}) and hence self-intersection
$-2$ or $-4$ by Lemma~\ref{lem-C_cpts}. We have $E \subset f_0$. 
If $E^2=-2$ we find $f_0=E_4 \cup \Gamma \cup E$, a contradiction.
Thus $E^2=-4$, we deduce $f_0=E_2 \cup E_3 \cup E_4 \cup \Gamma \cup E$.
Note that $E$ is a \emph{connected} component of $\tilde{C}$ of type (1).  

So, we have two possible types of degenerate fibre $f_0$ --- either $f_0=E_2 \cup E_3 \cup E_4 \cup \Gamma$ where $\Gamma$ meets
$E_3$, or $f_0=E_2 \cup E_3 \cup E_4 \cup \Gamma \cup E$, a chain of curves, where $E^2=-4$. 
In each case $E_1$ is horizontal --- it follows that there are no more degenerate fibres, since $E_1$ cannot intersect a
$\gamma$-curve. We see that $\tilde{Y}$ is obtained from $\bF_7$ by a sequence of blowups, and $E_1=B'$.
Also $C=0$ so $X=Y$.
A graded ring calculation shows that $X=X_{26} \subset \bP(1,2,13,25)$ in the first case.
In the second case we see $X=\bP(1,4,25)$ by toric methods.

Case (8) is very similiar: writing $E_0,E_1,\ldots,E_4$ for the components of $F$, we obtain the same possible degenerate fibres
$f_0$ as above, $E_1$  horizontal, and $E_0$ a fibre. 
$\tilde{Y}$ is obtained from $\bF_2$ by a sequence of blowups, and $E_1=B'$.
A graded ring calculation shows that $Y=X_6 \subset \bP(1,2,3,5)$ in the first case, where $C \sim 2H$.
In the second case we have $Y=\bP(1,4,5)$, where $C \sim 4H$.

Finally, for cases (1) and (2), we may assume that $-(K_Y+C)$ has index $2$, 
because I have already classified the $(Y,C)$ with an 
index $5$ singularity above. Then the same calculation as for $d=4$ gives $(Y,C)=(\bP(1,1,4),0)$ or $(\bP(1,1,2),H)$.
This completes the classification of the $(Y,C)$ for $d=5$.
We calculate that each surface $Y$ has $\rho(Y)=1$, and $(Y,C)$ is a log del Pezzo surface.

We now glue the $(Y,C)$ together to obtain the surfaces $X$. 
The normal surfaces $X$ are the surfaces $Y$ with $C=0$.
The non-normal surfaces $X$ are obtained by glueing components $(Y_1,C_1)$, $(Y_2,C_2)$ along $C_1$,$C_2$ so that
each $(\frac{1}{r}(1,a),\Delta)$ singularity on $(Y_1,C_1)$ is glued 
to a $(\frac{1}{r}(1,-a),\Delta)$ singularity on $(Y_2,C_2)$ 
to give a singularity of type $(xy=0) \subset \frac{1}{r}(1,-1,a)$ on $X$.
We also require that $K_X^2=9$, equivalently $(K_{Y_1}+C_1)^2+(K_{Y_2}+C_2)^2=9$ (in fact this is automatic in 
our cases $d=4$ and $5$).

Finally, we conclude that each surface $X$ constructed as above occurs in a stable pair $(X,D)$ of degree $d=4$ or $5$. 
We need to show that there exists a divisor $D$ on $X$ such that $dK_X+3D \sim 0$, 
$K_X+(\frac{3}{d}+\epsilon)D$ is slc --- this is easy to check.
Moreover, we require that $(X,D)$ admits a $\bQ$-Gorenstein smoothing $(\cX,\sD)/T$ such that 
$\cX_{\bar{\eta}} \cong \bP^2_{\bar{\eta}}$ and $\sD$ is $\bQ$-Cartier --- this follows from Theorem~\ref{thm-F_to_Def^QG_smooth},
using the existence of a $\bQ$-Gorenstein smoothing $\cX/T$ (Theorem~\ref{thm-Manetti} and Theorem~\ref{thm-B}).
\end{proof}

\section{Examples in the case $3 \mid d$} \label{examples}

We first construct explicit smoothings of some slc del Pezzo surfaces of type~C.
We thus obtain examples of stable pairs of degree $6$ --- in particular, we obtain an example where the surface has
 $18$ components.
We next give an example to show that $\cM_d$ is not smooth if $3 \mid d$ (Example~\ref{ex-elliptic_cone_defns}).
Our final example shows that the relative smoothability assumption \ref{defn-smoothable_deformation} 
in the definition of $\cM_d$ is necessary if $3 \mid d$. 
More specifically, we construct an slc del Pezzo surface $X$ which admits $\bQ$-Gorenstein smoothings to 
both $\bP^2$ and a surface which is not smoothable. 

\begin{con} \label{con-typeC_Gor}
Let $C$ be a cycle of smooth rational curves, i.e., $C=C_1 \cup \cdots \cup C_r$, where $C_i \cong \bP^1$ for all $i$, 
$C$ is nodal, and the dual graph of $C$ is a cycle.
Let $T$ be the spectrum of a complete DVR with generic point $\eta$.
Then there exists a smoothing $\cC/T$ of $C$ to an elliptic curve. 
Let $\sL$ be a relatively ample line bundle on $\cC$, write $n_i =\sL \cdot C_i > 0$.
Define $\tilde{\cX}=\underline{\bP}_{\cC}(\cO_{\cC} \oplus \sL^{\vee})$.
The special fibre $\tilde{X}$ is a $\bP^1$-bundle over $C$ with components $\tilde{X_i} \cong \bF_{n_i}$
ruled over $C_i$.
The generic fibre is a ruled surface of degree $\sum n_i$ 
over an elliptic curve. 
Contracting the negative sections, we obtain a family $\cX/T$ with special fibre
$X=\bP(1,1,n_1) \cup \cdots \cup \bP(1,1,n_r)$ (an slc del Pezzo surface of type C) 
and generic fibre an elliptic cone of degree $\sum n_i$.
If we fix $\sum n_i=9$ (equivalently $K_X^2=9$), we deduce that $X$ is smoothable, since an elliptic cone of degree $9$ is 
smoothable by Lemma~\ref{lem-elliptic_cone}.
\end{con}

\begin{con} \label{con-typeC_QG}
We generalise Construction~\ref{con-typeC_Gor}.
Let $C$ and $\cC/T$ be as above. Let $\sL$ be a sheaf on $\cC$ which is invertible in codimension $1$ and $S_2$. 
Locally at a node $P \in C$, we have
$$(C \subset \cC) \cong ((t=0) \subset ((xy=t^r) \subset \bA^3_{x,y,t})) \cong ((uv=0) \subset \bA^2_{u,v} / \frac{1}{r}(1,-1)).$$
Write $\cC' \map \cC$ for the local universal cover $\bA^2_{u,v} \map \bA^2_{u,v} / \frac{1}{r}(1,-1)$.
Then locally $\sL \cong (\cO_{\cC'})_a$, 
where the subscript $a$ is used to denote the eigensubsheaf where a generator $\zeta \in \mu_{r}$ acts as $\zeta^a$.
In particular, $\sL^{[r]}$ is invertible locally. Assume that $\sL$ is relatively ample, write $n_i=\sL \cdot C_i \in \bQ_{>0}$.
Define 
$$\tilde{\cX} = \underline{\Proj}_{\cC}(\oplus_{n \ge 0} S^{[n]}(\cO_{\cC} \oplus \sL^{\vee})),$$ 
where we write 
$S^{[n]}$ for the double dual of the $n$th symmetric power of a sheaf.
Where $\sL$ is invertible, $\tilde{\cX}/\cC$ is a $\bP^1$-bundle and thus easily understood.
Let $P \in C$ be a node where $\sL$ is not invertible. 
Locally at $P \in \cC$, with notation as above, we have
$$\oplus_{n \ge 0} S^{[n]}(\cO_{\cC} \oplus \sL^{\vee}) \cong k[u,v,W,Z]^{\mu_{r}}$$
where $u$, $v$, $W$ and $Z$ have weights $1$, $-1$, $0$ and $a$ with respect to the $\mu_{r}$ action, 
and $W$ and $Z$ have weight $1$ with respect to the grading of the ring. So, locally over $P \in \cC$, we have two affine pieces 
of $\tilde{\cX}$ as follows:
$$(\tilde{X} \subset \tilde{\cX}) \cong ((uv=0) \subset \bA^3_{u,v,w} / \frac{1}{r}(1,-1,-a)),$$
$$(\tilde{X} \subset \tilde{\cX}) \cong ((uv=0) \subset \bA^3_{u,v,z} / \frac{1}{r}(1,-1,a))$$
where $w=\frac{W}{Z}$ and $z=\frac{Z}{W}$. Let $(a,r)=h$, write $a=ha'$, $r=hr'$. Then the surface $\tilde{X}$ has singularities 
$(xy=0) \subset \frac{1}{r'}(1,-1,-a')$ and $(xy=0) \subset \frac{1}{r'}(1,-1,a')$, and the 3-fold $\tilde{\cX}$ has 
a $\frac{1}{h}(1,-1) \times \bA^1$ singularity generically along the double curve. 

We can now give a global description of $\tilde{X} \subset \tilde{\cX}$.
The surface $\tilde{X}$ has components $(\tilde{X_i},\tilde{\Delta_i})$ which are fibred over the components $C_i$ of $C$.
The fibres of $\tilde{X_i}/C_i$ are smooth rational curves (when given their reduced structure).
The double curve $\tilde{\Delta}_i$ is the sum of the (reduced) fibres of $\tilde{X_i}/C_i$ over the nodes of $C$ on $C_i$. 
Let $P$ be a node of $C$ on $C_i$ where $\sL$ is not invertible.
With notation as above, the fibre of $\tilde{X_i}/C_i$ over $P$  has multiplicity $r'$, and $(\tilde{X_i},\tilde{\Delta_i})$
has singularities $(\frac{1}{r'}(1,a'),\Delta)$ and $(\frac{1}{r'}(1,-a'),\Delta)$ at the fibre.
The bundle $\tilde{X_i}/C_i$ has a negative section of self-intersection $-n_i \in \bQ$.
We obtain $\tilde{X}$ by glueing the components $\tilde{X_i}$ along the fibres over the nodes of $C$.
The singularities of type $(\frac{1}{r}(1,a),\Delta)$ on the components glue to give 
singularities of type $(xy=0) \subset \frac{1}{r}(1,-1,a)$ on $\tilde{X}$, elsewhere $\tilde{X}$ is normal crossing.
The negative sections of $\tilde{X_i}/C_i$ glue to give a section of $\tilde{X}/C$.
The family $\tilde{\cX}/T$ is a $\bQ$-Gorenstein smoothing of $\tilde{X}$, with generic fibre a ruled surface
over an elliptic curve of degree $\sum n_i$.

We now contract the negative sections. 
We obtain a $\bQ$-Gorenstein family $\cX/T$ with special fibre $X$ an slc del Pezzo surface 
of type C and generic fibre an elliptic cone of degree $\sum n_i$. 
Again, assuming $K_X^2=\sum n_i =9$, we deduce
that $X$ admits a $\bQ$-Gorenstein smoothing to $\bP^2$.
\end{con}

\begin{prop} \label{prop-smoothability_typeC}
Let $X$ be a surface of type C. Assume that
\begin{enumerate}
\item $X$ has slt singularities of types $(xy=0) \subset \frac{1}{r}(1,-1,a)$, $(a,r)=1$.
\item $K_X^2=9$.
\item For every component $(Y,C)$ of $X$, $(Y,C)$ is toric, i.e., 
$Y$ is a toric surface and $C$ is a sum of toric strata of codimension $1$, and $\rho(Y)=1$.
\end{enumerate}  
Then $X$ admits a $\bQ$-Gorenstein smoothing to $\bP^2$.
\end{prop}

\begin{rem}
Note that conditions (1) and (2) are necessary.
\end{rem}

\begin{proof}
Let $(Y,C)$ be a component of $X$.   
Then there is a unique toric blowup 
$\tilde{Y} \map Y$ which blows up the node of $C$ such that the strict transforms of the components of $C$ give fibres of 
a fibration $\tilde{Y} \map \bP^1$. For, if $Y$ is given by the fan defined by the vectors $v_0,v_1,v_2 \in N_{\bR} \cong \bR^2$, 
where $v_1$ and $v_2$ correspond to the components of $C$, 
the map $\tilde{Y} \map Y$ corresponds to the subdivision of the fan obtained by 
adding the vector $-v_0$. Write $\tilde{C}$ for the strict transform of $C$.
Glueing the $(\tilde{Y},\tilde{C})$ together 
so that the negative sections form a cycle, we obtain a partial resolution $\tilde{X} \map X$. 
We claim that, for a suitable choice of $\cC/T$ and $\sL$ in Construction~\ref{con-typeC_QG}, we can 
recover $\tilde{X}$ as the special fibre of $\tilde{\cX}=\underline{\Proj}_{\cC}( S^{[n]}(\cO_{\cC}\oplus \sL^{\vee}))$.
Hence we obtain a smoothing of $X$ as above.

To prove the claim, note that $\tilde{X}$ is uniquely determined by the following data --- 
the self-intersections $-n_i$ of the negative sections of the components $\tilde{X_i}$, 
and the $(xy=0) \subset \frac{1}{r}(1,-1,a)$ singularities at the double curves. Let $\cC/T$ be any smoothing of $C$.
Possibly after base change, we can choose $\sL$ locally at the nodes of $C$ to obtain the required singularities 
(compare Construction~\ref{con-typeC_QG}).
These local choices extend to a global sheaf $\sL$. We require that $\sL \cdot C_i =n_i$ in order to obtain the correct 
self-intersections of the negative sections --- we can achieve this by twisting $\sL$ by a suitable line bundle.  
\end{proof}

\begin{ex} \label{ex-typeC}
Every candidate surface $X$  for $d=6$ (see \ref{sss-typeC&D}) of type C with components of Picard number $1$  
does indeed occur in a stable pair $(X,D)$ of degree $6$. 
For, by Proposition~\ref{prop-smoothability_typeC}, $X$ admits a $\bQ$-Gorenstein smoothing to $\bP^2$.
It only remains to show that there exists a divisor $D$ on $X$ such that $(X,D)$ is a stable pair of degree $6$.
Equivalently, we require $D \in |-2K_X|$ such that $K_X+(\frac{1}{2}+\epsilon)D$ is slc for some $\epsilon > 0$.
Let $\pi \colon \tilde{X} \map X$ be the partial resolution as in the Proof of Proposition~\ref{prop-smoothability_typeC}, 
and let $p \colon \tilde{X} \map C$
be the fibration. Let $\tilde{D}$ be a general double section of $p$ which is disjoint from the negative section.
Then $D=\pi_{\star} \tilde{D}$ has the required properties.
Note that smoothability of the \emph{pair} $(X,D)$ follows by Theorem~\ref{thm-F_to_Def^QG_smooth}.

In particular, we obtain an example of a stable pair $(X,D)$ of degree 6 such that the surface $X$ has $18$ components.
We construct $X$ by glueing 18 components $(Y,C) \cong (\bP(1,1,2), H+\sim 2H)$ together to form a surface of type C. 
At the degenerate cusp all the components of $X$ are smooth, and 
$X$ has 9 slt singularities of type $(xy=0) \subset \frac{1}{2}(1,1,1)$.
We note that 
$18$ is the maximum number of components of a surface $X$ occurring in a stable pair of degree $6$. For, by 
Lemma~\ref{lem-index_small_d} we have $2K_X$ Cartier, thus for each component $(Y,C)$ of $X$ we have 
$(K_Y+C)^2 \in \frac{1}{2}\bZ$, in particular $(K_Y+C)^2 \ge \frac{1}{2}$. 
Since $\sum (K_Y+C)^2 = K_X^2 = 9$, there are at most $18$ components as required.
\end{ex}

\begin{ex}\label{ex-elliptic_cone_defns}
The deformation space $\Def X$ of an elliptic cone $X$ of degree $9$ consists of $9$ smooth $1$-dimensional components, 
together with an embedded component at the origin. The components of $(\Def X)_{\red}$ meet as transversely as possible,
i.e., $(\Def X)_{\red}$ is isomorphic to the collection of coordinate axes in $(0 \in\bC^9)$.
Each component corresponds to a smoothing of $X$ to $\bP^2$ as constructed in Lemma~\ref{lem-elliptic_cone}.
Write $E$ for the section of $X$. 
There is an action of the $9$-torsion of $\Pic^0 E$  on $(\Def X)_{\red}$
which permutes the components transitively,
the $3$-torsion elements induce automorphisms of each component. For details see \cite{Ste}, p. 220, Example~4.5.

Now, we can add a divisor $D$ on $X$ such that $(X,D)$ is a stable pair of degree $d$, for any $d$ such that $3 \mid d$.
We describe the singularities of $\cM_d$ at the point $[(X,D)]$. 
By Theorem~\ref{thm-F_to_Def^QG_smooth}, $\cM_d$ is smooth over the closed subscheme of $\Def X$ corresponding to 
$\bQ$-Gorenstein smoothable deformations. First, since $K_X$ is Cartier, every deformation is $\bQ$-Gorenstein.
Second, the smoothable deformations correspond (by definition) to the scheme
theoretic closure in $\Def X$ of the locus where the geometric fibres are isomorphic to  $\bP^2$.
In our case this is exactly $(\Def X)_{\red}$ by the explicit description above.
Thus $\cM_d$ is smooth over $(\Def X)_{\red}$. In particular, $\cM_d$ is \emph{not} smooth at $[(X,D)]$.
\end{ex}

\begin{ex}\label{ex-rel_smoothability}
We define an slc del Pezzo surface $X$ of type C as follows: The normalisation $(X^{\nu},\Delta^{\nu})$ has two components
$(X_1,\Delta_1) \cong (\bF_4, A+B)$ and $(X_2,\Delta_2) \cong (\bF_1,A+B)$, where we write $A$ for the fibre and $B$ for the 
negative section of a rational ruled surface $\bF_n$. $X_1$ is glued to $X_2$ by identifying the negative section of $\bF_4$
with a fibre of $\bF_1$ and vice versa. Thus $X$ has a degenerate cusp, we remark that locally analytically at this point 
we have 
$$X \cong ((xy-z)z=0) \subset \bA^3_{x,y,z}.$$

We claim that $X$ smoothes to $\bP^2$. We construct an explicit smoothing as follows: Let $T$ be the spectrum of a DVR, let 
$\cY=\bP^2_T$. Let $\cY^1/T$ denote the blowup of $\cY/T$ with centre a smooth conic $Q$ in the special fibre $Y=\bP^2$. 
Then the special fibre 
$Y^1$ consists of the strict transform $Y' \cong \bP^2$ of $Y$ together with the exceptional divisor $E \cong \bF_4$, 
glued along the 
conic $Q$ and the negative section $B$. Let $\cY^2/T$ denote the blowup of $\cY^1/T$ with centre a fibre $A$ of $E$.
Then $Y^2=Y''+E'+F$, where $Y'' \cong \Bl \bP^2 \cong \bF_1$, $E' \cong \bF_4$, and $F \cong \bF_1$ is the exceptional divisor, 
glued to the fibre $A$ of $E'$ along its negative section, and glued to the $-1$-curve of $Y''$ along a fibre.    
We calculate that $-K_{\cY^2}$ is ample on $E'$ and $F$, and defines the fibering contraction $\bF_1 \map \bP^1$ of
$Y'' \cong \bF_1$. Thus some multiple of $-K_{\cY^2}$ defines a divisorial contraction $\cY^2 \map \bar{\cY}/T$
which contracts $Y'' \cong \bF_1$ to a curve. We observe that the special fibre $\bar{Y}$ is our surface $X$, so $\bar{\cY}/T$
is a smoothing of $X$ as required.

We next construct a partial smoothing of $X$ to a rigid normal crossing surface. We consider the components of $X$ separately.
First, the component $(X_1,\Delta_1) \cong (\bF_4,A+B)$ has a deformation $(\cX_1,\cC_1)/T$ with generic fibre 
$(\bF_2,B)$. The family $\cX_1/T$ can be obtained as a deformation of scrolls in some $\bP^N$.
Second, the component $(X_2,\Delta_2) \cong (\bF_1,A+B)$ has a deformation $(\cX_2,\cC_2)/T$ with generic fibre
$(\bF_1,C)$, where $C \sim A+B$ is a section of $\bF_1$ disjoint from $B$. Here $\cX_2/T$ is the trivial deformation.
Then, possibly after base change, we can glue $\cX_1/T$ and $\cX_2/T$ along $\cC_1/T$ and $\cC_2/T$ to obtain a partial smoothing
$\cX/T$ of $X$. The generic fibre is the normal crossing surface $Z$ obtained by glueing $(\bF_2,B)$ and $(\bF_1,C)$.
The surface $Z$ is rigid, in particular it is not smoothable. To see this, note that, by Lemma~\ref{lem-T^1}, we have
$\cT^1_Z \cong \cO_{\bP^1}(-1)$. So $H^0(\cT^1_Z)=0$, thus all deformations of $Z$ are locally trivial. It is easy to see that
$Z$ has no locally trivial deformations (by considering the components separately), hence $Z$ is rigid as claimed.

We can add a divisor $D$ such that $(X,D)$ is a stable pair of degree $d$, for any $d$ such that $3 \mid d$. Then, 
working locally at $[(X,D)] \in \cM_d$, the stack $\cM_d$ is smooth over the closed subscheme $\Def^{sm} X$ 
of $\Def X$ corresponding to smoothable deformations. We want to emphasise that, set-theoretically,
$\Def^{sm} X$ is strictly smaller than $\Def X$. For, by our construction above, there exists an irreducible component 
of $\Def X$
whose generic point corresponds to a surface $Z$ that is not smoothable. Thus this component is not contained in $\Def^{sm} X$.
In particular, we see that the relative smoothability condition (\ref{defn-smoothable_deformation}) for $\cM_d$ is necessary. 
For, let $(\cX,\sD)/T$ be a deformation of $(X,D)$ over the spectrum of a complete DVR such that $\cX/T$ has generic fibre $Z$. 
Write $\cM'_d$ for the groupoid obtained by dropping the relative smoothability assumption in the definition of $\cM_d$. 
Then $(\cX,\sD)/T \not \in \cM'_d(T)$ (since the generic fibre is not smoothable), 
but $(\cX,\sD) \times_T A \in \cM'_d(A)$ for any Artinian thickening $A$ of $0 \in T$. 
It follows that there does not exist a versal deformation $(\cX^u,\sD^u)/U \in \cM_d'(U)$ of $(X,D)$, 
thus $\cM'_d$ is not an algebraic stack.   
\end{ex}

\end{document}